\input amstex
 \documentstyle{amsppt}
 \magnification=\magstep1
 \vsize=24.2true cm
 \hsize=15.3true cm
 \nopagenumbers\topskip=1truecm
 \headline={\tenrm\hfil\folio\hfil}

 \TagsOnRight

\hyphenation{auto-mor-phism auto-mor-phisms co-homo-log-i-cal co-homo-logy
co-homo-logous dual-izing pre-dual-izing geo-metric geo-metries geo-metry
half-space homeo-mor-phic homeo-mor-phism homeo-mor-phisms homo-log-i-cal
homo-logy homo-logous homo-mor-phism homo-mor-phisms hyper-plane hyper-planes
hyper-sur-face hyper-sur-faces idem-potent iso-mor-phism iso-mor-phisms
multi-plic-a-tion nil-potent poly-nomial priori rami-fication sin-gu-lar-ities
sub-vari-eties sub-vari-ety trans-form-a-tion trans-form-a-tions Castel-nuovo
Enri-ques Lo-ba-chev-sky Theo-rem Za-ni-chelli in-vo-lu-tion Na-ra-sim-han Bohr-Som-mer-feld}

\define\rest#1{_{\textstyle{|}#1}} 

\define\Span#1{\left<#1\right>} 

\define\half{{\textstyle{1\over2}}}


\define\C{\Bbb C} 
\define\R{\Bbb R} 
\define\Z{\Bbb Z} 

\define\proj{\Bbb P} 

\define\sA{{\Cal A}} 

\define\al{\alpha}
\define\be{\beta}
\define\de{\delta}

\define\om{\omega}
\define\si{\sigma}
\define\De{\Delta}
\define\Ga{\Gamma}
\define\La{\Lambda}
\define\Om{\Omega}
\define\la{\lambda}









 \document

 \topmatter
  \title Geometric quantization in the framework
of algebraic Lagrangian geometry \endtitle
  \author    Nikolai Tyurin (Tiourine)
                \endauthor

   \address MPI, Bonn
   \endaddress
  \email ntyurin\@thsun1.jinr.ru
      jtyurin\@mpim-bonn.mpg.de
   \endemail

\abstract  This is a short version of the author's habilitation thesis.
The main results have been published but many details  are developed and clarified. As well
some new results are included: we additionally discuss here quasi classical limit of 
ALG(a) - quantization, mention some topological properties of the moduli space
of Bohr - Sommerfeld Lagrangian cycles of fixed volume and investigate some properties
of the Kaehler structure on it.
\endabstract

   \endtopmatter

\head Introduction  \endhead

The main theme of this paper is 
 quantization of the classical mechanical systems
in terms of  algebraic geometry. Thus here one relates
the questions of theoretical physics and mathematics.
To start with let us recall briefly the main problems and methods which turn us to study a new subject.

Quantization itself is the main topic of the theoretical
physics. Necessity of  its introduction and development
was dictated by the creators of the quantum theory.
According to the Copenhagen philosophy, the physical
predictions of a quantum theory must be formulated in terms
of classical concepts (the first phrase  of [33]; here
we quote the beginning of this survey). So in addition to the usual structures
(Hilbert space, unitary transformations, selfadjoint operators...) any sensitive quantum theory has to admit
an appropriate passage to a classical limit such that
the quantum observables are transfered to the classical ones.  
However as it was pointed out by Dirac at the beginning
of the quantum age  the correspondence between quantum theory
and classical theory has to be based not only on 
numerical coincidences taking place in the limit
$h \to \infty$ but on an analogy between their mathematical
structures. Classical theory does   approximate
 the quantum theory but it does do even more --- it supplies
a frame to some  interpretation of the quantum theory.
Using this idea we can understand quantization procedure
in general as a correspondence between classical theories
and quantum theories. In this sense quantization
of the classical mechanical systems is the moving in
 one direction while taking quasi classical limit we go
in the opposite direction of this correspondence. More
abstractly: the moduli space of the quantum theories
is a $n$ - covering of the moduli space of the classical ones
(one supposes that $n$ equals 2), and quantization
is the structure of this covering.

Quantization at all is a very popular subject. There are a number of different approaches to this problem. But one of them is honored as the first one in theoretical physics
and it is named as canonical quantization. In simple cases
the correspondence comes with some choice of fixed
coordinates. If  classical observable is represented
by a function $f(p_a, q^b)$ in these coordinates
then the corresponding quantum observable equals to the operator
$$
f(- \imath h \frac{\partial}{\partial q^a}, q^a).
$$
The canonical quantization of the harmonic oscillator    
is a standard of the theoretical physics: any alternative
approach should be compared with it and if the answer
is sufficiently different from the classical one then this approach is rejected. However this formal substitution
(when one puts some differential operators instead
of coordinates $p_a$) introduces a lot of problems.
Indeed, beyond of the simple cases during this process
the result of the quantization depends on the order
of $p$ and $q$ in the expression for the classical observable
$f$ and moreover the result strongly depends on the coordinate choice and it is not invariant under generic
canonical transformations. Nevertheless this canonical quantization supplied  by some physical intuition
together with its various generalization takes the central part of the modern theoretical physics.

One way to develop the canonical method and avoid the difficulty is provided by geometric
quantization. The geometric quantization (which is discussed in the text) has two
slightly different meanings as a term. One could understand this one
either as a concrete construction (see f.e. [13], [14], [21], [33], etc) well known
as Souriau - Kostant quantization or as a general approach to the problem
based on geometry. Nowadays the problem of quantization is used to be solved
by quite different methods: algebraic approach includes deformation quantization,
formal geometry, non commutative geometry, quantum groups; analytical approach
consists of  the theory of integral Fourier operators, Toeplitz structures  and
other ones. All the methods mentioned above have one mutual marking point ---
using these we almost completely forget about the  structure of the given system
(and the Dirac suggestion mentioned above) and the "homecoming" turns to be absolutely impossible. At the same time
going in the geometric quantization direction one at least tries to
keep (at least  in mind) the original system. The corresponding symplectic
manifold remains to be basic for all the constructions and takes real part in the definition
of all auxiliary geometrical objects which give us the result of the quantization.
At the same time the starting point of the geometric quantization is to avoid the choice of
any coordinates and this basic wish gives a possibility to deal with the complicated systems which do not admit any global coordinates at all. 
 But starting with a given classical phase space geometric quantization
should give us a result which has to be comparable
with the canonical  one for the simple systems. Thus in any case geometric quantization
is a generalization of the canonical quantization. To keep the relationship
one usually pays a cost loosing generality of the construction: from the whole space
of classical observables one takes only a subclass of "quantizable" functions,
and this subclass is sufficiently small. To separate such quantizable objects one should
choose a polarization of given symplectic manifold (see [20], [33]) then these objects are
distinguished by the condition that their Hamiltonian vector fields preserve the
polarization. The separation necessity is dictated by a classical von Hove result
(see f.e. [13]) which ensures us that it is impossible to realize so called
"ideal quantization" over given manifold in terms of usual geometry of this manifold
(functions, differential operators etc). Let us recall that one understands such ideal quantization
of classical mechanical system corresponding to given symplectic manifold $(M, \om)$ as a procedure 
giving us some appropriate Hilbert space $\Cal H$ together with an irreducible representation
$$
q: C^{\infty} (M, \R) \to O(\Cal H)
$$
(where the first space one takes as a Lie algebra endowed with the Poisson bracket and the second
one is the space of self adjoint operators on $\Cal H$) such that it satisfies a correspondence
  principle, namely
$$
i \hbar [q(f), q(g)] =  q(\{f, g\})
$$
for every $f, g$ (and, of course, one is interested in the cases when 
the operators on the right hand side satisfy a number of properties which makes
the corresponding quantum mechanical system meaningful and "computable") so $q$ has to be a representation
of the Poisson algebra. The requirement for $q$ to be irreducible led to the separation
mentioned above of subclasses of "quantizable functions" in the "classical" construction
which is known as Souriau - Kostant geometric quantization. As the result of this construction one
gets some representation of whole space $C^{\infty}(M, \R)$ which is manifestly reducible
being the direct sum of two isomorphic representations. This imposes consideration 
and quantization of a half part of the classical observable space. 

The known schemes of geometric quantization are unified by the fact that usually they takes
the space of regular sections of a prequantization bundle as the Hilbert space (and again
one imposes some additional conditions on these sections to be regular in our sense). In
original Souriau - Kostant construction one takes all smooth sections with
bounded  $L^2$ - norm  (with respect to a given hermitian structure on the fibers of
the prequantization bundle weighted by the Liouville form). Further specializations come
in  different ways:   Rawnsley - Berezin method (see [16]) uses only the sections
which are holomorphic with respect to a complex polarization (= fixed complex structure
on $M$) as well as in  Toeplitz - Berezin  approach (see [4]) while in the real polarization
case one collects only such sections (weighted by half weights) which are invariant with respect to
infinitesimal transformations tangent to the fibers of a real polarization.

One should say that the introduction of an additional structure ---
complex polarization --- turns the subject of geometric quantization
to  the most developed region of modern mathematics namely to algebraic geometry.
As it was mentioned above a number of methods uses  complex polarization. It imposes
the additional condition that our symplectic manifold $(M, \om)$ admits a Kaehler structure: there exists  some complex structure $J$ compatible with $\om$ which is integrable.
Together these two structures $\om, J$ give us the corresponding riemannian
metric $g$ such that complex manifold $M, J$ is endowed with
a hermitian metric, and since $\om$ is closed by the definition
it gives us a Kaehler structure on $M$. Moreover, one has usual for any
quantization method requirement for $\om$ to have integer cohomology class:
$$
[\om] \in H^2(M, \Z) \subset H^2(M, \R)
$$ 
(the charge integrality condition). This implies that the Kaehler metric described above is of the Hodge type and therefore
the Kaehler manifold is an algebraic variety. So one can quantize
a symplectic manifold if it admits the algebraic geometrical structure!
It is not so surprising if we take in mind  so called geometric formulation
of quantum mechanics. The basic idea is to replace the algebraic methods
of quantum mechanics by algebro geometrical methods. The author found
all these ones in [3], [18] but of course the original sources exist
as one thinks since the birth of the quantum theory itself.
Anyway [18] contains  the history of the question. Roughly speaking,
the starting point is that usually in quantum mechanics one
deals with a Hilbert space but the quantum states are represented
 by rays in the space since two vectors $\psi_1, \psi_2$
represent the same state iff they are proportional. 
Thus it is natural to consider the projectivization
$\proj (\Cal H)$ instead of $\Cal H$ as the space of quantum states.
This finite - or infinite - dimensional complex manifold is automatically
endowed with a hermitian metric (Fubini - Study) so one can regard 
it as a real manifold with Kaehler structure. This real manifold (finite
or infinite dimensional) is endowed automatically with  symplectic
structure and riemannian metric. Quantum states are represented just by points
of this manifold. Quantum observables are represented by smooth real functions
of special type which one should call Berezin symbols. The specialty
means that the desired functions are distinguished by the following condition:

{\it function $f$ is a Berezin symbol iff its Hamiltonian vector field
preserves the riemannian metric as well as the symplectic form.}

Further, instead of the commutator of two operators one has 
just the usual Poisson bracket for the corresponding symbols. Eigenstate means
critical point for the symbols, eigenvalue equals to the critical value.
Dynamics of the quantum system is described in terms of the classical Hamiltonian
dynamics. Probability amplitudes are given by the geodesic distances so the riemannian metric on the space responds to the probability aspects of the theory([3], [18]).

Applying these ideas to geometric quantization let us reformulate and slightly
generalize the main subject. First of all we give the following
\proclaim{Definition}
Let $\Cal K$ is a Kaehler manifold endowed with a Kaehler triple $(\Om, I, G)$.
Then real function $f$ is called quasi symbol iff its Hamiltonian vector field
preserve the riemannian metric:
$$
Lie_{X_f} G = 0.
$$
\endproclaim  
  
For any Kaehler manifold $\Cal K$ the set of quasi symbols is a Lie subalgebra
in the Poisson algebra. Let us denote it as $C^{\infty}_q(\Cal K, \R)$. If $\Cal K$
is a projective space then the notion is equivalent to the notion of
Berezin symbol.

Now we are ready to formulate what we will call
\proclaim{Algebro geometric quantization (AGQ)} Let $(M, \om)$ be a symplectic manifold
corresponding to a classical mechanical system. Then one says that AGQ results
with an appropriate algebraic variety $\Cal K$ together with a correspondence
$$
q: C^{\infty}(M, \R) \to C^{\infty}_q(\Cal K, \R),
$$
satisfying the following conditions:

1) $q(f+g) = q(f) + q(g)$;

2) $q(cf) = c q(f)$ where $c$ is a real constant;

3) $q(c) = c$ where $c$ is a real constant;

4) $q$ is irreducible;

5) the correspondence principle in the form
$$
\{q(f), q(g) \}_{\Om} = q(\{f, q\}_{\om})
$$
holds for every functions $f, g$.
\endproclaim

In this list we need to explain what the irreducibility (item 4) means. One knows
what does it mean for usual linear representations. But since we start with the case
of the projective spaces one sees that during the projectivization the usual irreducibility
condition is translated to the following double:

4a) $q$ has trivial kernel;

4b) the distribution $Vect_{\Om}(Im q)$ over $\Cal K$ either is non integrable or
spans whole the tangent bundle of $\Cal K$,

where $Vect_{\Om}(Im q) \subset Vect (\Cal K)$ is the Lie subalgebra which consists
of the Hamiltonian vector fields of the functions from the image of $q$.    
We will show in Section 1 that these two conditions form an analogy of
the usual "linear" irreducibility condition.

Thus the statement of the quantization problem is slightly changed.
For a given symplectic manifold we are looking for an appropriate
algebraic manifold (variety) together with an inclusion $q$. Following
the authors of [3] we require for the construction of this algebraic
manifold to avoid as an intermideate step the introduction of  Hilbert spaces
known from  usual methods of geometric quantization. 

The main aim of this text is to present an example of succesfull
algebro geometric quantization for compact simply connected symplectic
manifolds. We call this method ALG(a) - quantization. Of course,
it is an abbreviation. To decode this one we need to recall
some basic facts belonging to a new subject which was invented 
just on the border between algebraic and symplectic geometries
(if this border does exist). 

One could say that different subjects are mixed in modern mathematics.
For example, in connection with the mirror symmetry conjecture
one accepts the idea that  algebraic geometry of a manifold $X$
corresponds to   symplectic geometry of its mirror partner $X'$.
The ingredients of algberiac geometry over $X$  (bundles, sheaves,
divisors ...) are compared with some derivations of symplectic geometry
(Lagrangian submanifolds of special types). F.e. in so called homological mirror
symmetry one compares two categories came from algebraic geometry and symplectic
geometry respectively and in some particular cases (elliptic curve)
this approach gives the desired result. On the other hand one has a number of moduli spaces generated in the framework of algebraic geometry over $X$ and
another way is to find a number of moduli spaces in the framework 
of symplectic geometry. The development of this idea comes in different ways
and even now one could report about a number of promised results
and ideas clarifying the original one (see, f.e., [22], [24]). But these results are sufficiently far to be complete and to cover all the problems. But the main idea which proclaims
 to create some new synthetic (or at least synergetic) geometry  unifying
 algebraic and 
symplectic ones  remains to be very attractive and seems to be true.

One step in this way was done in 1999 when the moduli space of half weighted
Bohr - Sommerfeld Lagrangian cycles of fixed volume and topological type was proposed in [25] and constructed in [10]. Starting with a simply connected compact symplectic manifold with integer symplectic form (read "classical mechanical system with compact simply connected phase space which satisfies the Dirac condition") the authors 
construct a set of infinite dimensional moduli spaces which are infinite dimensional algebraic manifolds in dependence on   the choice of some topological fixing and a real number --- the volume of the half weighted cycles. Lagrangian geometry is mixed    in the construction with algebraic geometry and this construction itself belongs
to some new synthetic geometry. The authors called it ALAG --- abelian Lagrangian algebraic geometry (so it is wrong to think that they took their initials
and made a mistake). It was created as a step in the approach to mirror symmetry conjecture
generalizing some notions from  standard geometric quantization (prequantization data, Bohr - Sommerfeld
condition, etc.) so it is not quite surprising that this construction should play an important role
in geometric quantization. Namely, it was shown in [28] and [30] that these moduli spaces of
half weighted Bohr - Sommerfeld Lagrangian subcycles of fixed volume solve the problem
of algebro geometric quantization stated above for simply connected compact symplectic
manifolds. Briefly, if $(M, \om)$ is the symplectic manifold with integer symplectic
form then one can construct the moduli space $\Cal B^{hw,r}_S$ of fixed volume $r$
consists of pairs $(S, \theta)$ where $S$ is a Bohr - Sommerfeld Lagrangian cycle and $\theta$ is a half weight over it such that the volume $\int_S \theta^2$ equals to $r$ (see Section 2).
For every smooth function $f \in C^{\infty}(M, \R)$ (= classical observable) one has a smooth function
$F_f$ on the moduli space which is given by the map
$$
\aligned
\Cal F_{\tau}: C^{\infty}(M, \R) \to C^{\infty}(\Cal B_S^{hw,r}, \R), \\
\Cal F_{\tau}(f) (S, \theta) = F_f(S, \theta) = \tau \cdot \int_S f|_S \theta^2.\\
\endaligned
$$
It's easy to see that $\Cal F_{\tau}$ maps constants to constants:
$$
f \equiv s \implies F_f \equiv \tau \cdot r \cdot c.
$$
Moreover, we prove that this $\Cal F_{\tau}$ is a homomorphism of Lie algebras namely
$$
\{ F_f, F_g \}_{\Om} = 2 \tau F_{\{f, g\}_{\om}}
$$
for any $f, g \in C^{\infty}(M, \R)$ (we prove it in Section 2 by direct computations).
Moreover, in Section 3 we prove that for every $f \in C^{\infty}(M, \R)$ the corresponding
function $F_f$ is a quasi symbol (= quantum observable). Moreover, we prove that $\Cal F_{\tau}$ is
irreducible.

Thus, summing up all the facts one ensures that the pair $(\Cal B^{hw, r}_S, \Cal F_{\tau})$
satisfies all conditions of algebro geometric quantization if one takes
$$
\aligned
\tau = \half, \\
r = 2. \\
\endaligned
$$
So we have that
$$
\aligned
\Cal K = \Cal B^{hw, 2}_S, \\
q = \Cal F_{\half}
\\
\endaligned
$$
is a solution for AGQ problem.

This method, proved in [28], [30], was called ALG(a) - quantization. This new method gives new results
which are nevertheless quite consistent with the old ones if an appropriate polarization
on $(M, \om)$ is chosen ( see [30]). In Section 4 we show that if $(M, \om)$ admits
a real polarization then one can derive the result of the standard geometric  quantization described f.e.
in [20] from the result of ALG(a) - quantization. The same is true for the complex case. 
It means that ALG(a) - quantization is a natural generalization of the known methods
of geometric quantization.

At the same time we can change our parameters $r$ and $\tau$. Usually in geometric quantization
one has an additional integer parameter: if $(L, a)$ is the prequantization line bundle with prequantization $U(1)$ - connection deriving the Bohr - Sommerfeld subcycles from the space of all Lagrangian subcycles
then one can take its tensor power $(L^k, a_k)$ in dependence with  integer parameter $k$.
Then one understands Planck constant in this construction, following F. Berezin, as 
$$
\hbar = \frac{1}{k},
$$
and it is not hard to see that while $h$ tends to zero one gets a dense subset of the space of all Lagrangian
subcycles which consists of Bohr - Sommerfeld subcycles of any levels. The point is that
dealing with the parameters one gets the following limit: any Lagrangian subcycle 
plus zero half weight. Hence the  space of the quasi classical states is the space of all unweighted (or unhalfweighted) Lagrangian subcycles of
fixed topological type. We  show that this space is endowed with say canonical dynamical properties
while it hasn't any   canonical riemannian metric to define the corresponding probability
aspects. Really if we choose a Hamiltonian function $H$ on the given symplectic manifold then we
automatically get an induced  vector field on the space of all Lagrangian subcycles. It means
that we have dynamics on the space of Lagrangian subcycles which is induced by the classical dynamics
 of the given system and is compatible with it. At this point the story comes beyond of usual geometry.
Namely the thing which is generated by a smooth function on the given manifold
is not just a vector field on the space of Lagrangian submanifolds.  If $f$ is constant
on a Lagrangian submanifold $S$ then it defines a number which is evedentely equals to
this constant. This means that every $f$ generates a "thing" which is neither  vector field no 
function but it belongs to another geometry. Usually it is called super geometry. 
The "thing" should be called super function but of some special type.
If $\Cal L_S$ is the space of all Lagrangian subcycles of fixed topological type then
one takes odd supersymplectic manifold $\Pi T^* \Cal L_S$. Every vector field over $\Cal L_S$
now can be regarded as an odd super function on $\Pi T^* \Cal L_S$ so 
it has parity 1 while any numerical function has parity 0 as a super function.  "Super selection rules" usually
forbid to consider nonhomogineous super functions with different parity in different points but nevertheless
we've got such special nonhomogineous super functions. So the "thing" induced by any real function
$f$ is a vector field in  Lagrangian subcycles if it's not constant at these cycles and
it is a number   function on the cycles which belong to the level set of this function.
Of course, these super functions have very special type but nevertheless the fact takes place
and the further investigation of the subject should keep this way as well as the other possible ways.
While the mentioned observation belongs to odd supersymplectic geometry we have another point 
enforcing one to consider  even super symplectic case. One should notice that the definition of Bohr - Sommerfeld Lagrangian cycles has very natural and simple reformulation on the language of even super symplectic manifolds.
For this let us translate to the language the usual setup for geometric quantization.
One has a real symplectic manifold $(M, \om)$ and a prequantization bundle $(L, a)$
together with the corresponding prequantization $U(1)$ - connection.
Then one could consider this picture (see [17]) as follows: we have an even super symplectic manifold.
At each point of the total space 
$$
E = L \to M
$$
one has a decomposition of the tangent bundle $TE$ with respect to  connection $a$ and
the super symplectic form equals to usual $\om$ being restricted to  the horizontal part
and to the real part of the hermitian metric at the vertical part. The fact that $a$ is hermitian
implies that the Jacobi identity holds for the induced super Poisson bracket (Butten bracket).
 After this setup is understood one can take Lagrangian submanifolds of the even super symplectic
manifold $E$ having in mind the usual sense: $S$ is Lagrangian if it is isotropic with respect to
the super symplectic form and has maximal dimension. The point is that such Lagrangian submanifolds
project exactly to our Bohr - Sommerfeld submanifolds. And the quantization condition belongs
exactly to the super geometry. It gives a hint how the abelian case could be extended
to non abelian cases. 

 The text is organazied as follows. 

The first section briefly recalls and discusses two
sources of the new method: the geometric formulation of quantum mechanics and geometric
quantization itself. We start with the first one, almost copying appropriate parts of
[3] so we would like to strongly recommend this paper as  complete and perfect. 
Further we touch known methods of geometric quantization; at the same time 
we add some realizations  of the methods in terms of the geometric formulation
of quantum mechanics. 

In Section 2 we introduce a natural map from the space of smooth functions on a given
simply connected compact symplectic manifold to the space of smooth functions
over the corresponding moduli space of half weighted Bohr - Sommerfeld
Lagrangian cycles and prove that this map is a homomorphism of the Poisson algebras.
To do this first of all we recall the construction of the moduli spaces
(and again we recommend to use the original source [10] or its preliminary  version
[9]), describing the local geometry of the moduli space needed to the further computation.
Then we define the map
$$
\Cal F_{\tau}: C^{\infty}(M, \R) \to C^{\infty}(\Cal B^{hw,r}_S, \R)
$$
and prove by the direct computation that this map preserves both the Poisson structures.
First time it was done in preprint [7] and then published as [28].

In Section 3 we present the main technical result: we prove that the image of $\Cal F_{\tau}$
consists of quasi symbols. The proof is based on what we call "dynamical correspondence".
This means briefly speaking that for every $f \in C^{\infty}(M, \R)$ the differential 
of $\Cal F_{\tau}$ maps the Hamiltonian vector field $X_f$ to the Hamiltonian vector field
$X_{F_f}$ (up to a constant).  It means that the correspondence is "dynamical":
it is compatible with classical dynamics of the system. The results of this section
were established in preprint [29] and published as [30].

In Section 4 we apply the mathematical results to introduce ALG(a) - quantization
of classical mechanical systems. As the first two steps we show the reduction
of this method in presence of either complex or real polarization. Technically 
it requires the consideration of critical points of the quasi symbols. Further we discuss
the corresponding quasi classical limit.

In the last small section we add some observations linking ALG(a) - quantization with super geometry.
This section is the smallest one but it looks like quite profitable in a future.
While the previous sections have been finished (moreless) the last one consists of
remarks and observations which could be exploited in the further investigations.

\subheading{Asknolegements}  The author's work on this theme was started in 1999 at
Max - Planck - Institute for Mathematics (Bonn) when, just by a chance,
two preprints ([3] and [9]) were collected by him. The main results
of the text were established in MPI (Bonn), Korean Institute of Advanced Study
(Seoul) and Joint Institute for Nuclear Research (Dubna) so first of
all the author would like to thank all the people from the institutions
for the hospitality, the friendel attention and very good working conditions.
Personally  I would like to thank Prof. I.R. Shafarevich
as the leader of the scientific Al school (to which  the author
belongs) and as an attentive listener whose remarks were extremely
important and useful for the work. 
I would like to thank Prof. A.G. Sergeev and all the participants of his
seminar in Steklov Mathematical Institute (Moscow) for the constant attention
to the work during these years and valuable discussions and remarks.
I would like to thank Prof. A.S. Tikhomirov for the interest to the subject
and a number of important remarks. I would like to cordially thank
Professors Vik. Kulikov, A. Gorodentsev, V. Pidstrigach, B. Kim,
S. Kuleshov, D. Orlov, B. Karpov, A. Isaev, and
especially H. Khudaverdyan and P. Pyatov. Additionally I thank
my family for the hard assistance which I get all the time.
Thus the double appreciation (and the first one) 
should be expressed to Prof. A.N. Tyurin.

\newpage

\head Section 1. Geometric quantization and its geometric formulation
\endhead

$$
$$

Quantum states of a quantum mechanical system  are represented by rays in the corresponding Hilbert space (the space of wave functions, see f.e. [1]). The space of the rays is a projective space
endowed with a Kaehler structure. Turning to the projective geometry language one gets
the postulates of quantum mechanics in pure geometrical terms. Applications of this idea
in quantization problem are quite well known. Here we follow [3], [18] as
 modern and complete texts on the subject. The aim is to formulate the quantum mechanics postulates
without references on Hilbert spaces and linear structures. The basic notions of symplectic geometry can be found f.e. in [2].
Further in this section we recall three main constructions in geometric quantization together with
basic notions of it. We begin with the definition of the prequantization data and then 
remind the basic Souriau - Kostant construction.  It needs an additional specification since it gives
a reducible representation of the Poisson algebra. This specification comes with an additional choice of
an appropriate polarization (if the given manifold admits this choice). We list then two
"extremal" cases when it admits complex or real polarization. Recalling these known constructions we place some additional colors coming when one takes in mind the geometric reformulation.
$$
$$
$$
$$

\subhead 1.1. Hilbert space as a Kaehler manifold
\endsubhead

Let us consider a Hilbert space $\Cal H$ as a real space endowed with a complex structure $J$. It is a real linear operator such that
$$
J^2 = - id.
\tag 1.1.1
$$
At the same time the hermitian pairing can be decomposed on the real and imaginary parts:
$$
<\Phi, \Psi> = \frac{1}{2 h} G(\Phi, \Psi) + \frac{i}{2h} \Om(\Phi, \Psi).
\tag 1.1.2
$$
The standard properties of  hermitian products imply that $G$ is a symmetric positive pairing while
$\Om$ is skew symmetric and both the pairing are related by the complex operator $J$. The compatibility
condition reads as
$$
G(\Phi, \Psi) = \Om(\Phi, J \Psi)
\tag 1.1.3
$$
thus one has a Kaehler structure on the real space $\Cal H$. Therefore any Hilbert space 
can be described as a real space endowed with a Kaehler triple $(G, J, \Om)$.
So this space is simulteniously a symplectic space and a complex space. As a symplectic space it
has the corresponding Poisson bracket for the functions.
In classical mechanics the observables are smooth real functions which are defined (up to constant)
by the Hamiltonian vector fields (we are dealing with  connected manifolds only). In the quantum case 
one could represent the observables by vector fields as well: for any linear self adjoint  operator $\hat F$
(a quantum observable) let us call the vector field
$$
Y_{\hat F}(\Psi) = \frac{1}{h} J \hat F \Psi
\tag 1.1.4
$$
 as Schroedinger. The point is that since in the real setup the Schroedinger equation reads as
$$
\dot \Psi = - \frac{1}{h} J \hat H \Psi
\tag 1.1.5
$$
for a Hamiltonian operator $\Cal H$ then it can be shortly rewritten due to the given notation.
It's not hard to establish some appropriate properties of the Schroedinger  fields. First, we know that quantum observable $\hat F$ generates one - parameter family of unitary transformations of the Hilbert space. Since
$Y_{\hat F}$ is the generator of the family this implies that the Schroedinger vector field
preserves all three structures which the Kaehler triple comprises. Therefore field $Y_{\hat F}$ is locally Hamiltonian and taking in mind linearity of $\Cal H$ it is globally Hamiltonian as well. The corresponding function
is very well known in the standard quantum mechanics. It is the expectation value given by the formula
$$
F(\Psi) = < \Psi, \hat F \Psi>
\tag 1.1.6
$$
for operator $\hat F$. Then if $\eta$ is a tangent vector at point $\Psi$ we have
$$
\aligned
d F_{\eta} (\Psi) = \frac{d}{dt} <\Psi + t \eta; \hat F(\Psi + t \eta)> |_{t=0} \\
= <\Psi, \hat F \eta> + <\eta, \hat F \Psi> = \frac{1}{h} G(\hat F \Psi, \eta) \\
= \Om(Y_{\hat F}(\Psi), \eta) = \imath_{Y_{\hat F}(\Psi)} \Om (\eta),\\
\endaligned
\tag 1.1.7
$$
using self duality of $\hat F$ together with the definitions of $G$, $\Om$ and $Y_{\hat F}$.
This gives us a well known fact that the evolution in time of any quantum mechanical system 
can be described in terms of Hamiltonian mechanics; the corresponding Hamiltonian function
equals to the expectation value of the Hamiltonian operator.

Second, let $\hat F, \hat K$ be two quantum observables with expectation values $F, K$. A short computation
gives
$$
\{F, K\}_{\Om} = \Om(X_F, X_K) = <\frac{1}{i h} [\hat F, \hat K]>,
\tag 1.1.8
$$
where we use brackets $<\quad >$ for the expectation value (1.1.6). This means that the correspondence
$$
\hat F \mapsto F
\tag 1.1.9
$$
is a homomorphism of Lie algebras. Its kernel is trivial thus it is an isomorphism 
of the Lie algebra of self adjoint operators and a subalgebra of the Poisson algebra
over the Kaehler space. The condition on real functions distinguishing
this subalgebra  can be derived evidently what we'll do in the next subsection.

Third, one has an additional ingredient in the picture --- the riemannian metric $G$.
It defines a real pairing of Hamiltonian vector fields which corresponds to the Jordan product
in the standard quantum mechanics
$$
\{ F, K\}_+ = \frac{h}{2} G(X_F, X_K) = <\half [\hat F, \hat K]>.
\tag 1.1.10
$$
One could take the first equality in (1.1.10) as the definition 
of the riemannian bracket for $F$ and $K$. This riemannian bracket has not any classical
analog being an ingredient of pure quantum theory. F.e., it could be exploited to
define what one calls "uncertainty" of a quantum observable in a quantum state with unit norm namely
$$
(\De \hat F)^2 = <\hat F^2> - <\hat F>^2 = \{F, F \}_+ - F^2.
\tag 1.1.11
$$
And the famous Heisenberg uncertainty relation looks like
$$
(\De \hat F)^2 (\De \hat K)^2 \geq 
(\frac{h}{2} \{F, K\}_{\Om})^2 + (\{F, K\}_+ - FK)^2
\tag 1.1.12
$$
in these terms.

\subhead 1.2. The projectivization
\endsubhead

As it was already mentioned the real space of quantum states is represented by rays in the corresponding Hilbert space. Collecting the rays one get the projective space $\Cal P$ which is a Kaehler manifold
itself. The procedure giving us the desired structure over $\Cal P$ is usually called either 
"the Dirac theory of constraints" or  "the Kaehler reduction" (it depends on the context). One has
the natural $U(1)$ - action (phase rotations) on $\Cal H$, preserving all the structures.
As well it preserves the function
$$
C(\Psi) = <\Psi, \Psi> - 1,
\tag 1.2.1
$$
which is usually called  either "the Dirac constraint function" or "momentum map"
since it is an equivariant map acting  from the Hilbert space to the Lie co - algebra
of $U(1)$. After  this is understood we just apply the standard procedure of
Kaehler reduction getting as the result the projective space endowed with
the corresponding Kaehler structure. In other words since the time evolution preserves
the  levels of the constraint function we get a first constraint system; f.e.,  unit sphere
 $S \in \Cal H$ (the level $C = 0$) is preserved by the motion generated by
$C$ itself:
$$
L_{X_C} C = \{C, C \}_{\Om} = 0.
\tag 1.2.2
$$
For a first constraint system there is a gauge freedom and the corresponding gauge transformations
are generated by the Hamiltonian flow of function $C$. In our case these gauge directions
is defined by the Schroedinger vector field
$$
X_C(\Psi) = - \frac{1}{h} J \Psi.
\tag 1.2.3
$$
Thus the vector field
$$
\Cal J_C = h X_C|_S
\tag 1.2.4
$$
generates the phase rotations. Factorizing by the phase rotations one gets the real
space of quantum states (sometimes called as "reduced phase space"). To emphasize 
its physical role and geometrical structure we would like to call it
 quantum phase space. 

To define the symplectic structure on $\Cal P$ consider the pair of maps
$$
\aligned
\imath: S \to \Cal H, \\
\pi: S \to \Cal P, \\
\endaligned
\tag 1.2.5
$$
where the first one is the inclusion and the second --- the projection. On the unit sphere
$S$ we have closed 2- form $\imath^* \Om$ which degenerates exactly in the rotation directions. Moreover
this 2 - form is constant along the directions (since these directions are given by
the Hamiltonian vector field of the constraint function). Therefore
there exists a closed non degenerated 2 - form $\pi_* \imath^* \Om$
on $\Cal P$ (and moreover it is symplectic in strong sense for the infinite dimensional case). Let us remark
that we've described nothing but the mechanism of the symplectic reduction valid in
the infinite dimensional case (see f.e. [8]). 

Now what happens with quantum observables during this procedure? If $\hat F$ is a bounded self adjoint
operator then the corresponding expectation value is a good smooth function being
restricted to the unit sphere. First, the restriction $\imath^* F$ is evedentely invariant
under the gauge transformations therefore there exists a smooth function $f$ on $\Cal P$
such that
$$
\pi^* f = \imath^* F.
\tag 1.2.6
$$
Second, the relation between $F$ and $f$ is "dynamical" ---  Hamiltonian vector field 
$X_F$ projects exactly to Hamiltonian vector field $X_f$. We'll come back
to this dynamical correspondence later on, presenting now the following
\proclaim{Definition} Let $f: \Cal P \to \R$ be a smooth real function on the projectivization
of a Hilbert space $\Cal H$. Then we call it a symbol iff there exists a bounded self - adjoint
operator $\hat F$ on $\Cal H$ such that $f = \pi_* \imath^* <\hat F>$.
\endproclaim

\subheading{Remark} This notion was introduced by F. Berezin so we'd like to follow him
in the terminology.

Now coming back to the Hamiltonian vector fields let us note that any $\hat F$ commutes with
$Id$ thus any expectation value $F$ commutes with the constraint function so
$$
L_{X_C} F = 0.
\tag 1.2.7
$$
Therefore $X_F$ is constant along the integral curves of $\Cal J$ and one can push it down to
$\Cal P$, getting exactly $X_f$. Going further let's take the Poisson bracket induced by the reduced
symplectic form $\Om_p$. For any expectation values  $F, K: \Cal H \to \R$ one has
$$
\aligned
\pi^*\{f, k \}_{\Om_p} & = \pi^*(\Om_p(X_f, X_k)) \\
= \Om_p(\pi_* X_F, \pi_* X_K) & = \Om(X_F, X_K)|_S = \{F, K \}|_S \\
\endaligned
\tag 1.2.8
$$
by the projection formula where $f, k: \Cal P \to \R$ are the corresponding reduced functions.

Thus for any quantum observable there exists the corresponding symbol. Moreover, the correspondence is
one - to - one being an isomorphism of Lie algebras. At the same time we call this correspondence
"dynamical" because it transfers the Schroedinger dynamics exactly to the Hamiltonian dynamics
on the reduced quantum phase space.

\subhead 1.3. Riemannian geometry and the measurement 
\endsubhead 

As we have seen above the quantum phase space carries an additional structure
which is unusual in classical mechanics. The riemannian metric responds to the probability aspects
of the quantum theory.

The riemannian metric over $\Cal P$ is given during the same Kaehler reduction procedure. The restriction
$\imath^* G$ to $S$ is a non degenerate riemannian metric. Since every Schroedinger vector field preserves
whole hermitian structure then $\Cal J$ is a Killing vector field for $G$ over $S$:
$$
L_{\Cal J} (\imath^* G) = 0.
\tag 1.3.1
$$
Thus $\imath^* G$ is constant along the integral curves of $\Cal J$ being non degenerated. Taking
$$
\tilde G_p = (G - \frac{1}{2h}(\Psi \otimes \Psi + \Cal J \otimes \Cal J))|_S
\tag 1.3.2
$$  
we get a symmetric tensor equals to $\imath^* G$ on the transversal to $\Cal J$ subspaces and which degenerates
exactly in the directions of $\Cal J$. Therefore $\tilde G_p$ can be factorized by the phase transformations 
which gives us the corresponding riemannian metric over $\Cal P$. Now we get  the Kaehler structure
over $\Cal P$. The question is how one can recognize the symbols over $\Cal P$? One has the following
\proclaim{Proposition 1.1 ([3])} Function $f: \Cal P \to \R$ is a symbol iff its Hamiltonian vector field
preserves the riemannian metric:
$$
L_{X_f} G = 0
\tag 1.3.3
$$
(thus at the same time it is a Killing vector field).
\endproclaim

This means that as in classical dynamics the space of observables in the quantum theory
consists of the smooth functions whose Hamiltonian vector fields give the infinitesimal symmetry
of the given structure. At the same time one sees that in contrast with the classical case
the space of quantum observables is small: f.e., in  finite dimension one gets a finite dimensional
space of quantum observables. Necessity of (1.3.3) has been explained so it remains to
construct the corresponding operator on $\Cal H$ for every symbol $f$ on $\Cal P$.
The correspondence comes with Berezin - Rawnsley quantization scheme (see [16]) namely
let's take holomorphic line bundle $\Cal O(1)$ over $\Cal P$. The space of its holomorphic sections
is dual to $\Cal H$.  Taking any appropriate  hermitian structure on $\Cal O(1)$, compatible with the given riemannian metric,
one gets:

a self adjoint operator $Q_f$ as the result of Berezin - Rawnsley quantization;

an isomorphism between the dual spaces.

Using this isomorphism one transfers $Q_f$ from $H^0(\Cal O(1), \Cal P)$ to $\Cal H$ getting
what one needs. It's not hard to see that the result doesn't depend on the choice
of the hermitian structure. We'll recall the Berezin - Ransley quantization method
further in subsection 1.5.

The notion of uncertainty can be transfered in natural way to the projective space.
Let's define riemannian bracket of the form
$$
(f, k) = \frac{h}{2} G_p (X_f, X_k).
\tag 1.3.4
$$
Since these $f, k$ correspond to expectation values $F, K: \Cal H \to \R$ then (see (1.3.2))
$$
\aligned
\{F, K \}_+ = \frac{h}{2} G(X_F, X_K) = \frac{h}{2} \tilde G_p (X_F, X_K) + \\
\frac{1}{4} G(\Cal J, X_F) G(\Cal J, X_K) = \pi^*((f, k) + fk).\\
\endaligned
\tag 1.3.5
$$
Therefore the symbol which corresponds to the Jordan product of $\hat F$ and $\hat K$ is 
$$
\{f, k\}_+ = (f, k) + fk.
\tag 1.3.6
$$
It's reasonable to call this expression as "symmetric bracket" of $f$ and $k$. Note that the riemannian bracket
of two symbols is not in general a symbol but it carries a natural physical sense:
it is exactly  the function of quantum covariance. In particular, $(f, f)(p)$ is the square of the uncertainty
of observable $f$ at state $p$:
$$
(\De f)^2 (p) = (\De \hat F)^2 (\pi^{-1} p) = (f, f)(p).
\tag 1.3.7
$$
Uncertainty relation then reads as
$$
(\De f) (\De k) \geq (\frac{h}{2} \{f, k \}_{\Om_p})^2 + (f, k)^2.
\tag 1.3.8
$$
One gets its standard form 
$$
(\De f) (\De k) \geq (\frac{h}{2} \{f, k\}_{\Om_p})^2,
\tag 1.3.9
$$
at the states where the Hamiltonian vector fields $X_f, X_k$ are conjugated up to sign by the complex structure
on $\Cal P$. Therefore this quantum covariance $(f, k)(p)$ measures the "coherence" of state $p$ with respect to
observables $f$ and $k$.

Consider the probabilistic aspects. Let $\Psi_0$ be a vector with unit norm in the original
Hilbert space $\Cal H$ and $p_0 = \pi(\Psi) \in \Cal P$. Then one defines a function
$d_{\Psi}$ on $S$ in  natural way:
$$
d_{\Psi_0} (\Psi) = |<\Psi_0, \Psi>|^2
\tag 1.3.10
$$
for every $\Psi$ with unit norm. Since $d_{\Psi_0}$ doesn't depend on the gauge transformations 
it can be push down to $\Cal P$ which gives us a function $d_{p_0}$ on the projective space
defined by
$$
d_{p_0}(p) = d_{\Psi_0}(\pi^{1}(p)).
\tag 1.3.11
$$
If the quantum mechanical system stays at the state $p_0$ then the probabilistic
distribution is given by $d_{p_0}$. It's well known that this function has an expression in terms of geodesic lines, namely
\proclaim{Proposition 1.2 ([3])} For any $p_0, p \in \Cal P$ there exists a closed geodesic
line passing from $p_0$ to $p$ such that
$$
d_{p_0}(p) = cos^2(\frac{1}{\sqrt{2h}} \sigma(p_0, p)),
\tag 1.3.12
$$
where $\sigma$ is the corresponding geodesic distance.
\endproclaim

To prove it in the finite dimensional case it's sufficient to
recall what the Fubini - Study metric on a projective space is
(see f.e. [11]). To reduce the infinite dimensional case to the previous
one it's sufficient to note that for any points $p_0, p \in \Cal P$
there exists a projective plane $\proj^2 \subset \Cal P$
which contains both of them. Then it's clear that the desired geodesic line
in any case lies in this $\proj^2$.

The translation of all the measurement aspects to the projective language
 needs two more steps. First, 
\proclaim{Proposition 1.3 ([3])} Vector $\Psi$ is an eigenvector of quantum observable $\hat F$
with eigenvalue $\la$ iff the corresponding point $p = \pi(\Psi)$ is a critical one
of the corresponding symbol $f$ with critical value $\la$.
\endproclaim
The proof is evident. Second, let $\hat F$ be a quantum observable with an arbitrary
spectrum. To define the spectrum of the corresponding symbol $f$ we use (1.3.4) and
(1.3.6), getting the following
\proclaim{Definition} Spectrum $sp(f)$ of symbol $f$ consists of real numbers such that
 the function
$$
\aligned
n_{\la}: \Cal P \to \R \cup \{\infty\},\\
n_{\la}: p \mapsto ((\De f)^2 (p) +(f(p) - \la)^2)^{-1} \\
\endaligned
$$
is not bounded.
\endproclaim
The point $\la$ where $n_{\la} = \infty$ is a critical value of $f$.
Now turn to the spectral projectors. If $\La$ is a continuous component 
of $sp(f)$ and $P_{\hat F, \La}$ is the projector which corresponds to $\hat F$ and $\La$
then let's denote as $\Cal E_{f, \La}$ the following projectivization
$$
\Cal E_{f, \La} = \proj (Im P_{\hat F, \La}) \subset \Cal P.
$$
This submanifold corresponds to the critical set with critical volume 1 of the symbol,
induced by projector $P_{\hat F, \La}$ (notice that every such projector is a
bounded self adjoint operator). With submanifold $\Cal E_{f, \La}$ we can take
the geodesic projection instead of the spectral projector. In this terms
projector $P_{f, \La}$ maps point $p \in \Cal P$ to the point of $\Cal E_{f, \La}$ which is the nearest one with respect to the geodesic distances.

 \subhead 1.4. Postulates of quantum mechanics
\endsubhead

Summing up we get the following picture. One has a projective space $\Cal P$ which is a Kaehler manifold
equipped with the corresponding Kaehler structure. Thus $\Cal P$ has a fixed symplectic structure
which defines the corresponding Lie algebra structure on the function space and governs
the evolution of the system. But there are two major differences between our picture and the classical mechanics case.
First, the quantum phase space has very special nature as a Kaehler manifold being a projective space (we will
speak about possible generalizations later). Second, as a Kaehler manifold, it is equipped with a riemannian metric
and this last one governs the measurement process.  This ingredient was absent in the classical theory ---
in quantum theory it responds for such notions as uncertainty, state reduction and so on.

We list below short vocabulary just to summarize the translation.

\subheading{Physical states} Physical states of quantum system correspond to point of an appropriate Kaehler manifold
(which is in the basic example a projective space).

\subheading{Kaehler evolution} The time evolution of physical system is defined by a flow over $\Cal P$ which preserves 
whole the Kaehler structure. This flow is generated by a vector field which is dense everywhere on $\Cal P$. 

\subheading{Observables} Physical observables are given by special real smooth functions over $\Cal P$
whose Hamiltonian vector fields preserve the Kaehler structure. In other words physical observables
are presented by symbols.

\subheading{Probability aspects} Let $\La \subset \R$ is a closed subset of $sp(f)$ and the system stays at a 
state represented by point $p \in \Cal P$. The probability to get the result after the measurement process which
belongs to $\La$ is given by the formula
$$
\de_p(\La) = cos^2(\frac{\si(p, P_{f, \La}(p))}{\sqrt{2h}},
$$
where $P_{f, \La}(p)$ is the nearest to $p$ point of $\Cal E_{f, 
La}$.

\subheading{Reduction} The choice of an arbitrary closed subset  $\La$ in $sp(f)$ defines the ideal measurement which can be performed. This measurement answers to the question whether or not the volume $f$ belongs to $\La$.
After the measurement process is performed the state of the system is defined by either $P_{f, \La}(p)$
or $P_{f, \La^c}(p)$ in dependence on the measurement  result.

Now we need two direct citation from [3] to clarify the ways of possible generalizations. First, let's note that
the postulates of quantum mechanics can be formulated on pure geometrical language without any references to
Hilbert spaces. Of course the standard Hilbert space considerations and related algerbaic machinery
equip one with a good setup for concrete computations. But mathematically the situation is 
very similar to the usual geometric consideration of compact manifolds with nontrivial topology:
in practice it's convenient at the starting point of investigations to include this manifold
to an appropriate ambient space (Euclidean, projective, ...) and then study it. But the inclusion
is needed only for convenience; one could derive all what one needs directly from the geometry
of this manifold. Second, the quantum mechanics linearity could be an analogy of the inertial systems
in special relativity and the geometric formulation of quantum mechanics may be equivalent
to Minkowski formulation of special relativity and as the last one showed the way to general relativity
the geometric formulation would lead us to a new theory. One way generalizing quantum mechanics can be presented 
immideately after the geometric formulation is understood. One can assume that there are exist some other
Kaehler manifolds (not only projective spaces) which could carry quantum mechanical systems.  The dynamical properties could be  easily satisfied (since classical mechanics allows one to consider some other symplectic manifolds
than only projective spaces). The first question arising in this way is the question of observables.
Since there are exist the Kaehler manifolds which doesn't admit any real functions whose Hamiltonian vector field
would keep the Kaehler structure the question is really of the first range. F.e., one could require for
the Kaehler manifold to admit the maximal possibility --- it corresponds to the Kaehler manifolds of
constant holomorphic sectional curvature (see [3]). It is well known that in finite dimension this condition
is satisfied only by the projective spaces. In the infinite dimensional case this problem is still open
so one could conjecture that there are infinite dimensional Kaehler manifolds satisfying the condition
and different from the projective spaces. We would say here that this requirement is too strong:
it is sufficient to impose the condition that the space of quantum observables allowed over a tested
Kaehler manifold {\it is sufficiently large}. So the Kaehler structure over this manifold has to supply us
with a good reserve, which we could use in our investigations. So during  the text we will keep in mind 
this conception. Now we continue with the natural 
\proclaim{Definition} For a Kaehler manifold $\Cal K$ let's call a real smooth  function $f$ "quasi symbol"
if its Hamiltonian vector field preserves whole the Kaehler structure.
\endproclaim

Let's denote the space of all quasi symbols as $\Cal C^{\infty}_q (\Cal K, \R)$.
A manifest property of such functions comes immideately
\proclaim{Proposition 1.4} For any Kaehler manifold $\Cal K$ the space $C^{\infty} (\Cal K, \R)$ is a Lie subalgebra
of the Poisson algebra.
\endproclaim

To prove it one takes the Poisson bracket of any two quasi symbols, ensures that the corresponding Hamiltonian vector
field is proportional to the commutator of the Hamiltonian vector field of the given functions and then
differentiates the given riemannian metric in the commutator direction. The answer is obvious.

Further we will construct a nonlinear generalization of quantum mechanics. The nonlinearity means that
our exploited  Kaehler manifold is not a projective space. But in Section 4 we will show
that this Kaehler manifold admits a large space of quasi symbols sufficient to include
all smooth functions from some given symplectic manifold. Since we would like to speak concretely about
the quantization of classical mechanical systems it is quite enough for us. Therefore we omit in this text any discussions about pure generalizations of quantum mechanics which should be discovered in this way.
So from this moment we speak about quantization of classical mechanical systems.

\subhead 1.5. Souriau - Kostant quantization
\endsubhead 

Let's remind first of all what the problem is. If one considers a classical mechanical system
with finite degrees of freedom then the phase space is represented by a symplectic manifold $M$ with 
symplectic form $\om$. Classical observables are given by smooth real functions over $M$
and the corresponding  Poisson bracket $\{ , \}_{\om}$ defines the Lie algerba structure on the space of observables. One understands then the ideal quantization as a procedure which results
a Hilbert space $Q_M$ together with a correspondence
$$
q: C^{\infty} (M, \R) \to Op(Q_M),
\tag 1.5.1
$$
satisfying the usual properties (the terminology can vary: while one takes the complete list
of Dirac conditions then it's called "quantization"; if one relaxes the requirements leaving
the irreducibility condition then it is often called "prequantization" but we'd like
to reserve the last word for the notion "prequantization data" reminded below). According to a 
celebrated van Hove theorem ([13]) one can't realize the ideal quantization
using the geometry of $(M, \om)$ itself (since as we'll see usually one takes 
the space of functions or the space of sections of some appropriate bundle over $M$ as the Hilbert space).
As we will see the known examples of quantization leave unsatisfied one or two conditions
from the complete Dirac list. The constructions begin with so called "prequantization data" (see [10]).
First of all, the cohomology class of $\om$ has to be integer (some times 
one calls this condition "the Dirac condition of the charge integrality"):
$$
[\om] \in H^2(M, \Z) \subset H^2(M, \R).
\tag 1.5.2
$$
If this condition holds then there exists a topological complex vector bundle $L$ of rank 1
over $M$ characterized uniquely by  the condition $c_1(L) = [\om]$ where
$c_1$ is as usual the first Chern class. One can choose and fix
any appropriate hermitian structure on $L$, getting a $U(1)$ - bundle. Then
one gets another prequantization datum --- a hermitian connection $a \in \Cal A_h(L)$
satisfying the following natural condition
$$
F_a = 2 \pi i \om,
\tag 1.5.3
$$
where $F_a$ is the curvature form. If our given symplectic manifold $M$ is simply connected
then this prequantization connection is unique up to the gauge transformations.
In  other case the space of the equivalence classes of solutions is isomorphic to
$H^1(M, \R)$. In this text we are interested mostly in the simply connected case
so in the rest of the paper we'll work with compact simply connected 
symplectic manifolds only. All remarks about non simply
connected case are more than welcome.

Whatever the case we are trawling through there is an uncertainty in the choice of
the prequantization data: even if $M$ is simply connected then we are free to take
any hermitian structure on $L$ and this implies that the prequantization connections can be different. 
Recall briefly how two different hermitian structures could be compared. For any
pair $h_1, h_2$ where $h_i$ is a hermitian structure on $L$ one has a real positive function
$e^{\phi(x)}$ satisfying the following property
$$
\forall v_1, v_2 \in \Ga(M, L) \quad \quad <v_1, v_2>_{h_1} = e^{\phi(x)}<v_1, v_2>_{h_2}
\tag 1.5.4
$$
(since $M$ is simply connected!), and in terms of this function we can relate hermitian connections
compatible with $h_1$ and $h_2$. Consider the spaces $\Cal A_{h_1}(L)$ and $\Cal A_{h_2}(L)$
consist of all hermitian connections compatible with $h_1$ or $h_2$ respectively.
These are two affine subspaces in the ambient space of all $C^*$ - connections on $L$
which we denote as $\Cal A(L)$. One has the following
\proclaim{Proposition 1.5}  1) The affine spaces $\sA_{h_i}(L) \subset \sA(L)$ either do not intersect each other  or coincide;

2) two connections $a_i \in \sA_{h_i}(L)$ differ by a complex 1 - form $i \rho + \half d \phi$,
where  function $\phi$ is defined in (1.5.4).
\endproclaim

Really,  connection $a$ is compatible with hermitian structure $h$ if for any two sections
$s_1, s_2 \in \Ga(L)$ the following identity holds:
$$
d<s_1, s_2>_{h} = < d_a s_1, s_2>_h + <s_1 , d_a s_2>_h
\tag 1.5.5
$$
(see, f.e., [8]). If connection $a_0$ is compatible simultaneously with $h_1$ and $h_2$:
$$
a_0 \in \Cal A_{h_1}(L) \cap \Cal A_{h_2}(L),
$$
then using (1.5.5) one gets
$$
d<s_1, s_2>_{h_1} = d<s_1, s_2>_{h_2}
\tag 1.5.6
$$
for any sections $s_1, s_2 \in \Ga(L)$. But it implies that
$$
d e^{\phi} <s_1, s_2>_{h_1} = 0
\tag 1.5.7
$$
for any sections $s_1, s_2 \in \Ga(L)$. It means that $d \phi$ vanishes
and $\phi = const$. Coming back to (1.5.5) one sees that in this case 
$\Cal A_{h_1}(L) = \Cal A_{h_2}(L)$. This gives the first statement of the proposition.

Further, using once more (1.5.5), we have 
$$
\aligned
d<s_1, s_2>_{h_2} = d e^{\phi} <s_1, s_2>_{h_1} \\
= e^{\phi} d \phi <s_1, s_2>_{h_1} + e^{\phi} d<s_1, s_2>_{h_1} \\
= e^{\phi} <s_1, s_2>_{h_1} d \phi + e^{\phi}<d_{a_1}s_1, s_2>_{h_1}\\
+ e^{\phi}<s_1, d_{a_1} s_2>_{h_1} = e^{\phi}(<s_1, s_2>_{h_2} d \phi \\
<d_{a_1} s_1, s_2>_{h_2} + <s_1, d_{a_1}s_2>_{h_2}).\\
\endaligned
\tag 1.5.8.
$$
Comparing the last expression with (1.5.5) for $i = 2$, we get
$$
<\de s_1, s_2>_{h_2} + <s_1, \de s_2>_{h_2} = <s_1, s_2>_{h_2} d \phi,
\tag 1.5.9
$$
where $\de = d_{a_1} - d_{a_2}$ is a complex 1 - form. Its real part is fixed by
(1.5.9) while the imaginary part is arbitrary. The second part of the proof is over.

Therefore if one takes two affine subspaces in $\Cal A(L)$ corresponding to 
hermitian structures $h_1, h_2$ on $L$ then there exists a natural affine map
$$
\aligned
\Cal A_{h_1}(L) \to \Cal A_{h_2}(L) \\
a_1 \mapsto a_2 = a_1 + \half d \phi,\\
\endaligned
\tag 1.5.10
$$
and this is true only in the simply connected case. It reflects the fact that
in absence of additional structures over an arbitrary 
(not necessary simply connected) manifold $X$ there exists unique cohomology class
in $H^1(X, \R)$ 
which has a distinguished representation by closed 1 -form. This class is $[0]$
and the distinguished closed form is zero 1- form.
On the other hand, if we fix any $C^*$ - connection which solves equation (1.5.3)
then there exists a hermitian structure on $L$ such that this connection is compatible
with this structure (in our simply connected case).
Thus one can estimate the uncertainty in the choice of the prequantization data.
So let us choose an appropriate structure (or connection). In Souriau - Kostant quantization
one takes the space $L^2(M, L)$ as the space of wave functions: the fixed hermitian structure
on $L$ together with the Liouville volume form defines hermitian structure
$$
<s_1, s_2>_q = \int_M <s_1, s_2>_h d \mu_L
\tag 1.5.11
$$
on the space of sections $\Ga(L)$, and one takes the completion of the space with respect to
this hermitian structure $q$. 
 Then  
$$
\Cal H = \overline{\{ s \in \Ga(L) | \int_M <s, s>_h d \mu_L < \infty\}}.
\tag 1.5.11'
$$
Every function $f$ induces the following operator
$$
Q_f: \Cal H \to \Cal H \quad \quad | \quad Q_f s = \nabla_{X_f} s + 2 \pi i f \cdot s
\tag 1.5.12
$$
on the wave function space $\Cal H$. The commutator of two such operators can be directly 
computed:
$$
\aligned
[Q_f, Q_g] = (\nabla_{X_f} + 2 \pi i f)(\nabla_{X_g} + 2 \pi i g)\\
- (\nabla_{X_g} + 2 \pi i g)(\nabla_{X_f} + 2 \pi i f) =\\
\nabla_{X_f} \nabla_{X_g} - \nabla{X_g}\nabla_{X_f} + 2 \pi i \nabla_{X_f} g - 2 \pi i \nabla_{X_g} f \\
= \nabla_{X_f} \nabla_{X_g} - \nabla{X_g} \nabla_{X_f} + 4 \pi i \{f, g\}\\
= \nabla_{[X_f, X_g]} + R(X_f, X_g) + 4 \pi i \{f, g \} \\
= \nabla_{X_{\{f, g \}}} - 2 \pi i \{f,  g\} + 4 \pi i \{f, g \}\\
= \nabla_{X_{\{f, g \}}} + 2 \pi i \{f, g\}.\\
\endaligned
\tag 1.5.13
$$
This means that such operators could be exploited for the quantization. The point is that
such $Q_f$ is not self adjoint  --- conversely, it is a unitary operator.
Really, for any two wave functions we have
$$
\aligned
<Q_f s_1, s_2>_q = \int_M<\nabla_{X_f} s_1 + 2 \pi i f s_1, s_2>_h d \mu_L = \\
\int_M <d_a s_1, s_2>_h(X_f) d \mu_L - \int_M<s_1, s \pi i s_2>_h d \mu_L = \\
\int_M d<s_1, s_2>_h(X_f) d \mu_L - \int_M<s_1, d_a s_2>(X_f) d \mu_L - \\
\int_M<s_1, 2\pi i s_2>_h d \mu_L = \int_M \{f, g \} d \mu_L - \int_M<s_1, \nabla_{X_f} s_2>_h d \mu_L\\
- \int_M <s_1, 2 \pi i s_2>_h d \mu_L = - \int_M <s_1, Q_f s_2>_h d \mu_L = - <s_1, Q_f s_2>_q,\\
\endaligned
\tag 1.5.14
$$
where 
$$
g = <s_1, s_2>_h \in C^{\infty}(M, \R).
$$
To correct the picture let's take the dual bundle $L^*$ with connection $a'$ whose
curvature is equal to $- 2 \pi i \om$. Then it changes the sign in (1.5.13) so
$$
[Q_f, Q_g] = - Q_{\{f, g\}}.
$$
Now if we consider
$$
\hat Q_f = i Q_f,
\tag 1.5.12'
$$
then these operators should be self adjoint and the correspondence principle remains to
be satisfied. This way approaches us to the Souriau - Kostant quantization (see f.e. [13], [14], [21],
[33]).

What we could derive from this construction taking in mind some dynamical properties of the given system?
Really smooth functions on $(M, \om)$ correspond (up to constant) to infinitesimal
deformations of the symplectic manifold.  Add  now the prequantization data: let P be the principle
$U(1)$ - bundle, associated with line bundle $L$, and $A \in \Om^1_P(i \R)$ --- the corresponding connection.
Its differential equals just to $i \om$. For any smooth function $f$ consider the lifting
$$
Y_f = X_f + g \cdot \frac{\partial}{\partial t}
\tag 1.5.15
$$
of the Hamiltonian vector field to a vector field which preserves our hermitian connection $A$.
There $g$ is a smooth function on $M$ and 
$$
\frac{\partial}{\partial t} = \theta_t
$$
is the generator of the canonical circle action along the fibers of $\pi: P \to M$.
The condition that the lifting $Y_f$ preserves $A$ means that
$$
Lie_{Y_f} A = 0,
$$
thus one gets
$$
\aligned
Lie_{Y_f} A = d(A(Y_f)) + \imath_{Y_f} dA = \\
\pi^* id(g) + \pi^* \imath_{X_f}(i \om) = 0.\\
\endaligned
\tag 1.5.16
$$
Therefore one can relate $f$ and $g$, since (1.5.16) implies
$$
dg = - \imath_{X_f} \om,
\tag 1.5.17
$$
and consequently $g$ coincides with $f$ up to constant. This means that every
Hamiltonian vector field can be lifted to $P$ almost canonically.
So any symplectomorphism of the base $(M, \om)$ can be lifted to an automorphism
of the principal bundle $P$, preserving $A$, up to canonical $U(1)$ -
transformations. The constant which would be add to $g$ in (1.5.17) just responds to
an infinitesimal canonical transformation of $P$. To avoid  the uncertainty
let's take the projectivization of the wave function space
$$
\Cal P = \proj(\Cal H).
$$
Then any Hamiltonian vector field $X_f$ can be lifted uniquely which gives a vector field
on $\Cal P$. Denote this unique vector field $\tilde Y_f$ 
and consider what are the properties distinguishing this one. The main property is
that this field is defined {\it dynamically}. It means the following:
since any symplectomorphism of $(M, \om)$ preserves the projective space $\Cal P$
then its infinitesimal part corresponds to the infinitesimal part of the corresponding automorphism
of $\Cal P$ which is a vector field. Really such vector field can be defined for
any Hamiltonian vector field (again, we are working with compact simply connected
symplectic manifolds!). Then one has a correspondence between
Hamiltonian vector fields on $(M, \om)$ and some special vector fields on $\Cal P$.
This specialty means that since any automorphism of $\Cal P$ induced by a symplectomorphism
of $(M, \om)$ by the definition keeps whole the Kaehler structure on $\Cal P$ then the induced
vector fields preserve the Kaehler structure as well. Therefore for any $f$
the dynamically correspondent vector field  $\tilde Y_f$ is Hamiltonian
for the symplectic structure on $\Cal P$ and is a Killing field  for the riemannian metric.
Since our projective space is simply connected one could reconstruct
a function which Hamiltonian vector field equals to $\tilde Y_f$. Denote this function as
$Q_f$. By the definition this $Q_f$ is a symbol (in the original Berezin sense).   
Of course, this $Q_f$ is defined only up to constant. But for any such function
we can construct a self adjoint operator on $\Cal H$. It's clear that a normalization
rule is included in (1.5.12') and the operator defined in this formula
coincides up to constant with the one coming from the symbol $Q_f$.
Therefore one should say that the good success of the Souriau - Kostant construction
is based on the dynamical property, mentioned above.

So at this stage the idea of dynamical correspondence comes to geometric quantization.
Let's consider the Souriau - Kostant construction as the following procedure.
For $(M, \om)$ as above one takes the prequantization data $(L, a)$. Then one takes the projectivization
$$
\Cal P = \proj(\Cal H)
$$
of the completed space of sections with finite norm. Then one has the dynamical correspondence
$$
\Theta^p_{DC}: Vect_{\om}(M) \to Vect_K(\Cal P)
\tag 1.5.18
$$
such that the image consists of the fields which preserve the Kaehler structure on $\Cal P$.
Dynamical properties of $\theta_{DC}$ dictates the following "correspondence
principle":
$$
\Theta^p_{DC}([X_f, X_g]) = [\Theta^p_{DC}(X_f), \Theta^p_{DC}(X_g)]
\tag 1.5.19
$$
So the quantization question is the following: can we lift this correspondence
to a correspondence between smooth functions on both the spaces?
More rigouresly we will say that the dynamical correspondence is quantizable iff
there exists a linear map
$$
F: C^{\infty}(M, \R) \to C^{\infty}(\Cal P, \R)
$$
such that the following diagram
$$
\matrix  C^{\infty}(M, \R) & @>F>> & C^{\infty}(\Cal P, \R)\\
\downarrow &  & \downarrow \\ Vect_{\om}(M) & @>\Theta^p_{DC}>> & Vect_K(\Cal P) \\
\endmatrix
\tag 1.5.20
$$
commutes. And one can prove that in  the Souriau - Kostant case the dynamical correspondence
given by (1.5.18) is quantizable. (Hint: if one takes unit sphere $S \subset \Cal H$
which consists of 
$$
S = \{ s \quad | \int_M <s, s> d \mu_L = r \}
$$
then for every real function $f \in C^{\infty}(M, \R)$ one has
$$
\tilde F_f (s) = \tau \int_M f \cdot <s, s> d \mu_L
\tag 1.5.21
$$
which is a real function on $S$. It's clear that this function is invariant with respect
to the phase rotations therefore one  gets the corresponding function $F_f$ on the projective
space $\Cal P$. This function depends on two real parameters --- $r$ and $\tau$.
Choosing  appropriate values one gets the desired function which satisfies (1.5.20)
and moreover maps the constant function $f \equiv 1$ to the constant function
$F_f \equiv 1$.)

The Souriau - Kostant method was recognized as universal but not quite successful since
it gives a reducible representation of the Poisson algebra $C^{\infty}(M, \R)$.
This well known fact (see, f.e., [13]) we are going to explain
in terms of geometric reformulation. So the question is: what "reducibility"
and "irreducibility" mean when we are speaking on the projectivization language?
Usually they mean that some invariant subspaces either exist or do not exist.
After the projectivization one either gets some special submanifold of the projective space or
doesn't get it. This submanifold should possess the following property:
for every smooth function $f \in C^{\infty}(M, \R)$ the  Hamiltonian vector field of the corresponding symbol preserves this submanifold. In other words it is parallel to
the submanifold. On the other hand, we know that the image of $q$ (see the definition
of AGQ on pp. 4 and 5 above) should give a Lie subalgebra in the Lie algebra of
vector fields over the quantum phase space therefore if the submanifold exists
it were that the vector fields form an integrable distribution.
F.e., using the function map, given by (1.5.21) one can easily find
the invariant submanifold for the Souriau - Kostant quantization. But even more easily
one could find that the resulting Hamiltonian vector fields form
an integrable distribution on the reduced phase space of the quantization.
This is an explanation why we imposed condition 4b) in the definition
of AGQ. Condition 4a) was introduced to get some complete picture
(to be honest, it was imposed since we can make it satisfied  in what follows).

More generally, let us formulate the dynamical principle as a universal one in symplectic
geometry. Namely, if we associate any object with a given symplectic manifold $(M, \om)$
in  invariant way (so if it depends only on $(M, \om)$ itself) then every symplectomorphism
of $(M, \om)$ generates an automorphism of this associated object. Moreover, the generated automorphism
should preserve the structures on the associated object if they were defined in invariant way.
Turning to the infinitesimal level it gives us a correspondence between Hamiltonian vector fields
and some special fields on the associated object such that the induced vector fields preserve
the structures. Moreover the identity (1.5.19) holds for the dynamical correspondence.
Then if the associated object is a symplectic manifold then the dynamical correspondence 
maps Hamiltonian vector fields on the source manifold to Hamiltonian vector fields on the target
manifold and one  can ask whether or not this correspondence is quantizable so does there exist
a map from the space of smooth functions on the given manifold to the space of smooth functions
on the associated manifold such that the diagram (1.5.20) commutes. 
If the map exists then one can represent the Poisson algebra of the given manifold
in terms of the Poisson algebra of the associated symplectic manifold. 

Now we leave dynamics for a while to continue the recalling of the relevant constructions.

\subhead 1.6. Real and complex polarizations
\endsubhead

In this subsection we discuss what happens in known quantization picture
in presence of an additional structure --- a polarization. Here we focus only
on the simplest examples when this polarization is  "pure" --- real or complex.  
Moreover we simplify the story avoiding the introduction of metaplectic
structure and metaplectic correction since in the rest of this text
we will not meet these notions (and unfortunately we are still not in position
to introduce this correction in ALG(a) - quantization).

A complex polarization is a fixed integrable complex structure $I$ over $(M, \om)$,
compatible with $\om$. Thus our $M_I$ is a Kaehler manifold. The choice
of $I$ defines a riemannian metric $g$ on $M$. Moreover, since
our symplectic form is integer then $M_I$ is an algebraic manifold
(more precisely --- an algebraic variety). The prequantization data
$(L, a)$ over $(M, \om)$ can be specified in this complex case
as follows. Since $F_a$ is proportional to $\om$ it has type $(1,1)$
with respect to the fixed complex structure. This means that $(L, a)$
is a holomorphic line bundle: the holomorphic structure is defined
by   operator $\overline{\partial}_a$. This operator defines 
the space of holomorphic sections $H^0(M_I, L)$
which is naturally included in $\Cal H$ (our $M$ is compact). The celebrated Kodaira
theorem ensures us that there exists a power of $L$ such that the corresponding linear system
defines an embedding of $M_I$ to a projective space. This means that for
some appropriate $k \in \Bbb N$ the space of holomorphic sections
$H^0(M_I, L^k)$ hasn't based components; in other words for any point $x \in M_I$
there exists a holomorphic section $s \in H^0(M_I, L^k)$ non vanishing at the point
$$
s(x) \neq 0.
$$
Moreover, the subspace $H_x \subset H^0(M_I, L^k)$ consists of  vanishing  at $x$ holomorphic sections of $L^k$ has codimension 1 for any $x$. Therefore the following embedding
is correctly defined
$$
\aligned
\psi: & M_I \hookrightarrow \proj H^0(M_I, L^k)^*, \\
\psi: & x \mapsto (\proj H_x)^* \in \proj H^0(M_I, L^k)^*\\
\endaligned
\tag 1.6.1
$$
which is  the embedding by complete linear system $|k D|$ where $D$ is the Poincare dual
divisor. 

Usually one considers the space of holomorphic sections as the wave function space. The desired self adjoint operators can be defined in different ways.

\subheading{Berezin - Rawnsley quantization} In [16] (it has subtitle "geometric
interpretation of Berezin quantization") one finds an approach to quantize this Kaehler case.
First of all, one restricts  to a subclass of smooth functions which are called "qunatizable":
smooth function $f$ is quantizable over $M_I$ iff its Hamiltonian vector field preserves
the complex structure. In our notations it means that $f$ is a quasi symbol over the Kaehler manifold. For such a function the corresponding Souriau - Kostant operator $\hat Q_f$
(see subsection 1.5) maps the holomorphic sections to the holomorphic sections thus
it is an automorphism of the space $H^0(M_I, L)$. The correspondence principle
is satisfied since we just restrict the full picture over $\Cal H$ to the corresponding
subspace. We need not to explain what is geometric reformulation of this method
because it is geometric itself: all the paper [16] is done in this style and we could not
add anything to it. A Berezin idea to consider
on the projectivized space the symbols instead of self adjoint operators is exactly the same
geometric formulation of quantum mechanics applying to geometric quantization.

\subheading{Berezin - Toeplitz quantization} In [4] one can find
another way defining self adjoint operators on the same Hilbert space.
Let us denote
$$
\Cal H_k = H^0(M_I, L^k)
$$
and call $k$ the level of quantization. For any level we have the orthogonal projection
$$
S_k: \Cal H \to \Cal H_k
\tag 1.6.3
$$
which is called Szoege projector. Then every smooth function $f$ gives a self adjoint operator
on any level
$$
\hat A_f: \Cal H_k \to H_k
$$
coming as the composition of the multiplication by $f$ and the corresponding Szoege projector:
$$
\hat A_f(s) = S_k(f \cdot s) \in \Cal H_k
\tag 1.6.4
$$
where $s$ belongs to $\Cal H_k$. It's clear that for every $f$ the corresponding Toeplitz operator
is self adjoint. The choice of the holomorphic section spaces as  the subspaces to project on  is an appropriate one: in [4] one proves that in this case the correspondence principle 
is satisfied asymptotically.

Apply the geometric reformulation scheme to the Berezin - Toeplitz  method over compact symplectic manifolds (see [28]). Instead of finite dimensional Hilbert space $\Cal H_k$ we take
the corresponding projective space $\proj_k$. Any Toepltiz operator $\hat A_f$ can be replaced by
its symbol $Q_f$. Therefore one gets a linear  map:
$$
A: C^{\infty}(M, \R) \to C^{\infty}_q(\proj_k, \R).
\tag 1.6.5
$$
The last space is finite dimensional thus this map $A$ has a huge kernel.
The question about its cokernel can be simplified by the introduction of the Fourier - Berezin transformation. Let $M_I$ as above be a Kaehler compact manifold, (L, a) be the prequantization data read as a holomorphic line bundle equipped with a hermitian structure. Realize as above
the projective space as the factorization of unit sphere $S_k \subset \Cal H_k$ by the gauge
transformations. Then one has a real non negative smooth function $u_k(x, p)$ on the direct
product $M_I \times \proj_k$ defined as follows.  If section $s \in S_k$ represents
point $p \in \proj_k$ then
$$
u_k(x, p) = <s(x), s(x)>_h \in \R.
\tag 1.6.6
$$
This expression again is preserved by the gauge transformation so we get 
a correctly defined smooth function on the direct product. On the other hand
this function  strongly depends on the hermitian structure choice made before.
\proclaim{Proposition 2.3} For any smooth function $f$ on $M$ the corresponding
Toeplitz symbol $A_f$ is given by the formula
$$
A_f(p) = \int_M f \cdot u(x, p) d \mu_L.
\tag 1.6.7
$$
\endproclaim

This integral operator with kernel $u_k(x, p)$ performs the Fourier - Berezin transform
(of level $k$) mentioned above. At the same time there is a distinguished element
in the space of classical observables over $M$. Smooth function
$$
\la_k(x) = \int_{\proj_k} u_k(x, p) d \mu,
\tag 1.6.8
$$
where $d \mu$ is the Liouville volume form on $\proj_k$, depends, of course,
on  the choice of the prequantization hermitian structure:
when we vary the hermitian structure then both the integrand and the volume form
change in the right hand side of (1.6.8). This function $\la_k$ is an analogy
of the Rawnsley $\theta$ - function (or $\varepsilon$ - function) (see [16]).
To solve the cokernel problem for the map (1.6.5) one can 
choose such  a prequantization  hermitian structure on $L$ that 
$\la_k$ would have the best shape (and the bestest case is when $\la_k$ is constant,
but it is possible only if $M$ is a projective space itself, see [16]).
We don't touch it here but an interesting unsolved problem is to
find these distinguished functions just in the simplest cases
(f.e., for complete intersections) for algebraic manifolds in the framework
of algebraic geometry.
As well in symplectic geometry one should ask the questions about
a subspace in $C^{\infty}(M, \R)$ consists of all possible
$\la$: how large this subspace is and what is the geometrical meaning
of this subspace (note that this  subspace is defined absolutely
canonically without any additional choices over $(M, \om)$ with integer
$\om$).

To prove the statement it's convenient to choose an orthonormal basis
in the Hilbert space $\Cal H_k$. Assume that an appropriate one
is fixed, denoting it as $\{s_0, ..., s_d\}$. Thus
$$
\int_M <s_i, s_j>_h d\mu_L = \de_{ij}.
\tag 1.6.9
$$
Then the Toeplitz operator $\hat A_f$ is represented  in this basis
by a matrix whose elements are given by the formula
$$
(\hat A_f)_{i, j} = <f s_i, s_j>_q = \int_M <f s_i, s_j>_h d \mu_L.
\tag 1.6.10
$$
A unit vector $s \in S_k \subset \Cal H_k$ can be decomposed
$$
s = \sum_i \al_i s_i.
$$
If this vector represents point $p \in \proj_k$ then under the $\hat A_f$ - action it
comes to
$$
\hat A_f s = \sum_j(\sum_i \al_i \int_M < f s_i, s_j>_h d \mu_L) s_j.
\tag 1.6.11
$$
The expectation value can be computed in the form
$$
<\hat A_f s, s>_q = \sum_j \al_j(\sum_i \al_i \int_M <f s_i, s_j>_h d \mu_L).
\tag 1.6.12
$$
It gives us
$$
<\hat A_f s, s>_q = \int_M f <s, s>_h d \mu_L.
\tag 1.6.13
$$
Taking in mind the definition of $u_k(x, p)$ given in (1.6.6) we understand formula (1.6.13)
 the same one as given in the statement.

One should remark that formula (1.6.7) looks very similar to formula (1.5.21).
Really it means that $F_f$ defined by (1.5.21) and $A_f$ defined by (1.6.7)
are the same function but the first one is defined on the ambient space
while the second one equals to the restriction of the first one
to a finite dimensional subspace. The question arises: why in the first case the correspondence principle is satisfied completely while in the second it is satisfied only asymptotically?
The point is that if a function $f$ on $M$ generates the Hamiltonian vector field
preserving the complex structure $I$ as well as given symplectic structure then
the Hamiltonian vector field $X_{F_f}$ on $\Cal P$ is tangent to
the projective subspace $\proj_k \subset \Cal P$. For such functions ("qunatizable"
in terms of Berezin - Rawnsley method) the results of two different methods are the same
(it was shown f.e. by G. Tuynman in [23]). At the same time if $f$ is arbitrary
then $X_{F_f}$ is not tangent to $\proj_k$ and using the corresponding projections
for two such functions $f$ and $g$ we get that $A_{\{f, g\}}$ differs from
$\{A_f, A_g \}$. But taking the limit $k \to \infty$ we make this difference
smaller and smaller and this means that in the Berezin - Toeplitz method
the correspondence principle is satisfied asymptotically.

On the other side of the spectrum one has real polarizations. While the complex case
as we've seen belongs to algebraic geometry and the basic examples are given by the projective
spaces, the real polarization  case has the cotangent bundles as the main examples.
It was proposed exactly as an additional structure on a symplectic manifold
which could mimicry the case when one has two different kinds of variables:
position and momentum ones. Therefore the best example of real polarized manifold
is given by cotangent bundle: if any symplectic manifold is globally symplectomorphic to a cotangent bundle then one can divide the variables in any local chart into two 
sets such that the first set can be called position variables while the second one
--- momentum variables. One can see that the method reminded below is specified
exactly to that case. At the same time we will see that this method has no choice
to be well defined in the compact case --- but anyway it can be retranslated
in this case as we'll show in Section 4.

Avoiding the discussion of Lagrangian distribution we directly say that one can understand
real polarization as the case when over $(M, \om)$ there are exist 
$n$ smooth real functions $f_i$ (where $2n$ is the dimension of $M$) in involution:
$$
\{f_i, f_j\} = 0 \quad \quad \forall \quad i, j = 1, ..., n,
\tag 1.6.14
$$
which define a Lagrangian fibration
$$
\pi: M \to \De \subset \R^n,
\tag 1.6.15
$$
where $\De$ is a convex polytope in $\R^n$, namely
$$
\pi^{-1}(t_1, ..., t_n) = \cap_i \{f_i = const = t_i\}.
\tag 1.6.16
$$
The Hamiltonian vector fields $X_{f_i}$ at any point of any fiber $\pi^{-1}(t)$    
form a complete basis of the tangent to the fiber space.
As usual one fixes a prequantization data $(L, a)$ and gets the corresponding big Hilbert space
$\Cal H$. Further, following the same policy as in the complex case, one takes only the sections
which satisfy some compatibility condition with respect to the fixed real polarization.
Namely these sections must be invariant under the flows generated by $\{f_i \}$. Again one introduces the notion of level $k$ quantization taking all powers of $L$. Thus the Hilbert space of wave functions consists
of the following sections:
$$
 \Cal H_1 = \{ s \in \Cal H \quad | \quad \quad \nabla_{X_{f_i}} s \equiv 0 \quad \forall 1 \leq i \leq n\}.
$$
So the wave function is constant along the fibers
and it varies only along the base (thus it depends on just a half of the variables). But the based manifold $\De$
wasn't equipped with any canonical measure or volume form since the basic Liouville form
could not give any help when we turn to $\De$. At the same time in  general non compact case
the fibers could be  noncompact and then any section constant along the fibers should have unbounded norm
with respect to the Liouville norm. To avoid this difficulty in the noncompact case
(note that everybody takes in mind a cotangent bundle working in this way) one introduces
the following correction
$$
L' = L \otimes \sqrt{\La^n F}
\tag 1.6.17
$$
where $F$ is the subbundle of $TM$  tangent to the fibers of (1.6.15) at every point. In the real case
we are talking about this additional term is topologically trivial moreover it carries a canonical
connection gauge equivalent to ordinary $d$. Thus from the topological point of view 
this additional term doesn't change a choice on the base which we'll mention in a moment.
But now this term which is known as the bundle on half weights allows one to
define some  natural pairing  on $\Cal H_1'$ which is defined in the same vein as $\Cal H_1$. For any pair of the sections one takes the hermitian product for $L$ - parts and just the tensor product for the half weight parts.
This give a weight on the base $\De$ which can be integrated, giving a number. Now let us discuss what is the relationship between these special sections from $\Cal H_1$  and some data inside the based polytope.
Really if a section $s$ is constant along a fiber $S = \pi^{-1}(t)$ then either it vanishes at $S$ or
this Lagrangian submanifold $S$ must satisfy some special condition. We will discuss this condition
(which is called here the Bohr - Sommerfeld condition) in Section 2 but here one can formulate it easily:
if we take the restriction of the prequantization data
$$
(L, a)|_S
$$
to our Lagrangian fiber $S$ then it should admit a covarinatly constant section. It's clear that
by the definition of the prequantization data this restriction must be special, consisting
of trivial line bundle plus a flat connection. And it is true for {\it any} Lagrangian submanifold.
But this flat connection could be in general non equivalent to the ordinary one if
our $S$ has non trivial fundamental group. Thus in this case the requirement for the restriction
to be trivial is called the Bohr - Sommerfeld condition. The confusion  point is that there is a number of 
conditions having the same name; f.e. the integrality condition for the symplectic form is often
called by the same name. Generally speaking any Bohr - Sommerfeld condition
is just an appropriate quantization condition so f.e. in our main compact case 
the integrality condition is not sufficient to perform the quantization. Therefore we have to
deepen and specify this condition to get a real Bohr - Sommerfeld condition for
the quantization.

Now, it's clear that if $S$ is Bohr - Sommerfeld with respect to $(L, a)$ then it remains to be
Bohr - Sommerfeld with respect to $(L', a')$. Let us collect all the fibers which satisfy
the Bohr - Sommerfeld condition getting a submanifold
$B \subset \De$ which is called again the Bohr - Sommerfeld submanifold. Outside of $\pi^{-1}(B) \subset M$
all sections from $\Cal H_1$ have to vanish. The structure of the wave function space is the following: for any
connected component $S_{\al}$ of $\pi^{-1}(B)$ there is a subspace $\Cal H_{\al}$ in $\Cal H_1$ which consists
of the sections whose support is exactly $S_{\al}$. Then the whole space is decomposed
$$
\Cal H_1 = \sum_{\al} \Cal H_{\al}.
\tag 1.6.18
$$
(We will show in Section 3 that if $M$ is compact as well as the fibers then the Bohr - Sommerfeld
manifold $B$ either consists of a finite number of isolated points on $\De$ or coincides with $\De$.  Therefore
in the first case the wave function space is spanned on a discrete set of $\de$ - sections.) 
Again as in the complex case one should distinguish some set of functions in $C^{\infty}(M, \R)$
which could be quantized. These functions are quite rare: in some extremal cases
only our based $f_i$ are quantizable. For any such $f_i$ the corresponding 
self adjoint operator is diagonal: it acts by multiplication
$$
\hat Q_{f_i}(s) = f_i \cdot s.
$$
Since $f_i$ is constant along the fibers it's clear that $f_i \cdot s$ remains to be covariantly constant along
the fibers so
$$
\hat Q_i: \Cal H_1 \to \Cal H_1
$$
is an automorphism of the quantum space. In the compact case when the sections from $\Cal H_1$ have
bounded usual Liouville norms one could perform the same procedure as in Berezin - Toeplitz method,
defining the corresponding orthogonal projector to subspace $\Cal H_1 \subset \Cal H$.
Then a priori for every function $f \in C^{\infty}(M, \R)$ there exists a direct analog to
the Toeplitz operators discussed above. At this point the main difference between $\Cal H_1$ and
$\Cal H_1'$ comes: the quantum space $\Cal H_1'$ doesn't belong to $\Cal H$ and moreover
can't be included there. It carries manifestly different hermitian structure. Thus one must use
some additional technique to define quantum operators for a  wider set of smooth functions. 
Recall briefly, what one could do here,  following [20]. Let $F_1, F_2$ be two real polarizations with
transversal fibers. Then on the direct product of the corresponding quantum Hilbert spaces
$\Cal H_1^1$ and $\Cal H_1^2$ there is a sesquilinear pairing
$$
K: \Cal H_1^1 \times \Cal H_1^2 \to \C 
\tag 1.6.19
$$
which is called the Blattner - Kostant - Strenberg kernel ([20]). The hermitian structure  on $L$ 
plays the role in the definition
of this pairing as well as the half weight parts. For reducible sections
$$
s_i' = s_i \otimes \theta_i \quad | \quad s_i \in \Cal H^i_1, \theta_i \in \Ga(\sqrt{\La^n F_i}),
$$  
where $i = 1, 2$, one defines 
$$
<s_1', s_2'>_{\theta} = \int_M <s_1, s_2>_h \cdot \theta_1^2 \wedge \theta_2^2 \in \C.
\tag 1.6.20
$$
Indeed, since $F_1, F_2$ are transversal the wedge product of two elements from $\La^n F_1$ and
$\La^n F_2$ is just a volume form on $M$. This sesquilinear pairing together with
given hermitian structures on $\Cal H^i_1$
defines a $\C$ - linear map
$$
U_{21}: \Cal H_1^2 \to \Cal H_1^1.
\tag 1.6.21
$$
Assume now that a smooth function $f$ (which doesn't belong to the algebraic span of the set $\{f_i \}$)
generates the flow which deforms the given real polarization to a polarization transversal to the given one.
Then the quantum operator is defined as the composition of the infinitesimal moving of $\Cal H_1$ in the direction
of $X_f$ with the map defined in (1.6.21). One can see that in compact case  there is no such a function
at all since the definition means that for the function its Hamiltonian vector field should be nonvanishing
everywhere (otherwise we get a point where these two real polarizations would be tangent). On the other hand
even in non compact case it's an open problem because this derived operator is not in general self dual
(the detailed story can be found f.e. in [20], [33]).

Thus in the compact case we are speaking about the subject in geometric formulation
since as we've seen the real geometric information is carried by the set of the Bohr - Sommerfeld
fibers.

\newpage

\head 2. The correspondence principle in the algebraic Lagrangian geometry
\endhead

$$
$$

In this Section we construct a special representation of the Poisson algebra
over a simply connected compact symplectic manifold with integer symplectic form.
During this procedure the Poisson algebra is represented in the space of
smooth functions over an appropriate moduli space which is itself a Kaehler manifold.
This moduli space was constructed in [10], [25] and called the moduli space of half weighted
Bohr - Sommerfeld Lagrangian cycles of fixed topological type and volume. To establish
that this is a faithful representation we prove following [28] that the corresponding map
is a homomorphism of the Poisson algebras. This representation could be called canonical
since it doesn't require any additional choices for the definition. At the beginning
we recall briefly the basic notions and constructions from ALAG.

$$
$$
$$
$$

\subhead 2.1. Bohr - Sommerfeld condition
\endsubhead

Let $(M, \om)$ be again a compact simply connected symplectic manifold of dimension $2n$ with integer symplectic form. Consider again the prequantization data $(L, a)$ where $L$ is defined 
by the condition
$$
c_1(L) = [\om]
$$
and $F_a = 2 \pi i \om$. For an appropriate smooth oriented connected $n$ - dimensional
manifold $S$ consider the space of smooth Lagrangian embeddings of fixed topological type
so such smooth maps
$$
\phi: S \to M
\tag 2.1.1
$$
that 
$$
\phi^* \om \equiv 0
\tag 2.1.2
$$
and the images represent the same homology class $[S] \in H_n(M, \Z)$.
The choice of the prequantization data make it possible to impose an extra condition
on the Lagrangian embeddings which has been mentioned in the previous section. So we call
an embedding Bohr - Sommerfeld if the restriction of the prequantization data
to the image admits covariant constant sections. In other words flat connection $\phi^*a$
on trivial bundle $\phi^*L$ should have trivial periods with respect to 
the fundamental group of $S$. The condition has more elegant version recovering
its geometrical meaning more distinct what will be discussed in Section 5. If one takes
the corresponding $U(1)$ - principle bundle with the corresponding connection
1 - form $A$  then it is an example of a contact manifold (about these manifolds see, f.e., [2]).
Connection $A$, retwisted by $i$, satisfies the standard   
condition
$$
\al \wedge (d \al)^n = d \mu,
$$
where $d \mu$ is a volume form on $P$. In terms of the principle bundle
one can easily reformulate the Bohr - Sommerfeld condition. Lagrangian submanifold satisfies
the Bohr - Sommerfeld condition iff it can be lifted to $P$ along the fibers
of the canonical projection $P \to M$. Our connection $A$ decomposes the tangent to $P$
space at every point as the direct sum of the horizontal and the vertical parts and a map
$$
\tilde \phi: S \to P
$$
is called planckian  if $T(\tilde \phi)$ belongs to the horizontal
part at every point of $\tilde \phi$ and $\tilde \phi^* \pi^* \om \equiv 0$
where $\pi$ is the standard projection.

Define now the notions of Bohr - Sommerfeld and planckian cycles. Let $\tilde{ \Cal B}_S$
is the space of all Bohr - Sommerfeld Lagrangian embeddings (2.1.1) of fixed topological
type. Then the moduli space of Bohr - Sommerfeld Lagrangian cycles is given by the factorization
$$
\Cal B_S = \tilde{\Cal B}_S / Diff_0 S,
\tag 2.1.3
$$
where $Diff_0 S$ is the identity component in the diffeomorphism group of $S$. Recall
that $S$ is oriented thus $Diff_0 S$ could be understood as the parameterization group of $S$.
Points of the moduli space are called Bohr - Sommerfeld cycles of fixed topological type.
The moduli space of planckian cycles has almost the same definition: one just starts
with the space of all planckian embeddings of $S$ to $M$ described above.
We denote it as $\Cal P_S$, following [10], [25]. Every planckian cycle is represented
by a covariant constant lifting of a Bohr - Sommerfeld cycle therefore the natural map
$$
\pi: \Cal P_S \to \Cal B_S
\tag 2.1.4
$$
gives a principle $U(1)$ - bundle structure on $\Cal P_S$ such that the canonical $U(1)$
- action is generated by the canonical $U(1)$ - action on $P$.  This principle bundle is called Barry bundle. 

It's quite natural to include the notion of level $k$ here: so if one takes
the corresponding tensor power $(L^k, a_k)$ then one could define the notions with respect
to this power as well as in the original case. This gives us a set of moduli spaces
parameterized by $k$:
$$
\aligned
\Cal B_S = \Cal B_{S, 1}, ..., \Cal B_{S, k}, ... \\
\Cal P_S = \Cal P_{S, 1}, ..., \Cal P_{S, k}, ... \\
\endaligned
$$
In Sections 2 and 3 we formulate the main facts in terms of level 1, but all of these statements
can be easily generalized to the case of any level. We just have to
keep in mind that if one starts with the symplectic manifold $(M, k \om)$
then the pair $(L^k, a_k)$ is precisely a prequantization data for it. And two major differences come: the Poisson bracket for $k \om$ is slightly different from the original one
and the Liouville volume form for $k \om$ is slightly different to. All of these remarks
will be  valuable, important and meaningful in Section 4 in the discussion of the quasi classical limit.

Our first aim is to describe smooth structures on the moduli spaces $\Cal B_S, \Cal P_S$.
The first one is given by the following
\proclaim{Proposition 2.1.([10], [25])} The tangent space $T_S \Cal B_S$ at any point
$S \in \Cal B_S$ is isomorphic to $C^{\infty}(M, \R)$ modulo constant functions.
\endproclaim    

The proof can be found in [10]. Briefly, let $S$ be a regular point
of the moduli space $\Cal B_S$. We identify it in our discussion
with the image of the corresponding class of maps so we understand $S$ as
an oriented smooth Bohr - Sommerfeld Lagrangian submanifold of $M$. According to
the Darboux - Weinstein theorem (see [32]) there exists a tubular neighborhood $N(S)$ such that
it is symplectomorphic to an $\varepsilon$ - neighborhood of the zero section of the cotangent bundle:
$$
\psi: N(S) \to N_{\varepsilon}(T^*S), 
\tag 2.1.5
$$
where the last one is equipped with the restriction of the canonical symplectic form
on the cotangent bundle. Thus the work is reduced to the canonical case. Recall
that $T^*S$ is endowed with a natural 1 - form $\eta$ which is called by the natural
reason the canonical 1 - form. Any 1 - form $\al$ on the based manifold
$S$ can be regarded from two points of view: as a 1 - form on $S$ itself and as a section
$\Ga_{\al}$ of the cotangent bundle:
$$
\Ga_{\al}: S \hookrightarrow T^*S \to S,
\tag 2.1.6
$$
where the last arrow is the canonical projection. Then the canonical 1 - form on $T^*S$
is correctly defined by the following property satisfied for any 1 - form $\al$:
$$
\Ga^*_{\al} \eta = \al.
\tag 2.1.7
$$
The differential of this form $d \eta$ is a non degenerated everywhere closed 2 - form
defining the canonical symplectic structure over $T^*S$. Therefore one could understand
submanifolds of $M$ which belong to the neighborhood $N(S)$ as submanifolds
of $T^*S$ sufficiently close to $S$. If we are talking about lagrangian submanifolds
then all of them are described by sufficiently "small" closed  1 - forms over $S$.
It's clear that section $\be: S \hookrightarrow T^*S$ is lagrangian iff
$$
\om|_{\Ga_{\be}} = 0 = d \eta|_{\Ga_{\be}} = d \be.
\tag 2.1.8
$$
The prequantization data for the cotangent bundle are given by 
$$
L = T^*S \times \C, \quad \quad a = d + 2 \pi i \eta,
\tag 2.1.9
$$
where $d$ is the ordinary differential. The restriction of this flat connection
on any Lagrangian graph $\Ga_{\be} \subset T^*S$ corresponding to $\be \in \Om^1_S$
has the periods along the group elements from $\pi_1(S)$ equal to the periods
of the given 1 - form $\be$ (it follows from the definition of $\eta$). The last periods
are trivial iff the form is exact. Then the Bohr - Sommerfeld condition is equivalent to
the exactness of the corresponding 1 - forms. And this gives us
$$
T_S \Cal B_S = B^1(S) = \{ df \}.
\tag 2.1.10.
$$
Turning to the principal bundle (2.1.4) it's easy to see that the lifting corresponds
to the employing of the constant functions so one gets
$$
T_{\tilde S} \Cal P_S = C^{\infty}(M, \R),
$$
where $\tilde S$ is the corresponding planckian cycle. Moreover, the Darboux - Weinstein
theorem ensures that the representations for the tangent spaces are integrable
so one has distinct local coordinate systems on both the moduli spaces.
The set of the Darboux - Weinstein neighborhoods gives  atlases of the smooth structures.
The fact stated in the theorem underlies of  all our constructions.
One could say even more: this fact shows that Lagrangian submanifolds of  symplectic
manifold look like points of symplectic manifold. Really it's well known
(from the original Darboux lemma) that locally all symplectic manifolds are the same
(in contrast, f.e., with riemannian manifolds which have local invariants distinguishing
each from  others). Thus all the points of a given symplectic manifold
are indistinguishable having the same local structure. And according
to the generalization of the old classical result the same is almost true for
Lagrangian submanifolds: they differ only by the topological type. Therefore
one sees that "points" and "Lagrangian submanifolds" have quite similar behavior
from the kinematical point of view. We'll see that the same should be said
about the dynamical behavior. 

More generally, if one defines the moduli space of all Lagrangian cycles
in the same fashion as it was done for Bohr - Sommerfeld cycles
(applying the Darboux - Weinstein theorem and getting again
an atlas for the smooth structure and "canonical" local coordinate)
then there is a natural map
$$
\chi: \Cal L_S \to J_S = H^1(S, \R) / H^1(S, \Z),
\tag 2.1.11
$$
where $\Cal L_S$ is the moduli space of all Lagrangian cycles of fixed topological type.
$J_S$ is the Jacobian of $S$ understood as the set of flat connection classes modulo gauge transformations. Obviously the preimage of zero class is exactly $\Cal B_S$.
On the other hand one can take any fiber $\chi^{-1}(pt)$ as a moduli space
and this space would have the same description of the tangent bundle
as in Proposition 2.1. The main difference is that one couldn't define 
for this generic fiber the notion of the planckian lifting. Thus one couldn't
extend the description from $C^{\infty}(S, \R)/const$ (as for $\Cal B_S$)
to $C^{\infty}(S, \R)$ (as for $\Cal P_S$). Therefore our choice
of the fiber is based on the possibility of the lifting. At the same time
as we'll see  one can get almost the same fibration as (2.1.11)
in absolutely general symplectic situation. 

The description of the tangent bundles for $\Cal B_S$ and $\Cal P_S$ 
in terms of the smooth functions on $S$ has quite important consequences.
Before  the story is coming to the Kaehler setup we present here four remarks.

1. If $S$ has trivial fundamental group (or even
trivial first homology group) then each Lagrangian cycle is Bohr - Sommerfeld.

2. The linearization of the Bohr - Sommerfeld condition is exactly the same
as the so - called isodractic (or Hamiltonian) deformations. Indeed, any smooth function
$f$ on $S$ can be extended to a smooth function $\tilde f$ over $M$. Then the Hamiltonian vector
field $X_{\tilde f}$ generates some deformation of $S$. This infinitesimal deformation
preserves the Lagrangian condition. Moreover it preserves the Bohr - Sommerfeld condition
and the linearly deformed cycle is exactly that one in the neighborhood of
$S$ given by $df$. It means that the Bohr - Sommerfeld condition is a classical dynamical
condition over symplectic manifolds. Therefore we can introduce a kind of
fibration (2.1.11) in the case when $\om$ is not integer at all. Namely this analogy
is given by the flows of all complete Hamiltonian vector fields over $M$. But 
in the integer case we can first of all avoid the questions about the completeness
of the Hamiltonian vector fields and define the fibration on the "kinematical"
level. The same remark can be addressed to the planckian cycles as well.

3. Let's show some formulas to illustrate how we'll work in the setup. If you choose 
any function $f$ on $S \in \Cal B_S$ and extend it arbitrary to $M$,
getting a smooth function $\tilde f$ then one can decompose the corresponding
Hamiltonian vector field $X_{\tilde f}$ on the horizontal and vertical components
with respect to $TS$ and $\om$ and this decomposition is absolutely canonical.
So we have
$$
X_{\tilde f} = X_{ver} + X_{hor},
\tag 2.1.12
$$
where $X_{hor}$ belongs to $TS$ while $X_{ver}$ can be identified
with a section of the normal to $S$ bundle
$$
N_S = TM|_S / TS.
\tag 2.1.13
$$
It's clear that $X_{hor}$ corresponds to the part of deformation which preserves
the cycle $S$ (its flow generates some motion on $S$). Thus the deformation of
$S$ depends only on $X_{ver}$. Let's use the isomorphism
$$
\om: TM \to T^*M,
\tag 2.1.14
$$
getting the formula
$$
X_{ver} = \om^{-1}(d \tilde f|_S) = \om^{-1}(d(\tilde f|_S)) = \om^{-1}(df).
\tag 2.1.15
$$
Therefore the deformation depends only on the restriction to $S$. 

4. As we've seen there is a natural infinitesimal action of the Hamiltonian vector fields
on the moduli space of Bohr - Sommerfeld Lagrangian cycles. Indeed, 
every Hamiltonian vector field gives an infinitesimal
deformation of the based manifold so generates a vector field on the moduli space of Bohr - Sommerfeld Lagrangian cycles (dynamical correspondence) since the definition
is posed in invariant terms. The same is true for any fiber of (2.1.11).
The point is that the dynamical correspondent vector field $A_f$ for any
(global) smooth function $f$ on whole
$\Cal L_S$ is given by the following very simple formula
$$
A_f(S) = d(f|_S) \in T_S \Cal L_S.
$$
Thus this "quantum" vector field preserves the fibers of (2.1.11). The corresponding
foliation is integrable and we get (2.1.11) in another fashion, which has been hinted.

\subhead 2.2. Doubling circuit: Kaehler structure
\endsubhead

In this subsection we complexify the moduli space $\Cal B_S$, following [10], [25].
At the first step we take the moduli space of planckian cycles $\Cal P_S$. 
The source manifold $S$ is equipped with a space of half weights
(see [10], [25]). Since $S$ is orientable the determinant line bundle
$$
det T^*S = \La^n T^*S
\tag 2.2.1
$$
is trivial. Roughly speaking, a half weight is almost the same 
as a half form without zeros (at list in our case when $S$ is endowed
with a fixed orientation we will understand it so). For any pair
of half weights there are two derivations:
$$
\int_S \theta_1 \cdot \theta_2 \in \R
\tag 2.2.2
$$
and
$$
\frac{\theta_1}{\theta_2} \in C^{\infty}(S, \R).
\tag 2.2.3
$$
Moreover, the space of half weights admits a canonical involution
which could be written in the half form representation just as 
the multiplication by $-1$. The tangent space to the set of  half weights
over each point is modeled by $C^{\infty}(S, \R)$ ([10], [25])
and we consider the moduli space of half weighted planckian cycles
([10], [25]) consists of pairs
$$
(\tilde S, \theta) \in \Cal P^{hw}_S,
\tag 2.2.4
$$
where $\tilde S$ is a planckian cycle and $\theta$ is a half weight on the first element
which one understands as the image of the corresponding half weight on the source manifold.
The volume of this pair is given
by
$$
\int_S \theta^2 \in \R.
\tag 2.2.5
$$
By the definition the moduli space of half weighted moduli space is fibered
over the old one
$$
\aligned
\pi_{un}: \Cal P^{hw}_S \to \Cal P_S, \\
\pi_{un}:(\tilde S, \theta) \mapsto \tilde S. \\
\endaligned
\tag 2.2.6
$$
Moreover, there is another natural fibration
$$
\pi_c: \Cal P^{hw}_S \to \Cal B_S
\tag 2.2.7
$$
equals to the composition of (2.1.4) and (2.2.6). Therefore the moduli space of half weighted
planckian cycles inherits a $U(1)$ - principle bundle structure coming from (2.1.4).
And it was already mentioned that $\Cal P^{hw}_S$ carries the canonical volume  function
$$
\aligned
\mu: \Cal P^{hw}_S \to \R, \\
\mu(\tilde S, \theta) = \int_{\tilde S} \theta^2\\
\endaligned
\tag 2.2.8
 $$
which is obviously invariant under the $U(1)$ - action. The last remark is going to be
very important after the following fact is established:
\proclaim{Proposition 2.3 ([10], [25])} The moduli space of half weighted planckian cycles
$\Cal P^{hw}_S$ admits a Kaehler structure invariant under the $U(1)$ - action.
\endproclaim
The idea of the proof is to exploit the specialty of the tangent spaces
to the moduli space $\Cal P^{hw}_S$. Over a point it is the direct sum
$$
T_{(\tilde S, \theta)} \Cal P^{hw}_S = C^{\infty}(S, \R) \oplus C^{\infty}(S, \R),
\tag 2.2.9
$$
and both the summands are identified canonically. Moreover, as one has canonical Darboux - Weinstein local coordinates for the planckian "unweighted" cycles as well 
there are canonical {\it complex} Darboux - Weinstein local coordinates for the 
moduli space of half weighted planckian cycles (see [10]).
These coordinates were introduced in [10]. Thus in an arbitrary point $(\tilde S_0, \theta_0)$
belongs to the moduli space  the canonical local coordinates are given by a pair
of real smooth functions 
$$
(\psi_1, \psi_2), \quad \quad \psi_i \in \C^{\infty}(S, \R),
$$
where the first function responds to the deformations of the planckian cycle
while the second one reflects the deformations of the half weight part. One can
easily express  in these coordinates two natural tensors "living" on the moduli space.
The first one has type (1,1) being a linear operator:
$$
I|_{(\tilde S_0, \theta_0)} (\psi_1, \psi_2) =  ( - \psi_2, \psi_1).
\tag 2.2.10
$$
The next one has type (2,0) being a skew symmetric 2 - form:
$$
\Om_{(\tilde S_0, \theta_0)}(v_1, v_2) = \int_{\tilde S_0} [\psi_1 \phi_2 - \psi_2 \phi_1] \theta_0^2,
\tag 2.2.11
$$
where $v_1 = (\psi_1, \psi_2), v_2 = (\phi_1, \phi_2)$ are tangent vectors. One could easily
check that this 2 - form is nondegenerated everywhere and that $\Om$ is compatible
with $I$. The corresponding riemannian metric has the form
$$
G_{(\tilde S_0, \theta_0)}(v_1, v_2) = \int_{\tilde S_0} [\phi_1 \psi_1 + \phi_2 \psi_2]\theta_0^2.
\tag 2.2.12
$$
One can check directly that the form is closed and that the complex structure is integrable.
But the authors use in [10], [25] the following elegant trick to establish the  results.
The point is that for any smooth manifold $X$ its cotangent bundle admits a canonical
symplectic structure while its tangent bundle admits a canonical complex structure. 
Thus it's sufficient to consider two natural maps (see [10]): the first one is a local isomorphism
 of $\Cal P^{hw}_S$ and $T \Cal P_S$ and the second one is a global double covering 
$\Cal P^{hw}_S \to T^* \Cal P_S$ without ramification. Both the map are defined
canonically so they do not require any additional choices than were made in the beginning
of the story. Moreover, as it was checked the structures (2.2.10) and (2.2.11)
are isomorphic to the canonical ones. Moreover, since the canonical symplectic
structure on $T^* \Cal P_S$ is "strong" then the same is true for (2.2.10).
Therefore $\Cal P_S^{hw}$ is an infinite dimensional Kaehler manifold 
and one should mention that the Kaehler structure was constructed canonically
without any additional choices. As well one can easily check that this Kaehler structure is
invariant under the action of $U(1)$ described above. The function $\mu$ defined in (2.2.8)
is a moment map for this action (see [10]) thus one can produce a new Kaehler manifold
using the standard mechanism of Kaehler reduction. To get this new manifold
we choose a regular  value of the moment map function
$$
\mu(\tilde S, \theta) = \int_{\tilde S} \theta^2 = r \in \R.
\tag 2.2.13
$$
The reduced Kaehler manifold is called the moduli space of Bohr - Sommerfeld
Lagrangian cycles of fixed volume and denoted as $\Cal B_s^{hw, r}$. Thus the real parameter
$r$ shows the volume of the weighted cycles. This moduli space is fibered over $\Cal B_S$
so as a symplectic manifold it admits a canonical real polarization.
At the same time it admits a canonical complex polarization being a Kaehler manifold.
Moreover, it is algebraic since the Kaehler metric is of the so - called Hodge type
(the Berry bundle is related with the Kaehler class, see [10]).
So this algebraic manifold is ready to be quantized. But the main topic of our
text is to show that this moduli space is the quantum phase space for a system
which is strongly related to the original classical mechanical system, described
by our given manifold $(M, \om)$. Having this fact as the target we shall omit all
the related discussion (which could be found in [10], [25]).

\subhead 2.3. Induced functions on the moduli space
\endsubhead

We will perform the computations using the following short local description
of the moduli space of Bohr - Sommerfeld Lagrangian cycles of fixed volume 1.
So it is an infinite dimensional Kaehler manifold which consists of pairs $(S, \theta)$
where $S$ is a Bohr - Sommerfeld Lagrangian cycle in $M$ and $\theta$ is a half weight on it such that for any pair $(S, \theta) \in \Cal B^{hw,1}_S$ the volume is fixed:
$$
\int_S \theta^2 = 1.
\tag 2.3.1
$$
The computation is done for the moduli space $\Cal B^{hw,1}_S$ but
it can be easily rearranged for any other positive volume $r$ or any higher level $k$.
The tangent space to the moduli space $\Cal B^{hw,1}_S$ in a point $(S, \theta)$
is represented by pairs $(\psi_1, \psi_2)$ such that
$$
\int_S \psi_i \theta^2 = 0, \quad \quad \psi_i \in C^{\infty}(S, \R).
\tag 2.3.2
$$
For any two tangent vectors $v_1 = (\psi_1, \psi_2), v_2 =
(\phi_1, \phi_2)$ at the point $(S, \theta)$ the symplectic form
$\Om$ has the form
$$
\int_S[\psi_1 \phi_2 - \psi_2 \phi_1] \theta^2.
\tag 2.3.3
$$
We've spoke about the dynamical correspondences which could be either
quantizable or not. But we'd like to keep the order
of the investigated facts so here in this section we will
not say anything about the dynamical correspondence and just directly introduce
some functions over the moduli space induced by smooth functions
from $C^{\infty}(M, \R)$, following [28]. For any $f \in C^{\infty}(M, \R)$
one has
$$
F_f \in C^{\infty}(\Cal B^{hw,1}_S, \R),
\tag 2.3.4
$$
defined in absolutely natural way. Indeed, at each point $(S, \theta)$ it is given by
$$
F_f (S, \theta) = \tau \int_S f|_S \theta^2 \in \R,
\tag 2.3.5
$$
where $\tau$ is a real parameter. Formula (2.3.5) gives a map
$$
\Cal F_{\tau}: C^{\infty}(M, \R) \to C^{\infty}(\Cal B^{hw,1}_S, \R)
\tag 2.3.6
$$
which is obviously linear. The main aim of the rest of the present section
 is to show that this linear map is a homomorphism of Lie algebras. The original symplectic structure $\om$ defines the Poisson bracket on the source space
while the constructed symplectic structure (2.3.3) defines
the quantum Poisson bracket on the target space in (2.3.6).
And as we'll see the map $\Cal F_{\tau}$ transforms the classical bracket to 
the quantum one up to a constant which depends on our real parameter
$\tau$. Here we follow [28] keeping the notations.

Let $f \in C^{\infty}(M, \R)$ is an arbitrary smooth function on our given
symplectic manifold $M$. Then its differential $df$ gives a tangent vector being restricted to any Bohr - Sommerfeld cycle $S$. This gives a vector field
$A_f$ on the moduli space $\Cal B^{hw,1}_S$ which doesn't depend on the second coordinate. Therefore this $A_f$ is constant along the fibers
$$
\Cal B^{hw,1}_S \to \Cal B_S.
\tag 2.3.7
$$
It's not hard to see that the singular points of $A_f$ are given by the 
pairs $(S, \theta)$ such that $S$ belongs to a level set of $f$.
As well our smooth function $f$ generates naturally a 1 - form on $\Cal B^{hw,1}_S$
which we denote as $B^f$. In point $(S_0, \theta_0)$ this 1 - form is given
by
$$
B^f_{(S_0, \theta_0)}(\psi_1, \psi_2) = \int_{S_0} f \cdot \psi_2  \theta_0^2.
\tag 2.3.8
$$
 The direct substitution to the formula (2.3.3) ensures that these objects are related,
namely:
$$
B^f = \Om^{-1}(A_f).
\tag 2.3.9
$$
Indeed, for any vector field 
$$
v = (\psi_1(S, \theta), \psi_2 (S, \theta))
$$
one has
$$
\Om_{(S, \theta)}(A_f, v) = \int_S[f \psi_2 - \psi_1 \cdot 0] \theta^2 = \int_S f \psi_2 
\theta^2 = B^f(v),
$$
since our vector field in the local coordinates looks as $(f, 0)$.
Now we are in position to state and prove 
\proclaim{Proposition 2.4.([28])} For any smooth functions $f, g \in C^{\infty}(M, \R)$
the following identity holds
$$
\{F_f, F_g\}_{\Om} = 2 \tau F_{\{f, g \}_{\om}}
\tag 2.3.10
$$
where $F_f, F_g$ are the images of $f, g$ under $\Cal F_{\tau}$.
\endproclaim

Let's remark first of all that the map $\Cal F_{\tau}$ doesn't preserve the standard
algebraic structure on $C^{\infty}(M, \R)$, defined by usual pointwise
multiplication. It follows from the same property of the integral:
the integration of the product $f \cdot g$ is not usually the same as
the product of two integrations of $f$ and $g$ respectively. It means that
$$
F_f \cdot F_g \neq F_{f\cdot g}.
$$
At the same time according to Proposition 2.4. the image 
$$
Im \Cal F_{\tau} \subset C^{\infty}(\Cal B^{hw,1}_S, \R)
$$
is a Lie subalgebra. Assume that the given classical mechanical system is  integrable
hence $(M, \om)$ admits a set of  $n$ smooth functions in involution which are
algebraically independent. Since
$$
\{f_i, f_j \}_{\om} = 0
$$
the induced functions $F_{f_1}, ..., F_{f_n}$ commute with respect to
the quantum Poisson bracket over the moduli space $\Cal B^{hw,1}_S$. But the same is true
for the set consists of the functions of the following shape
$$
F_{f_1^{r_1} \cdot ... \cdot f_n^{r_n}}, \quad \quad r_i \in \Z.
$$
The corresponding preimages, of course, lie in the algebraic span
 of $\{ f_i\}$. But according to our remark the last function
doesn't belong to the algebraic span of $F_{f_1}, ..., F_{f_n}$.
It means that for any integrable classical system the corresponding
moduli space (read: quantum system) is integrable to. Now a question
arises: if our given classical system was completely integrable
(so $dim M = 2n$) is the same true for the quantum system?
Roughly the space of commuting functions over $\Cal B^{hw,1}_S$
has dimension $\Z^n$ while the moduli space itself has
dimension $2 \cdot C^{\infty}(S, \R) - 2$ thus it seems that
in general the question is not quite obvious.

Now we start the proof, performing the computation for the case when $\tau = 1$ for simplicity since the general case follows immideately. For a function $f \in C^{\infty}(M, \R)$
we take the corresponding induced function $F_f$ on the moduli space and compute
its differential. Perturbing the arguments one gets in the local coordinates that
$$
dF_f(S_0, \theta_0)(\psi_1, \psi_2) =
\int_{S_0} 2 f \psi_2 \theta_0^2 + \int_{S_0} d\psi_1(\om^{-1}(df)|_{S_0}) \theta_0^2,
\tag 2.3.11
$$
where the first 1 - form on the right hand side is constant with respect to the first
coordinate while the second 1 - form doesn't depend on the second coordinate.
It's easy to recognize our 1 - from $B^f$ (2.3.8), multiplied by 2, under the first summand in (2.3.11). As it was already pointed out (see (2.3.9)) the symplectically dual
to $B^f$ vector field is represented exactly by $A_f$. Therefore 
the Hamiltonian vector field for our induced function $F_f$ has the following form
$$
X_{F_f} = \Om^{-1}(dF_f) = 2 A_f + C_f,
$$
where $C_f$ is the vector field symplectically dual to
the second summand in (2.3.11):
$$
(C_f)^{sd} = \int_{S_0} d \psi_1(\om^{-1}(df)|_{S_0}) \theta_0^2.
\tag 2.3.12
$$
We don't need to compute any explicit formula  for $C_f$ taking in mind the following argument. It was mentioned that vector field $A_f$ is constant along the second coordinate while the symplectically dual 1 - form $B^f$ maps to zero any vector field which depends
only on the first coordinate. The vector field $C_f$ and its symplectically dual
1 - form possess the inverse property. Our "quantum" symplectic form is divided
in the variables hence one has
$$
\aligned
\{F_f, F_g\}_{\Om} = \Om(X_{F_g}, X_{F_f}) = \Om(2 A_g + C_g, 2 A_f + C_f) = \\
= 2 \Om(A_g, C_f) - 2 \Om(A_f, C_g),\\
\endaligned
\tag 2.3.13
$$
since it follows from (2.3.3) that
$$
\Om(A_f, A_g) = \Om(C_f, C_g) = 0
$$
for any $f, g \in C^{\infty}(M, \R)$.  
Further, from (2.3.13) one gets
$$
\{F_f, F_g\}_{\Om} = 2 (C_g)^{sd}(A_f) - 2 (C_f)^{sd}(A_g),
\tag 2.3.14
$$
where $(C_f)^{sd}$ is given by (2.3.12). According to (2.3.14)
we can compute the bracket, avoiding to write down any expilcite formula for
$C_f$. Substituting the expressions for the vector fields and 1 - forms
in (2.3.14) we get
$$
\aligned
\{F_f, F_g\}_{\Om} = 2 \int_{S_0} df|_{S_0}(\om^{-1}(dg)|_{S_0})\theta_0^2 \\
- 2 \int_{S_0} dg|_{S_0} (\om^{-1}(df)|_{S_0}) \theta_0^2. \\
\endaligned
\tag 2.3.15
$$
We claim that the total integrand
$$
2[df|_{S_0}(\om^{-1}(dg)|_{S_0}) - dg|_{S_0}(\om^{-1}(df)|_{S_0})]
\tag 2.3.16
$$
is the restriction of a smooth function to $S_0$. We'll see in a moment that this smooth function is exactly the Poisson bracket of $f$ and $g$ multiplied by 2. To check
this coincidence we take the Poisson bracket $\{f, g \}_{\om}$. It's easy to see
that
$$
\aligned
2 \{f, g\}_{\om} & = 2 df(\om^{-1}(dg)) \\
= - 2 dg(\om^{-1}(df)) & = df(\om^{-1}(dg)) - dg(\om^{-1}(df)).\\
\endaligned
\tag 2.3.17
$$
For simplicity let's choose any compatible riemannian metric $g$ on $M$
getting the corresponding almost complex structure $I$. We are interested in
the induced by this choice local decomposition of the tangent bundle $TM$
over our Bohr - Sommerfeld Lagrangian cycle $S_0$. At each point $s \in S_0$ one has 
$$
T_s M = T_s S_0 \oplus I(T_s S_0).
\tag 2.3.18
$$
First of all,  expression  (2.3.16) differs from (2.3.17) by the restrictions
of the differentials and the vector fields to $S_0$. Let's manage the corresponding 
decompositions in (2.3.17). We denote as "hor" - components the
restrictions to the tangent directions to $S_0$ and as "vert" components all
what belong to the transversal spaces (thanks to the riemannian metric).
Rearrange the first summand in the right hand side of (2.3.17)
as follows
$$
\aligned
(df_{ver}((Ig^{-1}(dg))_{ver}) + df_{ver}((Ig^{-1}(dg))_{hor})+\\
df_{hor}((Ig^{-1}(dg))_{ver}) + df_{hor}((Ig^{-1}(dg))_{hor})). \\
\endaligned
\tag 2.3.19
$$
In (2.3.19) only two summands are non trivial --- namely the first
and the fourth ones ("ver - ver" and "hor - hor"). Analogiuosly,
one has to the second summand in (2.3.17) that it equals to
$$
-(dg_{ver}((Ig^{-1}(df))_{ver}) + dg_{hor}((Ig^{-1}(df))_{hor})),
\tag 2.3.20
$$
thus again we have only two summands. However from the compatibility condition
of $g, \om$ and $I$ we have that
$$
df_{ver}((Ig^{-1}(dg))_{ver}) = - dg_{hor}((Ig^{-1}(df))_{hor})
\tag 2.3.21
$$
and
$$
df_{hor}((Ig^{-1}(dg))_{hor}) = - dg_{ver}((Ig^{-1}(df))_{ver}).
\tag 2.3.22
$$
Therefore we can rewrite the expression for the Poisson bracket $\{f, g \}_{\om}$
restricted to $S_0$ in terms of "hor" - components only:
$$
\{f, g \}_{\om}|_{S_0} = df_{hor}((\om^{-1}(dg))_{hor}) - dg_{hor}((\om^{-1}(df))_{hor}).
\tag 2.3.23
$$
It remains to mention that "hor" - components are exactly what one gets
taking restrictions to $S_0$ thus
$$
\{f, g\}_{\om}|_{S_0} = df|_{S_0}(\om^{-1}(dg)|_{S_0}) -
dg|_{S_0}(\om^{-1}(df)|_{S_0})
\tag 2.3.24
$$
and the result doesn't depend on the choice of any compatible riemannian metric.
Comparing (2.3.16) and (2.3.24) we get the desired identity.

\subhead 2.4. Adjunct: integer and real parameters
\endsubhead

In this small subsection we'd like to mention how the identity, stated in Proposition
2.4, changes when we vary integer and real parameters, contained in the picture.
Recall, that there are two real continuous parameters $r$ and $\tau$ and one integer
parameter $k$. We start with level $k = 1$. 

\subheading{The first level} In this case one has that the Possion brackets are proportional with coefficient 
  $2 \tau$ (see (2.3.10). It's clear that this coefficient doesn't depend on the volume
of cycles. On the other hand, $\Cal F_{\tau}$ maps 
$$
f \equiv const = c \quad \quad \implies F_f \equiv const = \tau \cdot r \cdot c.
\tag 2.4.1
$$
Therefore if one wants to establish the situation when all numerical
quantization requirements from the Dirac list hold (it means that
$$
\aligned
\tau \cdot r = 1, \\
2 \tau =  1) \\
\endaligned
$$ 
it would be necessary
to take
$$
\aligned
\tau = \half, \\
r = 2.\\
\endaligned
$$
At this step we see that anyway the product $\tau \cdot r$ should equal to 1 while $2 \tau$ can vary
with respect to the question about the Planck constant. We understand the Planck constant just as a proportionality coefficient. 
From the Berezin point of view the Planck constant depends only on the level of quantization.
Let's see what happens when we change the last one.

\subheading{General level} One could perform the same constructions for any generic level. Let's fix
any $k \in \Bbb N$ and construct the moduli space $\Cal B^{hw,r}_{S, k}$ in the same manner as $\Cal B^{hw,r}_S$
starting with  slightly different prequantization data $(L^k, a_k)$. Then one has a natural embedding
$$
\Cal B^{hw,r}_S \to \Cal B^{hw,r}_{S, k}
\tag 2.4.2
$$
(see [10]). The Kaehler structures on both the moduli spaces in (2.3.26) are slightly different;
this means that f.e. the symplectic form $\Om$ on the level 1 moduli doesn't coincide with the restriction of the symplectic form $\Om_k$ defined on the level $k$: here it's a crucial point which comes from
the difference between canonical Darboux - Weinstein coordinates for $\Cal B^{hw,r}_S$ and
$\Cal B^{hw,r}_{S, k}$. Indeed, if $(S, \theta)$ is originally Bohr - Sommerfeld Lagrangian cycle with respect
to $(L, a)$ then it's clear that it remains to be Bohr - Sommerfeld with respect to
$(L^k, a_k)$. But the canonical complex Darboux - Weinstein coordinates for
the level 1 moduli space are given by $\om^{-1}(df) \oplus df$ where $f$ lives on
$S$ while for the level $k$ moduli space they are given by $(k \om)^{-1}(df) \oplus df = k^{-1} \om^{-1}(df) \oplus df$.
This means that one rescales one half of the first coordinate system to get the second. Locally 
the difference can be recognized as follows. Let $T^*S$ is as usual the tangent bundle
to $S$. Then this tangent bundle admits not only one symplectic structure (the canonical one)
but a family of symplectic structures with respect to real parameter $\la$. Indeed, in the
basic definition of the canonical  1 - form (see (2.1.7)) one could slightly change the right hand side
$$
\forall \al \in \Om^1_S \quad \quad \Ga^*_{\al} \eta_{\la} = \la \al.
\tag 2.4.3
$$
Thus for any $\la \in \R$ one gets an "almost" canonical 1 - form $\eta_{\la}$ which is nondegenerated
and gives a nondegenerated 2 - form $\om_{\la} = d \eta_{\la}$. This "almost" canonical symplectic form
looks like the canonical one; it's not hard to find an appropriate symplectomorphism
$$
\Psi_{\la}: (T^*S, d \eta) \to (T^*S, d \eta_{\la}),
\tag 2.4.4
$$
it just multiplies every fiber of the canonical projection by $\la$.
The family $d \eta_{\la}$ is a possible degeneration of the canonical
symplectic structure on $T^*S$. At the same time an interesting effect
appears: due to the canonical form of symplectic structure in the canonical coordinates
$d \eta$ and $d \eta_{\la}$ are proportional while the same is not true
for the canonical coordinates. The proportionality coefficient is just $\la$.
Turning to the canonical Poisson brackets one sees that the corresponding
skew symmetrical pairings on the function space are proportional to;
the ratio is $\la^{-1}$. Coming back to the moduli spaces we get that
$\Om$ and $\Om_k|_{\Cal B^{hw,r}_S}$ are proportional with coefficient
$k$. Therefore   if one defines induced functions on $\Cal B^{hw,r}_{S, k}$
constructing the similiur map
$$
\Cal F^k_{\tau}: C^{\infty}(M, \R) \to C^{\infty}(\Cal B^{hw,r}_{S, k}, \R),
\tag 2.4.3
$$
given by the same formula
$$
F_f (S, \theta) = \tau \int_S f|_S \theta^2 \in \R
\tag 2.4.4
$$
then 
$$
\{F_f, F_g\}_{\Om_k} |_{\Cal B^{hw,r}_S} = \frac{1}{k} \{F_f, F_g\}_{\Om}.
\tag 2.4.5
$$
But it doesn't lead to a contradiction since
$$
\{f, g \}_{k \om} = \frac{1}{k} \{f, g\}_{\om}.
$$
And totally it gives
$$
F_{\{f, g\}_{\om}} =  \frac{k}{2\tau} \{F_f, F_g\}_{\Om_k}
\tag 2.4.6
$$
over $\Cal B^{hw,r}_{S, k}$. Now if we want to satisfy the Dirac conditions 
one needs
$$
\aligned
\frac{2 \tau}{k} = \hbar, \\
\tau \cdot r = 1.\\
\endaligned
\tag 2.4.7         
$$
Of course we could say that $2 \tau$ has to equal 2 and $r$ equals $\half$ as above.
On the other hand we can fix the ratio 
$$
\frac{2 \tau}{k} = const,
$$
then it implies that
$$
\tau \to \infty \quad \quad \implies r \to 0
$$
and in the limit one gets the moduli space of unweighted Lagrangian cycles.
Really, if $k$ tends to $\infty$ then the moduli space of unweighted Bohr - Sommerfeld
cycles covers the moduli space of all Lagrangian cycles as a dense set
(as rational points cover $b_1(M)$ - torus which is the Jacobian). At the same time
the weights are going to zero (since $r \to 0$) so in the limit one could forget
about the second components in the pairs $(S, \theta)$. Anyway one couldn't define,
say, a Poisson structure on the moduli space of Lagrangian cycles as the limit
of the symplectic structures on the moduli spaces of different levels since as we've
seen the symplectic structure on the moduli space $\Cal B^{hw,r}_{S, k}$ degenerates
when $k \to \infty$. 

We'll come back to the discussion in Section 4.

\newpage

\head 3. Dynamical correspondence in the algebraic Lagrangian geometry
\endhead

$$
$$

In this section we essentially develop the result of the previous one. 
We show that the functions given by map (2.3.6) on the moduli space 
are quasi symbols for any $f \in C^{\infty}(M, \R)$. To prove this main fact we use
a dynamical correspondence showing that this one is quantizable in the sense
of subsection 1.5, (1.5.20). Here we follow [29], [30].

$$
$$
$$
$$

\subhead 3.1. Quasi symbols over Kaehler manifolds
\endsubhead

Let $\Cal K$ is a Kaehler manifold. It could be considered from two different points of view:
in the setup of complex geometry it is a complex manifold equipped with
a positive polarization while in the framework of symplectic geometry it is a symplectic
manifold equipped with a complex polarization.  If we say that it is a real manifold equipped
with a Kaehler triple it were a boundary (or  the most complete) setup. Anyway we understand
a Kaehler manifold in this way: one has a Kaehler triple $(G, J, \Om)$ over the based real manifold
and each element in the triple is of the same importance. Here $G$ is a riemannian metric,
$J$ is a complex structure and $\Om$ is a symplectic form such that the usual compatibility conditions
hold. One could recover each element of any triple from the other two. The elements in
Kaehler triples play their own individual roles: the riemannian metric $G$ responds to some
measurement process (distances, volumes, etc.), the symplectic form generates some dynamical properties
while the complex structure imposes a different kind of geometry (holomorphic vector bundles,
complex submanifolds, etc.) growing up over the real one. Moreover, one could consider the Kaehler triples
over any real manifold as solutions of some natural $Diff$ - invariant equations on the space
of all possible hermitian structures over the base (see [26]). (As well one has to mention here
that symplectic geometry itself admits  generalizations of classical results from the algebro geometrical setup, f.e., S. Donaldson extended the famous Kodaira result to the symplectic case, see [7],
and we will exploit this generalization below).

We remind here the following definition from the previous part of the text:
\proclaim{Definition} Smooth real function over Kaehler manifold $\Cal K$ with Kaehler structure
$(G, J, \Om)$ is called quasi symbol iff its Hamiltonian vector field preserves
the riemannian metric:
$$
Lie_{X_f} G \equiv  0.
\tag 3.1.1
$$
\endproclaim
Notice that the constant functions satisfy the requirement above so we will mention that they are
quasi symbols. For any Kaehler manifold $\Cal K$ we denote the space of quasi symbols
as $C^{\infty}_q(\Cal K, \R)$. It's obvious (see Proposition 1.3.) that $C^{\infty}_q (\Cal K, \R)$ is a Lie subalgebra
of the Poisson algebra over $\Cal K$. The dimensions of this subalgebra 
(as a vector space and as a finitely generated subalgebra) are two integer characteristics of the Kaehler manifold. It's not hard to roughly estimate the first number in the framework of complex or algebraic geometry.
First of all it's easy to see that the low boundary always equals to 1. 
As well we have
\proclaim{Proposition 3.1} For any compact Kaehler manifold $dim C^{\infty}_q(\Cal K, \R)$ is finite.
\endproclaim

To prove it we just remark that if Hamiltonian vector field $X_f$ preserves the riemannian metric
then it preserves as well the complex structure. It means that $X_f$ is the real part of
a holomorphic vector field. Thus one has a map
$$
\si: C^{\infty}_q(\Cal K, \R) \to H^0(M_J, T^{1,0}M_J).
\tag 3.1.2
$$
This map is $\R$ - linear. It's clear that its kernel contains only constant functions. Indeed, every holomorphic vector field can be reconstructed from its real part uniquely. Since real holomorphic
vector field doesn't exist the map is an inclusion modulo constant.  On the other hand, if a holomorphic
vector field $v_h \in H^0(M_J, T^{1,0})$ has real part $Re v_h$
equals to a Hamiltonian vector field 
$$
X_f = Re v_h
$$
then $Im v_h$ is not Hamiltonian. Really, if $Re v_f$ is Hamiltonian then $Im v_f$
is proportional to the gradient vector field for $f$ with respect to $G$. But a gradient vector
field never coincides with a Hamiltonian one. Therefore one complex dimension in
$H^0(M_J, T)$ can give at least one real dimension in $C^{\infty}(\Cal K, \R)$. Thus 
$$
1 \leq dim C^{\infty} (\Cal K, \R) \leq h^0(M_J, T) + 1.
\tag 3.1.3
$$
Since $\Cal K$ is compact the right hand side in (3.1.3)
is finite (and can be estimated by, say, the Riemann - Roch formula) so the same is true for the "middle hand side".

Of course, the bounds (3.1.3) are too rough. Even in the case when $\Cal K$ as a complex manifold
admits a lot of infinitesimal automorphisms these holomorphic vector fields need not
possess the property that their real parts are represented by some Hamiltonian functions
with respect to some symplectic structure. From the point of view of holomorphic geometry
the question is quite meaningful and useful: if $X$ is a complex manifold
what is the Kaehler metric on it which admits the maximal symmetry? This classical question from
the complex analysis has as we've seen a projection to the problem of quantization.
On the other hand, in the set of Kaehler manifolds one could distinguish a subset consists of
such submanifolds for which the question can be reformulated. We mean the Kaehler metrics of
Hodge type. For a metric of Hodge type the corresponding Kaehler class is integer and
we turn to the framework of algebraic geometry. 
Here we meet first of all the most "algebraic" manifold --- the projective space. This algebraic
manifold admits holomorphic vector fields and moreover every complex dimension in the space of these fields
can be realized as a complexified real dimension in the space of quasisymbols with respect to
the Fubini - Study metric. This metric is the result of the projectivization procedure
(discussed in Section 1) starting from a Hilbert space. Every holomorphic vector field
on $\C \proj^n$ comes from a linear operator on $\C^{n+1}$. If a hermitian structure on $\C^{n+1}$
is fixed then the space of linear operators splits into the direct sum of real and imaginary parts
(self adjoint and skew symmetric operators). Every real part (= self adjoint operator)
gives a real function on the projective space which is  a quasi symbol. Thus the real part
of any holomorphic vector field is just the Hamiltonian vector field of
the quasi symbol   given by the real part of the linear operator which induced
our holomorphic vector field. In this case one has
$$
(n+1)^2
$$
quasi symbols (real dimension, of course) while $(n+1)^2 - 1$ holomorphic vector fields
(complex dimension) --- we count the identical linear operator in the first set
while we don't count it in the second one since it gives zero vector field.

Thus the projective spaces admit some maximality property with respect to
the question how many quantum symmetries one Kaehler manifold carries.
This fact is included in complex analysis in the following way:
for hermitian metrics there is a notion of sectional holomorphic
curvature (and it is discussed in [3] from the point of view,  appropriate for us) and the projective spaces possess the property
that the Fubini - Study metrics are of the constant sectional holomorphic curvature.
Thus the condition of the constant sectional holomorphic curvature
can be retranslated to the language of quantum theory:
a Kaehler manifold $\Cal K$ is of the constant sectional holomorphic curvature
if this quantum space carries maximal quantum symmetries. It's well known fact
that in finite dimension there is unique such a Kaehler manifold ---
the projective space. In the infinite dimensional case it is an open
question in complex analysis. Therefore if one constructs an infinite dimensional
Kaehler manifold which admits maximal quantum symmetries (and which is not
a projective space) then
it were a candidate to the solution of this classical problem in complex differential geometry.

Now, if we have any Kaehler manifold with integer Kaehler class then we can 
emdedd it to a projective space (by the complete linear system of some power of
the Kaehler class).  So the question can be reduced to the case when $Q$
is a submanifold in $\C \proj^n$. Of course, it's just a first step in the study
of the question but we should skeep this discussion leaving some additional details
(available, f.e., from [16]) outside of this text to keep our major theme.

Notice, that the arguments, applied above to prove Proposition 3.1,
give us an additional statement which will be useful in subsection 3.4. Hence
we place here
\proclaim{Lemma} If $f \in C^{\infty}(\Cal K, \R)$ is a quasi symbol
over some Kaehler manifold $\Cal K$ then the critical set 
$$
Crit(f) = \{x \in \Cal K \quad | \quad df(x) = 0 \}
$$
is a complex submanifold of $\Cal K$.
\endproclaim

Really, for any $f \in C^{\infty}_q(\Cal K, \R)$ it should be 
 a holomorphic vector field $\si(f)$, given by (3.1.2), such that by the definition
$$
Crit (f) \equiv (\si(f))_0 \subset \Cal K,
$$
where $(\si(f))_0$ is the vanishing set of $\si(f)$. This immediately gives us the statement
since $(\si(f))_0$ has to be a complex submanifold.

Let us emphasize that the last statement is true in the infinite dimensional case
as well as in the previous one when we estimate the dimension. But we have to go
to the infinite dimensional case since it's clear that the finite dimensional is
exhausted.   
In the previous section
it was constructed a map from the space of smooth functions over $M$ to
the space of smooth functions over the moduli space of half weighted
Bohr - Sommerfeld cycles of fixed volume. But as we've seen in Section 1 one is interested
only on smooth functions of special type in the quantization framework; quantum observables
have to possess the property which requires for their Hamiltonian vector field to
preserve all the kinematical data on the quantum phase space. Thus to continue the story
one has to show that all these induced smooth functions are quasi symbols over the moduli space.
We formulate the main result of this section in the following
\proclaim{Theorem 1 ([30])} Let $(M, \om)$ be a simply connected compact symplectic manifold
with integer symplectic class, $[S]$ be a homological class of middle dimension and
$\Cal B^{hw,r}_S$ be the corresponding moduli space of half weighted Bohr - Sommerfeld
Lagrangian cycles of fixed volume equipped with the corresponding Kaehler triple $(G, I, \Om)$.
Then for the linear  map $\Cal F_{\tau}$ defined in (2.3.5) one has that

1) for any $f \in C^{\infty}(M, \R)$ the induced function $F_f$ is a quasisymbol
over the moduli space;

2) the correspondence principle takes place in the form
$$
\{F_f, F_g\}_{\Om} = 2 \tau F_{\{f, g \}_{\om}};
$$

3) the map 
$$
\Cal F_{\tau}: C^{\infty}(M, \R) \to C^{\infty}_q(\Cal B^{hw,r}_S, \R)
$$ 
gives an irreducible representation of the Poisson algebra.
\endproclaim

The second item from the list is already known from the previous section but during
the construction which is given by a dynamical correspondence we'll get
the identity "charge-free".

\subhead 3.2. Dynamical correspondence
\endsubhead

We've discussed dynamical correspondences in subsection 1.5. The root idea is very natural in
the symplectic setup. The point is that symplectic geometry was created (and understood)
as a strong and convenient mathematical language describing Hamiltonian mechanical systems.
So dynamics should be imposed to any investigation in the symplectic setup. At the same time one should
emphasize that the Lagrangian condition 
$$
\om|_S \equiv 0
$$
looks like local and static  while the Bohr - Sommerfeld condition
is dynamical: local Bohr - Sommerfeld deformations (see subsection 2.1) precisely correspond to
Hamiltonian deformations induced by the Hamiltonian dynamics. Thus the dynamical property of the system
defines a correspondence between Hamiltonian vector fields on the based manifold
and some special vector fields on the moduli space. Let us discuss this point more concretely.

Let $(M, \om)$ be as usual a simply connected compact symplectic manifold of dimension $2n$ with integer symplectic form.
Then for a fixed homology class $[S] \in H_n(M, \Z)$ one has an infinite dimensional moduli space
$\Cal B^{hw,r}_S$ of half weighted Bohr - Sommerfeld cycles of fixed volume, endowed with
a canonical Kaehler structure, described in Section 2. The moduli space consists of pairs $(S, \theta)$ where
the first element is a Bohr - Sommerfeld Lagrangian cycle and $\theta$ is a half weight over it such that
$$
\int_S \theta^2 = 1.
$$
Of course, this moduli space could be empty for some classes from $H_n(M, \Z)$ and some source manifold $S$,
but we reject these vanishing cases (f.e. every symplectic manifold carries the moduli space corresponds to
the trivial homology class in $H_n (M, \Z)$ and an appropriate topological type of $S$). Since our $M$ is simply connected 
we identify the space of Hamiltonian vector fields with the smooth function space modulo constant.
One has a map
$$
\Theta_{DC}: Vect_{\om}(M) \equiv C^{\infty}(M, \R)/const \to Vect (\Cal B^{hw,r}_S),
\tag 3.2.1
$$
which is called a dynamical correspondence. It can be constructed as follows. For any function $f$
one has the corresponding dynamics, induced by the Hamiltonian vector field on $M$.
This dynamics preserves our symplectic manifold $(M, \om)$ and moreover if we choose a prequantization
data it could be lifted to this setup almost uniquely (up to the canonical gauge transformations). Thus this dynamics
preserves whole the data hence it defines a germ of automorphism of the moduli space $\Cal B^{hw,r}_S$. This just gives
us a vector field on the moduli space which admits some additional properties. The specialty reflects the fact that
this vector field preserve the Kaehler structure since it was defined in invariant way. Since $M$ is compact
every function $f$ defines a germ of symplectomorphism so every Hamiltonian vector field gives an infinitesimal
transformation of the moduli space. Generalizing over the space of all Hamiltonian vector fields we get
the map, which is obviously linear. Moreover, by the construction one gets
\proclaim{Proposition 3.2([30])} The image
$$
Im \Theta_{DC} \subset Vect_K(\Cal B^{hw,1}_S),
$$
where the last space consists of the fields which preserve the Kaehler structure on 
the moduli space.
\endproclaim

To make the story even more concrete we should  compute the coordinates of any such special field which is generated by
some smooth function $f \in C^{\infty}(M, \R)$. This function first of all gives us
the corresponding vector field $X_f$ over $M$. Let $(S, \theta) \in \Cal B^{hw,r}_S$ be 
a half weighted Bohr - Sommerfeld cycle. Over the  support $S$ the vector field $X_f$ can be decomposed
into the inner and the outer parts:
$$
X_f = X_{ex} + X_{in},
\tag 3.2.2
$$
where
$$
X_{in} \in TS
$$
is the tangent component. Of course, we've met this decomposition in subsection 2.3
where we denoted it as "ver - hor". However here we'd like to change the notations
emphasizing that decomposition (3.2.2) doesn't depend on any additional choices.
Indeed, it's a remarkable property of symplectic geometry that over any point of a Lagrangian submanifold
the tangent space canonically splits with respect to the symplectic form. Namely
over $p \in S$ one has $T_pS$ and $\om^{-1}(T^*_pS)$ as two direct summands
decomposing $T_pM$:
$$
T_pM = T_pS \oplus \om^{-1}(T_p^*S).
$$
And then $X_{in} \in TS$ while $X_{ex}$ belongs to $\om^{-1}(T^*S)$. Thus $X_{ex}$ gives a deformation
of the cycle $S$ itself while the inner part $X_{in}$ deforms $\theta$ (and  preserves $S$ being tangent to it).
Notice that if we reject the half weight parts then it leads to some degeneration of the picture:
in this case the deformation is defined only by the restriction of the source function to $S$ but
in the half weighted case two functions with the same restriction to $S$ give different deformations of
the pair $(S, \theta)$ if they differs in a small neighborhood of $S$.  

Therefore we can express the corresponding deformation given by $X_f$ in the canonical complex Darboux - Weinstein coordinates (see subsection 2.2) as follows:
\proclaim{Proposition 3.3 ([30])} The dynamically correspondent vector field, induced by $X_f$, has 
at a point $(S, \theta)$ coordinates $(\psi_1, \psi_2)$ such that
$$
\aligned
\psi_1 = f|_S - \int_S f|S \theta^2, \\
\psi_2 = \frac{Lie_{X_{in}}\theta}{\theta}.\\
\endaligned
\tag 3.2.3
$$
\endproclaim
The second equality in (3.2.3) is moreless evident, so let's ensure that the first one takes place.
The external part of the Hamiltonian vector field can be identified with a section of the normal 
bundle $N_S = TM|_S /TS$. On the other hand, we identified in subsection 2.1 the deformations of $S$ with
the function space over it modulo constant. Our symplectic form $\om$ can be regarded as an isomorphism
$$
\om: TM \to T^*M.
\tag 3.2.4
$$
Over $S$ both the spaces are decomposed in the direct sums
$$
TM|_S = TS \oplus V, \quad \quad T^*M|_S = T*S \oplus V^*.
\tag 3.2.5
$$
The map
$$
\om: V \to T^*S 
\tag 3.2.5
$$
(which has been used a moment ago) is an isomorphism, giving
$$
\om(X_{ex}) = df|_S = d(f|_S).
\tag 3.2.7
$$
The last equality gives the first equality in Proposition 3.3.
Thus we see that the dynamically correspondent to $X_f$ vector field
over the moduli space has very natural expression in terms of the canonical
"dynamical" coordinates. Dynamical properties of the given system give us another
condition satisfied by the vector fields from $Im \Theta_{DC}$ namely
\proclaim{Proposition 3.4 ([30])} For any pair of smooth functions
$f, g$ the following identity holds
$$
\Theta_{DC}([X_f, X_g]) = [\Theta_{DC}(X_f), \Theta_{DC}(X_g)],
$$
where in the right hand side one takes the standard commutator over the moduli space.
\endproclaim

Now the natural question arises: does the dynamical correspondence $\Theta_{DC}$
described above  admit a lift to the function level? In subsection 1.5
we called this as qunatability property hence in that terms
one asks whether or not $\Theta_{DC}$ is quantizable. We find the answer to the question
comparing for any smooth function $f \in C^{\infty}(M, \R)$ two vector fields over the moduli space, namely $\Theta_{DC}(X_f)$ and the Hamiltonian vector field for $F_f$, 
derived from (2.3.6). This matching gives us the following
\proclaim{Proposition 3.5 ([30])} For any smooth function $f \in C^{\infty}(M, \R)$
one has
$$
X_{F_f} = 2 \tau \Theta_{DC}(X_f).
$$
\endproclaim

The proof of this key statement is contained in the next subsection (we check this
relationship just by direct computations). At the rest of the present one we explain
how the statements of Theorem 1 follow from the propositions listed above.

The first statement of the theorem follows from the definition
of quasi symbols and Propositions 3.2 and 3.5: for any
function $f$ the induced function $F_f$ generates the Hamiltonian vector field
proportional to a vector field which belongs to $Im \Theta_{DC}$. It means that
$X_{F_f}$ preserves whole the Kaehler structure over the moduli space hence
$F_f$ is a quasi symbol.

The second item follows from Propositions 3.4 and 3.5: indeed,  one continues
the equality of Proposition 3.4 in both directions, substituting the equality of Proposition
3.5
$$
\aligned
\frac{1}{2\tau} X_{F_{\{g, f\}_{\om}}} = \Theta_{DC}([X_f, X_g]) =\\
[\Theta_{DC}(X_f), \Theta_{DC}(X_g)] = \frac{1}{4 \tau^2}[X_{F_f}, X_{F_g}]\\
\endaligned
$$
and the last term has the standard representation as a Hamiltonian vector field.
This gives the correspondence principle in the familiar form
$$
\{F_f, F_g\}_{\Om} = 2 \tau F_{\{f,g\}_{\om}}.
$$

The third statement is very important for us.  As it was explained in subsection 1.5 the irreducibility
 condition in our non linear algebro geometrical setup contains two items:
the first one says that $\Cal F_{\tau}$ has trivial kernel and the second says that
there is no any smooth proper submanifold in the moduli space such that for every
$f$ the corresponding Hamiltonian vector field
$X_{F_f}$ is tangent to the submanifold. We begin with the second requirement and we check it
even in much more stronger form. Namely we show that for every point $(S, \theta) \in \Cal B^{hw,r}_S$
and every tangent vector $v = (\psi_1, \psi_2) \in T_{(S, \theta)} \Cal B^{hw, r}_S$
there exists a smooth real function $f \in C^{\infty}(M, \R)$ such that
the Hamiltonian vector field of the induced function $F_f$ gives this vector at this point:
$$
X_{F_f}(S, \theta) = v.
$$
Of course, this stronger condition could be exploited  in the discussion on properties of our
Kaehler metric over the moduli space (f.e., is it a metric of constant holomorphic sectional
curvature or not) but these questions come outside of the main theme of our text.
One verifies this condition using Propositions 3.3 and 3.5 namely
for any pair of smooth functions $\psi_1, \psi_2 \in C^{\infty}(S, \R)$
one has to find a smooth function $f$ over whole $M$ such that
$$
\aligned
\psi_1 = f|_S - const, \\
\psi_2 = \frac{Lie_{X_{in}}\theta}{\theta},\\
\endaligned
\tag 3.2.8
$$
(see (3.2.3)).  
The first equation in (3.2.8) can be easily solved taking any extension of $\psi_1$
over $M$ as the function $f$. The second one is more delicate: this question
extracts from the set of possible extensions those which are appropriate and
a priori one could not say whether or not such extensions exist.
We reduce the question to the following 
\proclaim{Lemma} Let $S$ be any smooth compact real oriented manifold and $\eta$ a volume form.
Then for any real smooth function $\psi \in C^{\infty}(S, \R)$ with zero integral:
$$
\int_S \psi \eta = 0,
$$
there exists a vector field $Y$ such that
$$
\psi = \frac{Lie_Y \eta}{\eta}.
\tag 3.2.9
$$
\endproclaim

Indeed, since 
$$
\aligned
\int_S \psi \eta = 0,\\
d(\psi \eta) =0,\\
\endaligned
$$
by the hypothesis of the lemma our $n$ - form $\psi \eta$ belongs to the trivial cohomology class. It means that
there exists a $n-1$ form $\zeta$ such that
$$
d \zeta = \psi \eta.
\tag 3.2.10
$$
On the other hand, decoding $Lie_Y \eta$ one gets
$$
Lie_Y \eta = d (\imath_Y \eta),
\tag 3.2.11
$$
hence it remains to compare (3.2.9), (3.2.10) and (3.2.11) taking in mind that $\eta$ has no zeros.
The desired vector field thus can be constructed locally to solve the equation
$$
\imath_{Y} \eta = \zeta
$$
which is obviously possible (notice again that $\eta$ is nonvanishing everywhere).

Now one can apply this lemma in our context using the following trick
(which will be very useful in the next subsection as well). Since our
$S$ is oriented (see the first step in the construction   of the moduli space $\Cal B^{hw,r}_S$)
we can consider the corresponding volume form $\eta$ instead of the square $\theta^2$.
The Lie derivatives of $\theta$ and $\eta$ are related by the identity
$$
\frac{Lie_Y \theta}{\theta} = \half \frac{Lie_Y \eta}{\eta}
\tag 3.2.12
$$
(see the next subsection), thus the lemma above ensures us that for every $\psi_2$ there exists
a vector field $Y$ on the Bohr - Sommerfeld cycle $S$ such that the second equation in (3.2.8)
is satisfied. Now it remains to construct such an extension of $\psi_1$ to a neighborhood
of $S$ in $M$  that the corresponding Hamiltonian vector field should give
us our $Y$ as its $X_{in}$. Since the consideration are local we will construct such extension
for a small neighborhood of zero section in $T^*S$. The desired function $\tilde f$ has the form
$$
\tilde f(x, p) = \psi_1(x) + p_x(Y_x),
\tag 3.2.13
$$
where $x$ is the $S$ - coordinate, $p$ is the coordinate along the fiber, identified
simultaneously with the corresponding cotangent vector $p_x$, and $Y$ is the vector
field on $S$ defined by $\psi_2$  which exists due to the lemma above. The standard isomorphism
of a neighborhood of $S$ in $M$ and the neighborhood of zero section in $T^*S$ maps
$\tilde f$ to a function which we denote as $f$, claiming that it possesses the properties
one needs to solve (3.2.8).

Thus  we can deform any fixed point $(S, \theta)$ in any direction along the moduli space,
acting by an appropriate induced quasi symbol. On the other hand,
here we'd like to enforce the statement which has been proved in [30] for homogeneous
symplectic manifolds. We mean that the map (2.3.6) is an inclusion. Here we 
establish it in full generality just slightly extending the arguments of [30].
Namely, it's not hard to see that $\Cal F_{\tau}$ could have a kernel if it were a point $x \in M$ with a neighborhood $\Cal O(x)$  such that
for every Bohr - Sommerfeld Lagrangian cycle $S \in \Cal B_S$  one has
$$
S \cap \Cal O(x) = \emptyset.
$$
This means that if one takes a bump smooth function concentrated in $\Cal O(x)$ then it has
trivial restriction on any Bohr - Sommerfeld cycle and therefore it belongs to the kernel.
But it is not the case if $\Cal B^{hw,r}_S$  is nonempty. Indeed, for any point $x$ we can
arrange a Bohr - Sommerfeld Lagrangian cycle passing any fixed neighborhood of $x$.
If $x$ is any fixed point and $S$ is a Bohr - Sommerfeld cycle
(which does exist since we assume that $\Cal B^{hw,r}_S$ is nonempty) then it's not hard
to construct a smooth function $f$ such that the corresponding flow generated by the Hamiltonian vector field
$X_f$ moves $S$ to the place of $x$. Due to the dynamical property of the Bohr - Sommerfeld condition
(see subsection 2.1) the image of the Bohr - Sommerfeld cycle should be Bohr - Sommerfeld again,
thus for any compact smooth symplectic manifold if $\Cal B^{hw,r}_S$ is nonempty then  
its "points" cover whole the based manifold. This remark is extremely important
if we are going to exploit a "universal cycle" induced by the construction
(and we are really going to do that in a future), --- this remark hints that
such a cycle does exist. Now, after this fact is understood, let us note that
if for every smooth function $f \in C^{\infty}(M, \R)$ there exists a Bohr - Sommerfeld
Lagrangian cycle $S$ such that the restriction $f|_S$ is non trivial then
one can choose a half weight $\theta$ over $S$ such that
$$
\int_S f|_S \theta^2 \neq 0.
$$
It's not hard to see that one can do the choice. 

Thus we've completed the proof of Theorem 1 modulo the statement of Proposition
3.5. The next subsection contains direct verification of the desired identity.

\subhead 3.3. Computations for Proposition 3.5
\endsubhead

First of all we remind again that the tangent space to the moduli space $\Cal B^{hw,r}_S$ at point $(S_0, \theta_0)$ is represented by pairs of smooth functions $(\psi_1, \psi_2)$ over $S_0$ satisfying the norm condition
$$
\int_{S_0} \psi_i \theta_0^2 = 0.
\tag 3.3.1
$$
Consider infinitesimal action of the Hamiltonian vector field $X_f$ generated by a smooth function $f \in C^{\infty}(M, \R)$. Again we mention the possibility to  split the tangent space $TM$ at the points of $S_0$ into two parts:
$$
X_f = V_f + W_f
$$
where 
$$
V_f = \om^{-1}(d(f|_{S_0})),
\tag 3.3.2
$$
and the second part 
$$
W_f = X_f - V_f
\tag 3.3.3
$$
is tangent to $S_0$ (it contains in $TS_0$ at points of $S_0$).
We've made an agreement to understand $V_f$ as the "outer" part of
the Hamiltonian vector field with respect to $S_0$ while
$W_f$ is the "inner" part, which preserves $S_0$ itself. Thus the inner
part acts on the objects which live over the fixed base $S_0$.
The infinitesimal deformation of $(S_0, \theta_0)$ under $X_f$ is reflected by
the coordinates of the dynamically correspondent field $\Theta_{DC}(X_f)$.
This deformation splits as well into two parts: the inner part and the outer part. As we've seen (see Proposition 3.3) the outer part represented by  vector field $V_f$ precisely corresponds to
$$
\psi_1 = f|_{S_0} - const,
\tag 3.3.4
$$
where the constant is fixed by (3.3.1) so
$$
const = \int_{S_0} f|_{S_0} \theta_0^2.
$$
The second, inner, component $\Theta_{DC}(X_f)$ has the shape
$$
\psi_2 = \frac{Lie_{W_f} \theta_0}{\theta_0},
\tag 3.3.5
$$
where the Lie derivative is defined by the inner part of the Hamiltonian vector field and is applied to a inner object - our half weight $\theta_0$.
To perform our computations successfully let us note that since
the source manifold $S$ from the definition of the moduli space (see 2.1)
is oriented we can consider $S_0$ as equipped with the corresponding orientation. In this case the square $\theta_0^2$ represents a volume form
over $S_0$ and we consider this volume form instead of $\theta_0^2$.
We denote it as $\eta_0$. Then one has the following smooth real function 
$$
L_f = \frac{Lie_{W_f} \eta_0}{\eta_0},
$$
corresponding to the logarithmic Lie derivative. At the same time
$$
Lie_{W_f}(\theta_0^2) = 2 \theta_0 \cdot Lie_{W_f} \theta_0,
$$
hence
$$
L_f = 2 \frac{Lie_{W_f} \theta_0}{\theta_0} = 2 \psi_2.
$$
Therefore we reformulate the statement of Proposition 3.3 substituting
$$
\psi_2 = \half L_f.
$$

Now let's compute the coordinates of the second vector field induced by
$f$ over the moduli space. So we take the corresponding induced function
$F_f$ and write down the coordinates of its Hamiltonian vector field.
Its differential $dF_f$ has been computed in Section 2, namely
$$
dF_f (\al_1, \al_2) = \tau \int_{S_0} f|_{S_0} 2 \al_2 \theta_0^2 +
\tau \int_{S_0} d \al_1 (W_f) \theta_0^2
\tag 3.3.6
$$
over point $(S_0, \theta_0)$ for tangent vector $v = (\al_1, \al_2)$.
Due to the simplicity of the expression for the symplectic form $\Om$
(see Section 2, (2.3.3)) one can easily convert the first summand
in (3.3.6), getting the first coordinate of the Hamiltonian vector field:
$$
\psi_1'  = 2 \tau f|_{S_0} - const',
\tag 3.3.7
$$
where the constant is defined by the same ruling condition.
The second summand requires a little bit more time. First of all,
$$
\aligned
\int_{S_0} d \al_1(W_f) \theta_0^2 = \int_{S_0} d \al_1 \wedge
\imath_{W_f} \eta_0 = \\
- \int_{S_0} \al_1 \cdot d(\imath_{W_f} \eta_0) = - \int_{S_0} \al_1
\frac{Lie_{W_f} \eta_0}{\eta_0} \theta_0^2,\\
\endaligned
\tag 3.3.8
$$
where at the second step we use the integration by parts. Substituting
this result to the expression (2.3.3) for the symplectic form we see that
the second coordinate of the Hamiltonian vector field is the following
$$
\psi_2' = \tau \frac{Lie_{W_f}\eta_0}{\eta_0} = 2 \tau L_f.
\tag 3.3.9
$$
It remains to compare (3.3.4) and (3.3.7), (3.3.5) and (3.3.9)
and as the result of the matching one gets the statement of
Proposition 3.5.

\subhead 3.4. Critical points of $F_f$
\endsubhead

The next question we'd like to discuss is the problem of critical
points and critical values for the induced quasi symbols over the moduli space.
This question so quite important for the quantum theory. Really, to perform
the measurement process the quantized observables should have enough 
eigenstates (= critical points, see Section 1). As well in any quantum mechanical
model one needs some requirement for the corresponding spectrum. As usual the case when the
spectrum is discrete, simple, etc. is understood as a good case.
 In our case there are two ways to arrange what one needs: first, one could find
an appropriate compactification of the moduli space; second, one can establish that
the moduli space itself carries the property that for sufficiently
generic smooth function on the based manifold the corresponding 
induced quasi symbol has a number of good  critical points. The situation when one has these two alternative ways are rather simuliar, say, to the crossroad which one meets in the Donaldson theory: from the first view point there exists a good compactification of the moduli space of instantons, constructed by K. Ulenbeck, but S. Donaldson computes the numbers, taking
such appropriate representatives for $\mu$ - classes which belong to
the moduli space itself (see [8]). Of course, this analogy is too rough and
pretentious but nevertheless we have to make a choice coming to the consideration
of the critical points. In this text we'll study the second way
(although the first one is  even more interesting).
To begin with we assume that for a sufficiently good function $f \in C^{\infty}(M, \R)$ (we'll deal with
Morse functions) the corresponding quasisimbol $F_f$ has critical points. Then we'll list
a number of good facts about these critical points. At the end of this subsection
we discuss the most important part --- the existence theorem.

So the first fact about the critical points is rather trivial:
\proclaim{Proposition 3.6} Pair $(S, \theta) \in \Cal B^{hw,r}_S$ is a critical point for
function $F_f$ iff the Hamiltonian vector field of the original function $f$
preserves this pair or in other words
$$
\aligned
f|_S = const,\\
Lie_{W_f} \theta \equiv 0.\\
\endaligned
\tag 3.4.1
$$
\endproclaim
The association "in other words" above has very important geometrical sense:
\proclaim{Lemma} Hamiltonian vector field  $X_f$ preserves Lagrangian submanifold $S$ if and only if
$S$ lies on a level set of $f$. 
\endproclaim
Thus in any case one should look for critical points on the level sets of a given $f$.
On the other hand, for any Morse function $f \in C^{\infty}(M, \R)$ one has the following
\proclaim{Proposition 3.7} Let $f$ be a sufficiently generic function over $M$.
 Then for  each generic value 
the corresponding critical set for $F_f$ is descrite.
\endproclaim
Here we assume that $\tau = \frac{1}{r}$ such that 
$$
f \equiv const \implies F_f \equiv const
$$
(see subsection 2.4). Thus we need to prove that if $(S, \theta) \in \Cal B^{hw,r}_S$
is a critical point for $F_f$ with critical value $v$ outside of the set of some special
values (we will specialize what these values are in a moment)
then it is isolated. 
The specialty of values means, firstly, that $v$ doesn't
belong to the set of critical values of $f$. Further, we divide the verification into two steps. Firsts, let's show that
the support $S$ is isolated. Indeed, if $S_{\de}$ lies sufficinely close to $S$
then it comes to some Darboux - Weinstein neighborhood of $S$. But all Bohr - Sommerfeld Lagrangian cycles  in this neighborhood have to intersect our given $S$. These cycles are represented by
exact  1 - forms but every exact 1 - form over a compact manifold has zeroes (at least two
since every real smooth function over a compact manifold has minimum and maximum  which give
at least two points where the differential of this function vanishes). Therefore $S_{\de}$
intersects $S$. At the same time $S$ and $S_{\de}$ belong to level sets of our function
$f$ being supports of critical points. If these level sets are different then
we come to the contradiction since two submanifolds from different level sets
never intersect each other. If $S$ and $S_{\de}$ belong to the same level set
$f^{-1}(v)$ then let's note that 

1) $X_f$ is a nonvanishing vector field being restricted to
non critical level set $f^{-1}(v)$;

2) $X_f$ is tangent to both $S$ and $S_{\de}$;

and we establish, comparing items 1) and 2),  that the intersection
$$
S \cap S_{\de} 
$$
either is trivial or its dimension is at least 1. If it is trivial then applying again
the arguments above we show that $S$ and $S_{\de}$ are distinct and separated. If
there is an intersection of positive dimension then the level where $S$ and $S_{\de}$
is special: it contains a limit cycles. Indeed, the intersection
should be invariant under the flow given by $X_f$. The theory of limit cycles
predicts that in general
the number of these limit cycles is finite. Thus if we
reject the level sets with these special values then for the remaining
generic value the level set could carry only a discrete set 
of the Bohr - Sommerfeld Lagrangian cycles. On the other hand
if $S$ is not isolated one gets even more than the existence of a limit cycle.
Namely one should have a family of deformations of $S$ inside of the level set;
this would mean that for our generic function $f$ and for the fixed level set
there is a Bott integral (or may be even more - there is an integral
for the dynamical system defined by $(M, \om)$ and Hamiltonian $f$).
Here we would like to decode prefix "sufficiently" in the statement of
Proposition 3.7: $f$ is "sufficiently generic" if the corresponding system
is non integrable. It was already mentioned in [30] that
the continuous deformations of a stable cycle $S$ with respect
to $f$ comes from the functions which commute with $f$. Of course,
if $g$ satisfies
$$
\{f, g\}_{\om} = 0
$$
it doesn't imply that $g$ gives a deformation of $S$ since it can happen
that
$$
g|_S = const.
$$
But in general the converse takes place: if we can deform 
a stable cycle on a generic level set that we can derive from this
deformation some integral for the system. But one should say
that  integrable system is not the general case  in classical mechanics.

The second step is very simple: we just need to apply Lemma from subsection 3.1
above together with Theorem 1. Indeed, if the support is fixed then the half weight part
could not vary continuously. Otherwise one gets that the critical set of
quasi symbol $F_f$ (since Theorem 1 ensures that it is) is Lagrangian at some point while according to
the lemma it has to be a complex submanifold.

\subheading{Remark} From the proof we derive  first constraints on the topology of
$S$. Indeed, as we've seen if $(S, \theta)$ is a critical point with critical value
$v$ which doesn't belong to the critical value set of our given $f$ then
the euler characteristic of $S$ should vanishes. Really, in this case $S$ admits
a nonvanishing everywhere vector field (namely $W_f$ which coincides with $X_f$
if $S$ belongs to the level set). Thus we get a first "vanishing"
result: 
if  the euler characteristic of $S$ is non trivial then
for every Morse function $f$ one has
$$
Spec F_f \subset Spec f \subset \R,
$$
where $Spec$ is the set of critical values.
On the other hand, one can easily deduce that either one could apply the same argument as above
to the case of the critical levels or one gets that there are no "critical" smooth Bohr - Sommerfeld
cycles belonging to the critical levels of our Morse function $f$. Indeed, let
$(S, \theta)$ be a critical point of $F_f$ such that $S \subset f^{-1}(v_i)$ where $v_i \in Spec f$.
Let $p_1, ..., p_d \in f^{-1}(v_i)$ are the critical points. Then it appears two possibilities:
either 
$$
S \cap (p_1, ..., p_d) = \emptyset
$$
or $S$ contains some of $p_i$'s. The first case is reduced to the case of noncritical levels
since $X_f$ vanishes only at $(p_1, ..., p_d)$. The second one is impossible in general since
it means that $S$ is not smooth. Therefore if we speak about  simply connected $S$
(say, $n$ - dimensional sphere) then we have to insert on the discussion
the question of an appropriate completion (or compactification) of the moduli space.
Here the simplest step should be the following: one has to consider the maps 
$$
\phi: S \to M
$$
which are smooth outside of  finite sets of points and continues everywhere
together with pure smooth maps (as it was in the definition of the moduli space, see
subsection 2.1). But we'd like to postpone the discussion on the completion problem
taking in mind that a priori  the case when one takes $S$ with zero euler characteristic
is much more attractive for us. Indeed, $n$ - torus is the classical (and primordial) example of
compact Lagrangian submanifold. So it seems that in this case we can avoid 
any completion procedure keeping the moduli space itself as the central geometric object.
At the end of the discussion we just mention that it remains a global problem in our considerations
--- the existence problem. This problem should be stated first of all
for the moduli space itself: for any compact simply connected symplectic manifold with integer form
one asks whether or not some middle homological class can be realized as a Lagrangian Bohr - Sommerfeld cycle of a fixed topological type. If  the answer is  the  affirmative one then one asks
for a sufficiently wide class of smooth functions over $M$ are there sufficiently many
critical points.  Some preliminary result in this way is given by
the following construction.

Let $(M, \om)$ is a simply connected compact symplectic manifold
with integer symplectic form. Then according to the Donaldson result (see [7])
there exists a level $k$ such that the power of the prequantization line bundle $L^k$ has a smooth section
$s$ such that its zero set $M_1 \subset M$ is a smooth symplectic submanifold.
Thus the pair $(M_1, \om_1 = \om|_{M_1})$ is again a symplectic manifold
with integer symplectic form. Hence we can apply again the argument.
Since $M$ is finite dimensional one can choose such level $k$ that
starting with
$$
M_0 = M, \om_0 = k \om
$$
one constructs a flag of symplectic submanifolds
$$
pt = M_n \subset M_{n-1} \subset ... \subset M_0
$$
such that every $M_i$ is presented by zeros of the corresponding section
$s_{i-1}$ of the prequantization line bundle $L_{i-1}$ over 
$$
(M_{i-1}, \om_{i-1} = \om_0|_{M_{i-1}}).
$$
Let $f$ be a generic smooth function over $M$ such that $pt$ belongs to
a non critical level set. Fix this level set and denote it as $N \subset M$.
The intersection of $N$ with $M_{n-1}$ is a smooth circle. Notice, that
the symplectic setup is much more flexible in contrast with the holomorphic case
so we can deform slightly every stepping stone in our pyramid ---
f.e., the second one --- to get the desired property. Namely if
the intersection $N \cap M_{n-1}$ is not smooth we can slightly deform
$M_{n-1}$ to get some smooth intersection. This smooth intersection
is a set of loops and we just choose one single circle to
be our first constructed pip. It is Lagrangian in $(M_{n-1}, \om_{n-1})$
just by the dimensional property. We denote it as $N_1$. Further, consider
the next stepping stone, including the picture to the symplectic manifold
$(M_{n-2}, \om_{n-2})$. Notice that by the definition
of the flag the normal bundle of $M_{n-1}$ in $M_{n-2}$ is trivial.
This means that there exists a neighborhood $\Cal O_{n-1}$
of $M_{n-1}$ in $M_{n-2}$ which is presented topologically as
the direct product $M_{n-1} \times S^2$. Then (perhaps, after 
small deformation) the intersection $N_2 = N \cap M_{n-2}$
is smooth submanifold of $M_{n-2}$ fibered over $N_1$ (we take, if it's necessary,
just an appropriate smooth component of the intersection). The point is
that in general $N_2$ is Lagrangian in $(M_{n-2}, \om_{n-2})$.
Indeed, the restriction of $\om_{n-2}$ to $\Cal O_{n-1}$
splits as the direct product of $\om_{n-1}$ and a symplectic form on
$S^2$. Since the intersection of any fiber $S^2$ with $N$ has dimension
1 it should be Lagrangian. Thus whole the intersection $N_2$ is Lagrangian,
being the product of two Lagrangian cycles. Let's remark that topologically
$N_2$ is  2 - torus which is topologically trivial. Repeating
the procedure toward the top of the pyramid one gets $n$ - torus $N_n$
in $(M, k\om)$  such that:

it represents trivial homology class;

it belongs to a level set of our function $f$.

Now, we didn't say anything about some connection on the prequantization bundle.
At the same time the ingredient could be included in the pyramid
such that the resulting Lagrangian cycle will be Bohr - Sommerfeld. It remains
to construct an appropriate half weight, invariant under the flow of
$W_f$ in the terminology of subsection 3.3. To get such a half weight
we can take any half weight $\theta_0$ on the constructed $S = N_n$
and then average it switching on the flow of $W_f$. Since
our $S$ is compact and $W_f$ is regular (without zeros)
the resulting half weight is correctly defined. By the construction
it should be invariant under the flow.

\newpage

\head 4. ALG(a) - quantization
\endhead

$$
$$

In this section we use all the results which were established above to perform 
 algebro geometric quantization of simply connected compact symplectic manifold
with integer symplectic form. We call this method "ALG(a) - quantization" since it is
a quantization in the framework of algebraic Lagrangian geometry (abelian case).
In subsection 4.1 we project all the results to the subject of the first part
of Section 1 and show that our construction looks like a quantization
in that terms. However the result of the quantization is quite consistent 
with the results of known geometric quantization constructions: in subsections 4.2 and 4.3 we
study how ALG(a) - quantization can be reduced in the case when the based manifold $M$ is equipped
with an appropriate polarization. We begin with the real case in 4.2 and then consider the case
of complex polarization in subsection 4.3. The rest subsection is devoted to
some applied questions: we discuss a natural quasi
classical limit of ALG(a) - quantization.

$$
$$
$$
$$

\subhead 4.1. The geometry of quantization
\endsubhead

Let's come back to Section 1. There we postulate that symplectic manifolds can be quantized
in terms of geometric formulation of quantum mechanics. Thus we ask for any symplectic manifold
about an appropriate algebraic manifold which should play the role of the corresponding quantum 
phase space. This quantum phase space should admit sufficiently many functions whose Hamiltonian vector fields
would generate its symmetries: we call these functions quasi symbols.  Moreover, we have to define
an appropriate map from the space of classical observables to the space of quantum observables
and this map has to be an irreducible representation of the given Poisson structure.
The irreducibility condition has been discussed there; here we just remind that irreducibility
means that the kernel of this map is trivial and there is no a smooth submanifold of the quantum phase space
which is invariant under the flow generated by quantized observables unless whole the quantum space.
After all this reformulation is understood let's turn to Algebraic Lagrangian Geometry, proposed
in [10], [25]. For any simply connected compact symplectic manifold with integer symplectic form
this programme gives  series of infinite dimensional algebraic manifolds which are called the
moduli spaces of half weighted Bohr - Sommerfeld Lagrangian cycles of fixed volume and topological type.
It looked like a very natural proposal to consider such  moduli space as the quantum phase space for
some algebro geometric quantization. One gets the same feeling just comparing the materials of Section 1
and Section 2. Further, in Section 2 we construct  a map
$$
\Cal F_{\tau}: C^{\infty}(M, \R) \to C^{\infty}(\Cal B^{hw, r}_S, \R)
$$
and prove that it is a homomorphism of Lie algebras (it maps Poisson bracket to Poisson bracket).
Moreover,   in Section 3 we prove that the image of $\Cal F_{\tau}$ belongs to the subalgebra of quasi symbols and that $\Cal F_{\tau}$ is irreducible. On the other hand, Theorem 1
shows that $\Cal B^{hw,r}_S$ is a good candidate to be the quantum phase space of some
quantum mechanical system. Indeed, we prove (see there) that for every point $(S, \theta) \in
\Cal B^{hw,r}_S$ and every tangent vector at this point $v \in T_{(S, \theta)} \Cal B^{hw,r}_S$
there exists a quantum observable (moreover, a quantized classical observable)
whose Hamiltonian vector field  coincides at this point with this vector. It means roughly  that
this quantum phase space possesses maximal quantum symmetries.

Now to perform some additional construction, making the quantization more familiar for the audience, we need the fact which has been mentioned many times above:
\proclaim{Proposition 4.1 ([10], [25])} The symplectic form on the moduli space
$\Cal B^{hw, r}_S$ has integer cohomology type.
\endproclaim

To prove the statement it sufficient to remind the definitions from Section 2.
The symplectic structure on the moduli space of half weighted planckian cycles
is a double cover of the canonical symplectic structure on a cotangent bundle
(see subsection 2.2) thus it is integer (however it corresponds to
the trivial topological line bundle). Then we apply the symplectic reduction
getting our moduli space $\Cal B^{hw,r}_S$ but during the procedure the structure
remains to be integer. From these defining constructions one gets the following
\proclaim{Proposition 4.2 ([10], [25])} The prequantization line bundle
 corresponding to symplectic form (2.3.3) is exactly the complex line bundle,
associated to $U(1)$ - principle Berry bundle, defined in (2.1.4).
\endproclaim

Notice that an integer symplectic manifold, given by our moduli space $\Cal B^{hw,r}_S$
is almost ready to be quantized: it is canonically equipped
with a real polarization 
$$
\pi: \Cal B^{hw,r}_S \to \Cal B_S
$$
and a complex polarization, given by the complex structure. It remains to 
choose an appropriate prequantization connection and then start a known programme.
But we would like to mention here an additional aspect coming
when one applies Berezin - Rawnsley method. Namely
if $F \in C^{\infty}(\Cal B^{hw,r}_S, \R)$ is a quantized function by Berezin - Rawnsley
(and it is a quasi symbol in our terminology) then it gives a self adjoint operator
on the space of holomorphic sections of the Berry bundle. Indeed, if one fixes
an appropriate prequantization connection then one could perform
the Berezin - Rawnsley construction getting desired correspondence.
But we can restrict the investigations to the subset in $C^{\infty}_q (\Cal B^{hw,r}_S, \R)$
which consists of the images of $F_{\tau}$. Then one has the following composition:
$$
\aligned
C^{\infty}(M, \R) \hookrightarrow & C^{\infty}_q(\Cal B^{hw,r}_S, \R) \to Op(\Cal H),\\
f \mapsto & F_f \mapsto \hat Q_{F_f},\\
\endaligned
\tag 4.1.1
$$
where $\Cal H$ is the space of holomorphic sections of the Berry bundle and
$\hat Q_f$ is the corresponding self adjoint operator, given by the Berezin
- Rawnsley procedure. The possibility to construct this composition
is the reason why we emphasized everywhere that the Kaehler manifold
desired in AGQ should be algebraic. In this case there is a superstructure
which makes it possible to come back, getting some geometric quantization
procedure. Indeed, the composed correspondence given by (4.1.1)
is an example of a representation of the Poisson algebra in some
Hilbert space. It satisfies a number requirements from the Dirac list;
f.e., it's obviously linear; at the same time we can manage
that constants would come to the same constants; it is a homomorphism
of the Lie algebra (since both the maps in the first row of (4.1.1)
possess the property). As well one could prove that it is irreducible,
but it is another story and we will return to it in a future.

Now we are going to discuss how ALG(a) - quantization is related to the known
methods, recalled in Section 1.

\subhead 4.2. Real polarization
\endsubhead

Assume now that our symplectic manifold $(M, \om)$ admits an appropriate real polarization.
It means that $(M, \om)$ can be equipped with a Lagrangian distribution that is a field of Lagrangian
subspaces in the complexified tangent bundle which is integrable. The complex and the real case
are differ by the nature of these Lagrangian subspaces: in the real case all of them
are real while in the complex one all of them are pure complex. We've recalled both the notions in Section 1. Turning to that we understand the real polarization case as that one when on $M$ there are
a set of smooth functions $f_1, ..., f_n$ such that
$$
\{f_i, f_j \}_{\om} = 0, \quad \quad \forall i, j,
\tag 4.2.1
$$
defining some Lagrangian fibration
$$
\pi: M \to \De,
\tag 4.2.2
$$
where $\De \subset \R^n$ is a convex polytope. For each inner point 
$$
(t_1, ..., t_n)  \in \De \backslash \partial \De
$$
the corresponding fiber
$$
\pi^{-1}(t_1, ..., t_n) = f^{-1}_1(t_1) \cap ... \cap f^{-1}_n(t_n)
\tag 4.2.3
$$
is a smooth Lagrangian cycle. The degenerations at the faces of $\De$ are regular so
the inner part of each $n-k$ - face corresponds to $n-k$ - dimensional
isotropic submanifolds. Usually (see [20]) the quantization scheme in this case
distinguishes some special fibers of (4.2.2) namely Bohr - Sommerfeld fibers
(we've discussed this subject in Section 1). However as we prove below
there are finitely many such fibers for the compact case hence one has just
a  finite set $S_1, ..., S_l$. Then the Hilbert space of the quantization is
$$
\Cal H = \sum_{i=1}^l \C <S_i>,
\tag 4.2.4
$$
and $S_i$ play the role of a canonical basis. The functions $\{ f_i \}$ are represented then by
diagonal operators; for some other functions one can define operators which are not in general
unitary, see Section 1.

Now we wish to apply ALG(a) - programme to our completely integrable system taking 
the homology class of the fiber as $[S]$. Then one gets the corresponding moduli space
$\Cal B^{hw,1}_S$ where
$$
[S] = [\pi^{-1}(pt)] \in H_n(M, \Z)
$$
and we take volume 1 just for simplicity. The set $(f_1, ..., f_n)$ defines
 quasi symbols $F_{f_1}, ..., F_{f_n}$  which are again in involution (and moreover,
one can manage an infinite set of functions in involution taking all finite products
of $f_1, ..., f_n$ and then mapping them by $\Cal F_{\tau}$). Denote as
$Crit (F_{f_i})$ the set of critical points of $F_{f_i}$. Consider the following intersection
$$
P = Crit(F_{f_1}) \cap ... \cap Crit(F_{f_n}) \subset \Cal B^{hw,1}_S
\tag 4.2.5
$$
so the mutual critical set. One has
\proclaim{Proposition 4.3 ([30])} The set $P$ is a double cover of the set $\{S_i\}$
consists of the Bohr - Sommerfeld fibers.
\endproclaim

 Thus in general we can recover the well known method from [20]: one just takes the mutual critical
set for the distinguished functions which preserve given real polarization and the supports of
these critical points correspond to the basis in $\Cal H$. Moreover, the proposition
is true in  general non compact case (but we do not consider it here) so one could try
to expolite this correspondence in some different context.
On the other hand, one can deduce what happens for some other function $f_0$ which doesn't belong
to the algebraic span of $\{f_1, ..., f_n\}$. Say, one can deform the given system $\{f_1, ..., f_n\}$
in involution using such a function if the set
$$
g_i = \{f_0, f_i\}_{\om}
$$
consists of functions in involution. Then this set defines some other real polarization which is given
by the action of the flow of $X_{f_0}$ on the original one. Then if we are lucky 
one can define an operator, corresponding to $f_0$ in terms of $\Cal H$. But here we would
like to discuss only general facts so at this step we know how to quantize
only the functions which preserve the polarization.

The proof of Proposition 4.3 uses the equations for  critical points of quasisymbols given
in subsection 3.4.  Since for every Bohr - Sommerfeld Lagrangian fiber $S_i$ all $f_i$
are constant along it then the first equation from (3.4.1) is satisfied automatically. To
complete the proof we need to find an appropriate half weight $\theta_i$ on $S_i$ such that
the pairs $(S_i, \pm \theta_i)$ should be invariant under all Hamiltonian vector fields $X_{f_j}$.  
Note, that since $S_i$ lies on a level set for any $f_j$ then
$$
X_{f_j} \equiv W_{F_j}
$$
at the points of $S_i$. Really we are looking for an appropriate volume form
on $S_i$ using again the same argument as before.
Now take over $S_i$ the set of differentials $df_i$ and combine $n$ - form
$$
\tilde{\eta} = df_1 \wedge ... \wedge df_n.
\tag 4.2.7
$$
From the properties of the completely integrable system one sees that this $n$ - form
is non degenerated everywhere over $S_i$ (the order of $f_i$'s is fixed by the given
orientation, see again the definition of the Bohr - Sommerfeld cycles). Then there exists
such $n$ - form $\eta$ that
$$
d\mu_L = \eta' \wedge \tilde{\eta}
\tag 4.2.8
$$
over $S_i$, where $d\mu_L$ is as above the Liouville form. Of course, this $\eta'$ defined by (4.2.8)
is not unique but its restriction
$$
\eta = \eta'|_{S_i}
\tag 4.2.9
$$
is defined uniquely. By the construction this volume form is invariant under the action
of every $X_{f_i}$. Indeed, both the Liouville form and the intermideate $\tilde{\eta}$
are invariant under the flows thus the last $\eta$ possesses the same property. It remains
to normalize $\eta$ comparing given volume $r$ with the volume
$$ 
\int_{S_i} \eta_i.
$$
Then one takes the square roots 
$$
\pm \theta_i
$$
such that
$$
\theta_i^2 \equiv \eta_i.
$$
Let show that there are no any other invariant half weights over $S_i$.  Indeed, let $\theta_i'$ is another half weight
over $S_i$ which is invariant under each $X_{f_j}$. Then the ratio
$$
\psi = \frac{\theta_i'}{\theta_i}
$$
is a smooth function over $S_i$, satisfying
$$
Lie_{W_{f_j}} \psi = 0 \quad \quad \forall j = 1, ..., n.
\tag 4.2.10
$$
But the set $\{X_{f_j}\}$ form some local basis of the tangent space  at each point of $S_i$ 
 thus $\psi$ has to be constant. The normalizing
condition implies that this constant is either plus or minus 1 hence it remains
only one pair $\pm \theta_i$ as possible solutions.

Conversely, let $(S_0, \theta_0)$ be a mutual critical point for all $F_{f_j}$'s. Then, again
applying (3.4.1) we get that all functions $f_j$'s are constant over $S_0$. Forgetting
about the second component one gets
$$
f_j|_S = const = t_j \quad \quad \forall j = 1, ..., n,
$$
consequently $S$ is a Bohr - Sommerfeld fiber:
$$
S = \pi^{-1}(t_1, ..., t_n).
$$
It completes the proof of Proposition 4.3.

Further, we've mentioned that in the compact case the set of Bohr - Sommerfeld fibers
is finite. Although it doesn't lie on the mainstream of our  consideration
we present here
\proclaim{Theorem 2 ([30])} Let $(M, \om)$ be a symplectic manifold with integer symplectic form which admits
Lagrangian fibration
$$
\pi: M \to \De
$$
with compact fibers. Then the set of smooth  Bohr - Sommerfeld fibers is discrete.
\endproclaim

The proof is very short in the case when $(M, \om)$ is a completely integrable system: the statement follows from Proposition 3.7 and Proposition 4.3. 
In our more general case one uses the argument which has been exploited in the proof of
Proposition 3.7: if $S_0$ is a Bohr - Sommerfeld fiber so
$$
S_0 = \pi^{-1}(p_0), \quad \quad p_0 \in \De
$$
then there exists a neighborhood $\Cal O(p_0)$ of the point in $\De$ such that
$\pi^{-1}(\Cal O(p_0)$ is a Darboux - Weinstein neighborhood of $S_0$
in $M$. Then if we suppose that there is another Bohr - Sommerfeld fiber $S$
projecting to $p \subset \Cal O(p_0) \subset \De$ then it should be 
a smooth function $\psi \in C^{\infty}(S_0, \R)$ such that $S$ coincides
with the graph of $d\psi$ in this Darboux - Weinstein neighborhood.
Since $S_0$ and $S$ have zero intersection being two different fibers
the differential $d \psi$ has to be nonvanishing everywhere. But any
smooth function on a compact set has at least two extremum points: the minimal and
the maximal ones. This means that $d \psi$ has to vanish somewhere which leads
to the contradiction. Therefore if our $S_0$ is Bohr - Sommerfeld then
there exists a neighborhood of $\pi(S_0) = p_0$ in $\De$ such that
$p_0 \in \Cal O(p_0)$ is unique "Bohr - Sommerfeld point" in this neighborhood.
Thus, globally, every Bohr - Sommerfeld fiber of $\pi$ is separated by such a neighborhood
and hence the set of Bohr - Sommerfeld fibers is discrete. Moreover if $M$ is compact
it follows that the set is finite.

\subheading{Remark} It's quite natural and reasonable to continue here the observation
given at the end of subsection 2.3 above. There we spoke about the case of completely
integrable systems: in this case one could construct, starting with the given
set of first integrals $\{f_1, ..., f_n\}$, an infinite set of commuting quantum observables
over the moduli space $\Cal B^{hw,r}_S$. Indeed, one just takes 
$$
\{F_{f_j^k} \}, \quad \quad j = 1, ..., n, k>0
\tag 4.2.11
$$
and for every pair from this set the quantum Poisson bracket vanishes (see subsection 2.3).
However it's clear that the mutual critical set $P$ is the same for every degrees
$k$ and any combinations of the first integrals (while the corresponding quantum
observables are not longer algebraically dependent). Really, the conditions
$$
f_j|_S = const
$$
and
$$
f_j^k|_S = const 
$$
are absolutely equivalent (since our functions are real and smooth). And despite of the fact
that quasisymbols of type (4.2.11) are algebraically independent, their critical values
{\it are} algebraically dependent in mutual critical points. Indeed, every first integral
$f_j$ gives the following critical values (via the powers of $f_j$):
$$
c=f_j|_S, c^2, ..., c^k, ...,
$$
and it's clear that this set is algebraically dependent. Hence we could not
derive some additional geometric information for the completely integrable systems
(at least in the present discussion) using out method.

\subhead 4.3. Complex polarization
\endsubhead

This case belongs exactly to algebraic geometry. Complex polarization is a choice on the symplectic manifold $(M, \om)$ of any compatible with $\om$ integrable complex structure $I$
(so one suggests that such a structure exists), transforming our $M$ to an algebraic manifold.
While any symplectic manifold admits an infinite set of almost complex structure
the possibility to choose in this set some integrable one makes the horizon of the examples
much less wider. This condition implies that $(M_I, \om)$ is a Kaehler manifold
and the integrability condition for $\om$ ensures that the Kaehler metric 
has the Hodge type hence $(M_I, \om)$ is an algebraic variety. 

Known quantization methods have been discussed in Section 1. It was mentioned
there that they are based on some reductions of the basic Souriau - Kostant
method. Thus we'll use here both the reductions: Berezin - Rawnsley method
is more appropriate to claim some dynamical coherence while
Berezin - Toeplitz method is described by some explicit
formulas, see Proposition 2.3. In any case the corresponding Hilbert space
is the same --- the space of holomorphis sections of the prequantization line bundle
with respect to  the prequantization connection. We porjectivize the space
following the strategy of Section 1. The first step is to
relate the quantum phase spaces of the known method and ALG(a) - quantization.
The desired relationship is given by so-called BPU - map ("BPU" means
the first letters of the authors names, see [5]). The moduli space
$\Cal B^{hw,r}_S$ is fibered over the projective space:
$$
BPU: \Cal B^{hw,r}_S \to \proj H^0(M_I, L).
\tag 4.3.1
$$
Recall the construction following [10], [25]. Let $s \in H^0(M_I, L)$
be a holomorphic section of the prequantization bundle. Restrict it
to on any half weighted Bohr - Sommerfeld cycle $S, \theta)$. The restriction
is represented by a smooth complex function. Indeed, it is true for any smooth section of $L$: the pair of prequantization data $(L, a)$ is restricted to 
$S$ as trivial bundle with flat connection admitting a covarinatly constant
trivialization by the definition hence every section over $S$ is presented by
a covarinatly constant section multiplied by a smooth complex function.
Since the trivialization over $S$ is defined up to scale one should lift
the consideration to the set of half weighted planckian cycles to kill the ambiguity. Then the restriction of any section to $\tilde S$ is exactly
a complex function. Then we define a map
$$
\Cal P^{hw}_S \to H^0(M_I, L)
\tag 4.3.2
$$
by the following condition
$$
(\tilde S_0, \theta_0) \mapsto s_0
$$
  iff
$$
 \int_{\tilde S_0} s|_{\tilde S_0} \theta_0^2 = \int_M <s, s_0> d\mu_L
\tag 4.3.3
$$
for any holomorphic $s \in H^0(M_I, L)$. Since $H^0(M_I, L)$
is finite dimensional (we work over a compact symplectic manifold) for
any $(\tilde S_0, \theta_0)$ such holomorphic section exists.
Indeed, every $(\tilde S_0, \theta_0)$ defines a linear functional
on $H^0(M_I, L)$ given by the left hand side of (4.3.3) hence
one has a fibration
$$
\Cal P^{hw}_S \to (H^0(M_I, L))^*
\tag 4.3.2'
$$
and then one can get (4.3.2) applying a standard identification since
$H^0(M_I, L)$ is equipped with "quantum" hermitian form $<,>_q$ (see
Section 1).  

It's clear that the canonical $U(1)$ - actions on the source 
and on the target in (4.3.2) are compatible.
Hence one could take the corresponding Kaehler reductions
of the both spaces:
then we get a map from $\Cal B^{hw,r}_S$
to the projective space $\proj H^0(M_I, L)$ which is the result
of the Kaehler reduction of $H^0(M_I, L)$.
This is BPU - map (4.3.1)
which gives a reduction of the quantum phase space
of ALG(a) - quantization and the one of the known method.
 
 Further, let smooth function $f$ be a quantizable 
observable (as in Section 1) with respect to the given complex polarization. Then the Hamiltonian
vector field $X_f$ preserves the complex structure $I$ hence the dynamically correspondent
vector field on $\proj(\Ga(M, L))$ preserves the finite dimensional piece $\proj(H^0(M_I, L))$.
Moreover, the field $\Theta^p_{DC}(X_f) \in Vect(\proj(H^0))$ preserves
whole the Kaehler structure and corresponds (see Section 1) to some smooth function $Q_f$ (let's emphasize
again that it is true since $f$ is quantizable). On the other hand, the Hamiltonian vector field
defines infinitesimal transformations on both the ingredients of the BPU - map. 
 For each quantizable function the corresponding dynamical actions
on the source space and the target space have to be compatible. Thus for any quantizable function
(in the sense of Rawnsley - Berezin method) one has:

1) a pair of quasisymbols $F_f$ and $Q_f$ on the source and the target spaces respectively;

2) a pair of dynamically correspondent vector fields $\Theta_{DC}(f)$ and
$\Theta^p_{DC}(f)$ on the source and the target spaces respectively. 

Taking into account the dynamical arguments we get that the differential
of BPU - map translates our special vector field $\Theta_{DC}(f)$
to the special vector field $\Theta^p_{DC}(f)$. Propositions 1.6. and 3.5. ensure that there is
\proclaim{Proposition 4.4 ([30])} For any quantizable function $f$ the Hamiltonian vector fields
of quasisymbols $F_f$ and $Q_f$ are related as follows
$$
dBPU(X_{F_f}) = c \cdot X_{Q_f},
\tag 4.3.4
$$
where $c$ is a real constant.
\endproclaim

Thus one reduces ALG(a) - quantization to a well known method in the complex polarization case.
As a consequence of Proposition 4.4 we get how to find the eigenstates of quantum observable
$Q_f$ having the eigenstates of $F_f$. The relationsheep is very similar to the answer
in the previous real case:
\proclaim{Corollary 4.5 ([30])} BPU - map projects the set of eigenstates of the quantum observable $F_f$
to the set of eigenstates of $Q_f$.
\endproclaim

Note that we get the statement of Proposition 4.4 just directly
from the dynamical arguments. It's quite hard to compute explicitly the differential,
but it's clear that the definition of BPU - map and the definition of our
quasisymbols over the moduli space are sufficiently close. Moreover, one could see
that as one defines a function $F_f$ over the moduli space starting with a smooth function
on the given symplectic manifold one can define a section of a vector bundle
on the moduli space starting with any section of the prequantization bundle.
We've almost seen what the bundle over the moduli space is: for every section
$s \in \Ga(L)$ there is a submanifold $K_s \subset \Cal B^{hw,r}_S$. 
The topological type of the submanifold is obviously fixed. Thus the Chern classes
of the bundle could be defined using the topology of these submanifolds together
with  their intersection theory. But one could define the corresponding bundle
taking into account just local considerations in the Darboux - Weinstein
neighborhoods. We just hint here some topological relationships given in the picture:
namely for any symplectic manifold $X$ with integer symplectic form there two distinguished
2 - cohomology classes, the class of the symplectic form and the associated
canonical class. What are the classes for our moduli space $\Cal B^{hw,r}_S$?
The symplectic class is represented by the Berry bundle. At the same time,
one hints that the associated canonical bundle is given by our prequantization
bundle $L$. Then from BPU - map one could relate the Berry bundle with the canonical bundle
over given $M$ associated with given symplectic form. It gives an interesting
duality: during the process the canonical bundle turns to be
a prequantization bundle while the prequantization bundle turns to be
the canonical bundle over the moduli space.

The construction of BPU - map outlines a way how to understand some properties of the moduli space
which we derived above. Namely, let us consider instead of the ruling map (4.3.2) the following one
$$
\Cal P^{hw}_S \to \Cal H,
\tag 4.3.5
$$
where $\Cal H$ is the Souriau - Kostant space of wave functions (see subsection 1.5).
Indeed, apply the same scheme: every planckian cycle gives a linear map
from $\Cal H$ to complex numbers. Thus one gets
$$
\Cal P^{hw}_S \hookrightarrow \Cal H^*,
\tag 4.3.6
$$
and the question is whether or not one can use here the corresponding duality property
to get (4.3.5). If it the case then we again perform the Kaehler reduction and get
a fibration
$$
gBPU: \Cal B^{hw,r}_S \to \proj (\Cal H).
\tag 4.3.7
$$
Then one expects that

1) this map is a double cover;

2) this map is $\al$ - holomorphic, so the differential maps
holomorphic directions to  directions
with constant Kaehler angle;

3) for every smooth function $f$ the quasi symbols $F_f$ and $Q_f$ on the corresponding quantum phase spaces are related by 
$$
d(gBPU)(X_{F_f}) = c \cdot X_{Q_f}.
$$

This result would give very good compactification of the moduli space.
This compactification would possess very good property:
over it any induced quasi symbol should have a lot of
eigenstates. And at last it would explain three things: why the moduli space carries 
a hermitian metric of constant holomorphic sectional curvature; why the Souriau
- Kostant quantization is reducible; and why the moduli space 
is a good candidate to solve AGQ - problem.

At the same time the complex polarization case has another face: the point is that for any choice of
a compatible (almost) complex structure over our given symplectic manifold gives simultenously 
the corresponding riemannian metric, compatible with $\om$. This gives us the following specialization
for some half weighted Bohr - Sommerfeld Lagrangian cycles: we denote as $\Cal B^g_S$ the subset in the moduli space which consists of such pairs $(S, \theta) \in \Cal B^{hw,r}_S$ that 
$$
d\mu(g, S) \equiv \theta^2
\tag 4.3.8
$$
where $d\mu(g, S)$ is the volume form, given over $S$ by the restriction of our riemannian metric $g$ to
our cycle $S$. Turning to the moduli space of "unweighted" Bohr - Sommerfeld Lagrangian cycles
one can see that in the presence of $g$ the moduli space is fibered over real numbers:
$$
\aligned
W_g: \Cal B_S \to \R_+,\\
W_g(S) = Vol_g S \\
\endaligned
\tag 4.3.9
$$
(since $S$ is oriented). Then one has the following simple
\proclaim{Lemma} The subset $\Cal B^g_S \subset \Cal B^{hw,r}_S$ is represented as the double cover
of  $W_g^{-1}(r)$.
\endproclaim

Indeed, $(S, \theta)$ belongs to $\Cal B^g_S$ if and only if (4.3.8) holds. But
if $S \in W_g^{-1}(r)$ then there are exist exactly two half weights
such that (4.3.8) takes place. And if $S$ doesn't belong to $W_g^{-1}(r)$ then
the equality (4.3.8) is impossible for any half weight over it.
Thus one gets the picture
$$
\Cal B^g_S \to W^{-1}_g (r) \subset \Cal B_S.
$$
Then for this subset $\Cal B^g_S \subset \Cal B^{hw,r}_S$
one could formulate a number of valuable statements and include it to the quantization picture.

\subhead 4.4. Quasi classical limit
of ALG(a) - quantization
\endsubhead

Here we list some remarks on the limit of our method when an appropriate parameter
goes to infinity. We discuss how the picture
changes with respect to the level $k$ in subsection 2.4. Here we
translate the mathematical dependences on the parameters to
quasi classical background.

Following Berezin we think that level $k$ is inversely proportional
to the Planck constant. We distinguish some dependence
above (see subsection 2.4) introducing parameter $\tau$
in the definition of $\Cal F_{\tau}$. Now the time
is to recover that our parameter $\tau$ has to be proportional to
$k$. Then the formulation of the Dirac principle (see
Propositions 2.4) turns to be more familiar
in the framework of quantization. At the same time
one has that during the limit the volume of 
Bohr - Sommerfeld Lagrangian cycles decreases. Really, according to
condition (2.4.7) $r$ tends to zero while $\tau \to \infty$.
Since for every half weight $\theta$ the square $\theta^2$
is "positive" it means that during the limit one kills 
the half weight part. Hence in the limit one gets the moduli space
of unweighted cycles. At the same time the result is based  on not only Bohr -
Sommerfeld Lagrangian cycles. The point is that when level $k$ tends to infinity
the moduli space $\Cal B_{S, k}$ covers whole the set
of Lagrangian cycles. Let explain it more carefully.

Recall that since a prequantization data are fixed then one has the following map
$$
\chi: \Cal L_S \to J_S.
\tag 4.4.1
$$
Indeed,  the restriction of the prequantization line bundle equipped with
a prequantization connection gives a pair (trivial line bundle, 
flat connection) on the embedded Lagrangian cycle $S$. The flat connections
modulo gauge transformations are described by so - called Jacobian of $S$:
$$
J_S = H^1(S, \R) / H^1(S, \Z).
$$
It can be seen more explicitly if we fix a basis in $H_1(S, \Z)$. 
Then for any flat connection one has a set of numbers which are just
the result of the integration of our flat connection in any trivialization
along the basic submanifolds. The changing of the trivialization
gives nothing since we are working modulo integer points of the lattice
$H^1(S, \Z)$. Thus the map
$$
\chi: \Cal L_S \to J_S
\tag 4.4.2
$$
is defined. And the moduli space of Bohr - Sommerfeld Lagrangian cycles 
is just the preimage 
$$
\Cal B_S = \chi^{-1}(0).
$$
 Further, what happens when we go up to the higher levels?
Then it's clear that for any level $k$ one gets
$$
\Cal B_{S, k} = \chi^{-1}(0) \cup \chi^{-1}(p_1) \cup ... \cup \chi^{-1}(p_{k^{b_1} - 1}),
\tag 4.4.3
$$
(here it is a formal expression since $B_{S, k}$ is connected) where $p_i$'s are  the points of order $k$ on the Jacobian 
(and there are $k^{b_1}$ such points including zero) and $b_1$ is the first Betty number
of $S$. It's clear that when $k$ tends to infinity the order $k$ point set 
covers densely the torus $J_S$ and consequently when we arrange the procedure
the moduli space $\Cal B_{S, k}$ covers  densely the moduli space of
Lagrangian cycles.

Now the question is: what are the limits of the quantum observables
$F_f$ during the process? And what about the limiting Poisson bracket
on the moduli space $\Cal L_S$? One gets the answer to both the questions
considering the following objects over the moduli space of Lagrangian cycles.
Namely each smooth function $f \in C^{\infty}(M, \R)$ generates a special
object on $\Cal L_S$ which possesses two different natures. From
the first viewpoint one has just a vector field $Y_f$ defined by the restrictions
of $f$ to the Lagrangian cycles (we denoted  the same field as $A_f$ in subsection
2.3). Let's remind, the restriction of $f$ to $S \in \Cal L_S$ 
gives the corresponding  Hamiltonian (isodrastic) deformation of $S$
(given by $d(f|_S)$ in the Darboux - Weinstein neighborhood);
hence $Y_f$ at the point just equals to this tangent vector. On the other
hand, $Y_f$ is not a single vector field: the point is that
at the points where $Y_f$ vanishes as a vector field one has some numerical values.
Indeed, $Y_f$ vanishes as a vector field at $S$ if and only if
our function $f$ is constant being restricted to $S$. But it means that
it gives some number equals to this constant. Thus
the induced object $Y_f$ is described by a pair
$$
Y_f = (Y_f^0, Y_f^1),
\tag 4.4.4
$$
where $Y_f^0$ is a real function (sufficiently singular, of course)
and  $Y_f^1$ is a vector field (absolutely smooth, of course).
Let us denote the set of all such objects (given by smooth functions
on $M$) as $C^q(\Cal L_S)$. Then one has the following
\proclaim{Proposition 4.5} The set $C^q(\Cal L_S)$
is a Lie algebra.
\endproclaim

To ensure that the fact takes place  we
note firstly that the correspondence 
$$
f \mapsto Y_f
$$
is obviously linear. The operation  $[ , ]$ is given
just by the formula
$$
[Y_f, Y_g] = Y_{\{f, g\}_{\om}}.
\tag 4.4.5
$$
The Jacobi identity is satisfied just by the definition.

On the other hand, we should emphasize that:
\proclaim{Proposition 4.6} 

1) For every $f$ the corresponding object $Y_f$
is the natural result of the limiting procedure, applying
to quasisymbol $F_f$;   

2) the Lie bracket (4.4.5) is the natural result of
the limiting procedure, applying to the quantum Poisson bracket $\{ , \}_{\Om}$.
\endproclaim

Further, we see that the system based on $\Cal L_S$ is equipped with some
dynamical properties coming from the classical dynamics of the given
classical mechanical system. Indeed, if one choose a Hamiltonian
$H \in C^{\infty}(M, \R)$ then it generates a dynamics on $\Cal L_S$,
preserving the Lie bracket on the space of objects over the moduli space.
It's just a simple exercise in the computational technique
which we adopt during whole the text. Therefore we understand
the process as an appropriate quasi classical limit of ALG(a) -
quantization method: during the limiting procedure we lose
the measurement aspects but we keep the dynamical properties
compatible (more then compatible, we'd say) with 
dynamics  of the given system.

At the end of the story we have to mention that
we just present here some introductory part
of a new quantization method, skeeping a lot of additional questions
and details which will be clarified and established (if any)
in a future. Overing the discussion we claim that we set up
the problem having in mind some (may be accidental) coincidences listed above
but we are pretty sure that the study of algebraic Lagrangian geometry,
introduced in [10], [25], will lead  to new and interesting results.

\newpage

\head Addendum: supergeometry
\endhead

$$
$$

In this small additional part we mention that some geometrical notions and objects
from the main text can be naturally understood in terms of supergeometry.
Here we use the supergeometry just as a convenient language
which allows us to simplify  the definitions but at the same time
to outline some ways how one could generalize the constructions.

$$
$$
$$
$$

\subheading{Bohr - Sommerfeld condition} We start with the basic
for the method  condition. As it was pointed out 
it is a dynamical condition compatible with
Hamiltonian deformations from the "classical" position. At the same time it could be understood in terms of
even supersymplectic manifold. Again, let $(M, \om)$ is a symplectic manifold
with integer symplectic form. Again, let's take some prequantization data $(L, a)$
(see Section 1) such that 
$$
c_1(L) = [\om]
$$
and
$$
F_a = 2 \pi i \om.
$$
Consider the associated principle $U(1)$ - bundle
$$
\pi: P \to M;
$$
it is equipped with the same connection, represented now
by a pure imaginary 1 - form $A$. This principle bundle equipped with the connection
(compatible with the fixed hermitian structure on $L$) is an example of
 even supersymplectic manifold. The supersymplectic form
is divided at each point into two parts: the connection defines 
a splitting
$$
T_pP = T_{hor} \oplus T_{ver}
$$
of the tangent space  and then the supersymplectic form $\Om$ is defined
as the direct sum of the usual symplectic pairing on the horizontal part
plus the natural symmetric $U(1)$ - invariant pairing
defined by the natural metric on the vertical component.
It can be checked directly that the Jacoby identity for
the corresponding super bracket is implied just by the compatibility condition
on the hermitian structure and our prequantization connection. To be more familiar
the picture can be drawn over the prequantization bundle $L$ which is considered
as a real rank 2 bundle. Then the riemannian metric in the fibers  is given
just by the real part of the hermitian structure, and it's clear that the prequantization
connection keeps this real part as well as whole the hermitian structure.
More generally, Rothstein theorem (see [17]) distinguishes the case
as the main one (for Batchelor trivial manifolds). Let us consider
the prequantization principle $U(1)$ - bundle as an even supersymplectic manifold  $(P, \Om_s)$. The following statement
gives an elegant reformulation of the Bohr - Sommerfeld condition:
\proclaim{Proposition A} A cycle $S \subset M$ is Bohr - Sommerfeld Lagrangian
if and only if it is the projection of a Lagrangian cycle of  even
supersymplectic manifold $(P, \Om_s)$.
\endproclaim

The proof is rather routine: if $S$ is a Bohr - Sommerfeld Lagrangian
cycle then it can be lifted to $P$ ( we called the lifting planckian cycles) such that
the resulting $\tilde S$ is horizontal at each point. Hence the restriction of
$\Om_s$ to $\tilde S$ is trivial. Conversely, if a submanifold $K \in P$ is Lagrangian with respect to
$\Om_s$ (so it is isotropical and has maximal dimension) then

1) it has to be horizontal at each point; otherwise the restriction of $\Om_s$
to $K$ should be nontrivial,

and

2) it has to give a Lagrangian submanifold being projected to
the base $(M, \om)$; otherwise the restriction should be nontrivial.

Here (to be honest) we mention that the simplification in the definition
could be rewritten in terms of contact geometry: the principle bundle
$P$ with connection $A$ (rescaled to be real) gives an example of
contact manifold with contact form $A$. Indeed,
$dA$ by the definition is projected to $\om$ thus the wedge product
$$
A \wedge (dA)^{\wedge n}
$$
is a volume form on $P$ which is (locally) the product of the Liouville volume
form lifted from $M$ and our nonvanishing vertical 1- form $A$. Then 
\proclaim{Proposition A'} A cycle $S \subset M$ is Bohr - Sommerfeld Lagrangian
iff it is the projection of a maximal isotropical cycle $\tilde S \subset P$
such that
$$
A|_{\tilde S} = dA|_{\tilde S} \equiv 0.
$$
\endproclaim

But we can go further, where the contact geometry doesn't appear: we mean applications
of supersymplectic settings in the nonabelian cases. In the main text we discuss
ALG(a) - programme but what about any ALG(n) - generalization; "ALG(n)" means that
we would like to consider some nonabelian algebraic Lagrangian geometry.
The subject of ALG(n) now exists just as a number of suggestions and examples.
First steps in this way could be done using the superproperties of
the Bohr - Sommerfeld condition. Namely, let $(\Cal M, \Om_s)$
is an even supersymplectic manifold which is projected by some
appropriated Batchelor trivialization to a "classical" symplectic
manifold $(M, \om)$. Let $\Cal L_S$ is the moduli space of usual
Lagrangian cycles on $M$ of a fixed topological type. It's natural to consider
the set of super Lagrangian cycles in $(\Cal M, \Om_s)$ so the cycles
which satisfy
$$
\Om_s|_{\tilde S} \equiv 0.
$$
The local geometry of these cycles is not well understood
(one has some version of the Darboux theorem in the super case, see f.e. [12], but
what about the Darboux - Weinstein theorem?) but it can be done in appropriate cases.
The projections of these super Lagrangian cycles would lie in $\Cal L_S$
and should be analogous to the Bohr - Sommerfeld Lagrangian cycles
in the abelian case. To project $\tilde S$ to $M$ one  uses
some Batchelor trivialization (in general) therefore (in general)
there is some dependence on $H^2(M, \Z)$ (since every two Bachelor
trivialzations differ by a line bundle, see [17]): the homology class of
projected cycles depends on the Batchelor trivialization choice. To illustrate
the generalization way as a sufficiently new method we study the following simple
\subheading{Example} We discuss a classical construction relating
riemannian and holomorphic geometries namely the twistor construction
(see f.e. [19]). We take a symplectic 4 - dimensional manifold
$(M, \om)$ and to perform the construction we choose a riemannian metric
$g$, compatible with the symplectic form. Remark that in contrast with
the material of subsections 1.6, 4.3 we don't impose the integrability condition
on the corresponding almost complex structure. Consider the twistor space, defined
by the conformal class of the riemannian metric. It is a 6 - dimensional real manifold
denoted as $Y$, which is fibered over $X$ with fiber $\C \proj^1$. One could construct it
as follows: let $W^-$ be the spinor bundle defined by our metric, canonically equipped with
a hermitian structure (see, f.e. [26]). The adjoint $SO(3)$ - bundle $ad W^- \to M$
is equipped with: 

1) fiberwise standard riemannian metric

and 

2) a connection, induced by the Levi - Civita connection, compatible
with the riemannian metrics. 

Topologically over a point $x \in M$ the fiber $adW^-_x$ is 3 - sphere, fibered
over the twistor 2 - sphere $\C \proj^1_x$ as in the Hopf bundle picture. 
(It can be seen as the construction of subsection 1.2: we construct
the projective line, factorizing unit sphere $S^3$ in $\C^2 = W^-_x$.)
The total space of the adjoint bundle is the twistor space $Y$. 

Further, our riemannian metric $g$ defines 
over the twistor space a hermitian triple $(G, J, \Om)$, see [26] (notice, that
the triple is defined by the riemannian metric, not by its conformal class).
Here $J$ is usual twistor complex structure, which depends only on
the conformal class, but $G$ and $\Om$ feel the deformations
inside of a fixed class. Suppose that $\Om$ is closed (so our $g$ is a special,
satisfying some appropriate condition) and that the class of $\Om$ is integer.
Thus one can take the pair $(Y, \Om)$ and perform a few first steps
of ALAG - programme, fixing an appropriate topological class
in $H_3(Y, \Z)$. Then we get the moduli space $\Cal B_S^Y$ of Bohr - Sommerfeld 
Lagrangian cycles in $Y$ with respect to some prequantization data. It's easy to see
that every Bohr - Sommerfeld Lagrangian cycle $\tilde S \in \Cal B^Y_S$
admits smooth projection to the based manifold $M$. But it's not longer true
that the projection is a Lagrangian submanifold in $M$ since
our symplectic form $\Om$ is twisted along the fibers as well as the classical
twistor complex structure $J$. But since we started with a symplectic manifold
one has in the picture some additional ingredient: namely, the corresponding
to $\om$ and $J$ almost complex structure $I$ defines a smooth embedding
(see [26]):
$$
i_I: M \to Y,
$$
and we can separate a component from the moduli space $\Cal B^Y_S$ consists of
the Bohr - Sommerfeld cycles intersecting the image of $i_I$
in maximal dimension 2. Then these intersections are transformed
by $i_I^*$ back to $M$ and it's clear that the resulting
submanifolds are Lagrangian with respect to
our given symplectic form $\om$. Remind that we can perform this construction
if  the given symplectic manifold admits a compatible riemannian metric
which possesses the property that the corresponding induced 2- form
$\Om$ is closed. But in any case the total space $ad W^-$ for any riemannian metric
can be considered as an even supersymplectic manifold, endowed with
the corresponding to $\om$ and $g$ supersymplectic even form ($g$ here defines
simultaneously the bundle $W^-$, the hermitian structure and the desired
connection), and hence one could
consider superlagrangian submanifolds of it. The answers are different for these
two constructions: indeed, the first case is twisted while the second is
"constant" along the fibers. The  first construction
belongs to riemannian geometry (f.e. one can simply remove the symplectic structure
and consider the picture just in the riemannian framework;
f.e. such a construction can be (and has been) done for $S^4$ with standard
self dual metric) but the second lives in symplectic geometry.
In [6] it was shown that these "pseudolagrangian" submanifolds
for $S^4$ (which doesn't admit at all any symplectic structure)
just are minimal in the corresponding topological classes. 
(However in [19] one finds an interesting way
to quantize in terms of the twistor spaces.)
But in the second case we get something new. For instance 
it isn't quite clear what is the dimension
of the projected superlagrangian cycles in $M$. For some
more general cases it could be that the dimension is less then
half of the dimension of the based manifold. Thus
one can introduce to the picture considerations
of any rank bundles and some special isotropic submnaifolds,
mixing the algebraic Lagrangian geometry for both the abelian and nonabelian
cases. All together they give us a completion of the quantization picture,
presented in the main text.

$$
$$

\subheading{Special superfunctions}

From the geometrical point of view (see [3]) classical and quantum mechanical 
systems are distinguished by the presence in the last one
of some additional kinematical ingredient which was unknown in
the classical case. The quantum phase space in contrast with the classical one
is endowed by a riemannian metric, compatible with given symplectic structure
such that the corresponding complex structure is integrable. However
in geometric formulation of quantum mechanics, which underlies 
of all our quantization programme, a natural question about the measurement process arises.
Indeed, in the formulation any quantum observable is a special function
(quasisymbol), which has exact numerical values in the points of
the quantum phase space (and this value is precisely the same as the expectation
value of the original operator, see Section 1). Of course, the presence of
the riemannian metric gives us the corresponding measurement process,
but why we should perform it when we have these exact values? From this point
of view it would be natural to look for a realization of the quantization programme
where some geometrical objects, corresponding to the classical observables,
carry more "uncertant" properties than to be just smooth functions.
At the same time we've seen in subsection 4.4 that in quasi classical limit
one gets some strange objects $Y_f$ which are neither pure functions no pure
vector fields. Thus another question is: what is $Y_f$?
 
Coming back to the basic constructions we see that these two questions are related
if we restrict the discussion of subsection 4.4 to
the moduli space of unweighted Bohr - Sommerfeld Lagrangian cycles.
Then for every smooth function $f$ one has over the moduli space
$\Cal B_S$ the object $Y_f$. One could say that $Y_f$ is a super function
in some super symplectic setting. Recall, that as $T^*M$
is a "classical" symplectic manifold for any smooth manifold $M$, one has at the same time an
odd super symplectic manifold $\Pi T^* M$ (see f.e. [12]). The manifold $\Pi T^* M$
is given by the reversing of the parity of the fibers of $T^*M$
(of course, there is a huge list of references for the "super" subject, but I
would like to refer here  survey [12] since my interest to the subject
was intended by its author), and the corresponding odd bracket is called  
Butten bracket. Now what is a superfunction for $\Pi T^*M$? It is a sum
of polyvector fields of different degrees on the given manifold $M$ (f.e.
pure numerical function has degree zero) and the Butten bracket could be reduced
to the standard Schouten bracket for multivector fields. 
Therefore one can regard our object $Y_f$ as a superfunction
on odd supersymplectic manifold $\Pi T^*\Cal B_S$.  This function is rather
special: it has different degrees in different points and it isn't
smooth. But the representation of the classical observables by such type
superfunctions looks very natural and interesting. This representation
is dynamical (we discussed it above). One has a Lie bracket on
the set of such special superfunction hence it would be natural
to suppose that the Butten bracket, defined by smooth functions
only, can be extended to a more general class of superfunctions and
this extension gives us our Lie bracket.
Remark, that while our quaisymbols $F_f$'s have exact numerical values
at the quantum states our special superfunctions
do not have exact values in all the points --- only in 
some "eigenstates", where $f$ is constant.  Thus from some point of view
(expressed at the beginning of this last remark) a way to use
these supersetting is more attractive: here we have some intriguing uncertainty
as one accustomed to have. Here we have an appropriate dynamical property
but we don't have any invariant measurement process. If we find
this desired ingredient (say, in this supersetting) it were
an important movement in the geometric quantization theory.

\enddocument